%% file: ellipsoids.tex
\documentclass{article}

\usepackage{algorithm}
\usepackage{algpseudocode}
\algnewcommand\algorithmicforeach{\textbf{for each}}
\algdef{S}[FOR]{ForEach}[1]{\algorithmicforeach\ #1\ \algorithmicdo}

\usepackage{booktabs}

\usepackage[
  isbn       = false,
  giveninits = true,
  url        = false,
  style      = alphabetic,
]{biblatex}

\usepackage{ amssymb, amsmath, amsthm, amsfonts,mathrsfs, verbatim }
\usepackage{mathtools}
\usepackage{url}
\usepackage{xcolor}
\usepackage{pgfplots}
\usepackage{tikz}
\usetikzlibrary{decorations.fractals, arrows, calc}

\usepackage{mathtools}
\usepackage{float}
\usepackage{xcolor}
\usepackage{comment}
\usepackage{xy}
\xyoption{all}
\usepackage{cancel}
\usepackage{graphicx}
\usepackage{float}
\graphicspath{ {images/} }
\usepackage{fullpage}
\usepackage{ stmaryrd }
\usepackage{enumerate}
\usepackage{hyperref}\usepackage{wrapfig}
\usepackage{booktabs,tabularx}

\usepackage{stmaryrd} 
\SetSymbolFont{stmry}{bold}{U}{stmry}{m}{n} 
\usepackage{bm}

\usepackage{paralist}

\DeclareMathOperator{\PCA}{PCA}
\DeclareMathOperator{\diam}{diam}

\DeclarePairedDelimiter{\set}{\{}{\}}
\DeclarePairedDelimiter{\norm}{\|}{\|}

\DeclarePairedDelimiter{\absval}{|}{|}

\DeclareMathOperator{\distance}{d}
\DeclareMathOperator{\dist}{dist}

\newtheorem{theorem}{Theorem}[section]
\newtheorem{corollary}[theorem]{Corollary}
\newtheorem{proposition}[theorem]{Proposition}
\newtheorem{lemma}[theorem]{Lemma}
\newtheorem{example}[theorem]{Example}
\newtheorem{remark}[theorem]{Remark}
\newtheorem{definition}[theorem]{Definition}

\newtheorem{notation}[theorem]{Notation}

\newcommand{\dgm}{\operatorname{dgm}}

\newcommand{\uu}[1]{\underline{x}}

\newcommand*{\R}{\mathbb{R}}

\newcommand{\NN}{\mathbb{N}}

\newcommand{\RR}{\mathbb{R}}

\newcommand{\id}{\mathrm{id}}
\newcommand{\isom}{{\ \cong\ }}
\newcommand{\define}[1]{{\bf \boldmath{#1}}}
\newcommand{\st}[1]{\mathcal{#1}}
\newcommand{\gmnf}{\st{M}}
\newcommand{\ad}{n}

	\newcommand{\df}[1]{{\bf{#1}}}
 \newcommand{\gp}{x}
 	\newcommand{\dfeq}{\mathrel{\mathop:}=}
\newcommand{\md}{m}
\newcommand*{\eps}{\varepsilon}

\DeclareMathOperator{\crit}{crit}
\renewcommand*{\crit}{\rm {crit}}

\newcommand*{\Ev}{\mathbb{E}}

\newcommand{\zz}[1]{}

\usepackage{cleveref}

\usepackage{subfiles}
\usepackage{caption}
\usepackage{subcaption}

\pgfplotsset{compat=1.18}

\newcommand{\anachange}[1]{\textcolor{red}{#1}}

\begin{document}

\title{Persistent Homology via Ellipsoids}
\author{
 Niklas Canova\thanks{
ETH Zurich,
\texttt{niklas.canova@math.ethz.ch}}, 
 Sara Kali\v{s}nik\thanks{
Pennsylvania State University,
\texttt{skalisnik@psu.edu}},  
Aaron Moser\thanks{
MIT,
\texttt{maaron@mit.edu}},  
Bastian Rieck\thanks{
University of Fribourg,
\texttt{bastian.grossenbacher@unifr.ch}},  
 Ana \v Zegarac\thanks{
ETH Zurich,
\texttt{ana.zegarac@math.ethz.ch}}
}

\date{}
\maketitle

\begin{abstract}
Persistent homology is one of the most popular methods in topological data analysis. An initial step in any analysis with persistent homology involves constructing a nested sequence of simplicial complexes, called a filtration, from a point cloud. There is an abundance of different complexes to choose from, with \v{C}ech, Rips, Alpha, and witness complexes being popular choices. 
In this manuscript, we build a different type of geometrically informed simplicial complex, called a \emph{Rips-type ellipsoid complex}. This complex is based on the idea that ellipsoids aligned with tangent directions (with respect to the data) better approximate the data compared to conventional (Euclidean) balls centered at sample points, as used in the construction of Rips and Alpha complexes, for instance. We use Principal Component Analysis to estimate tangent spaces directly from samples and present an algorithm as well as an implementation for computing \emph{Rips-type ellipsoid barcodes}, i.e., topological descriptors based on Rips-type ellipsoid complexes.
Additionally, we show that the ellipsoid barcodes depend \emph{continuously} on the input data so that small perturbations of a $k$-generic point cloud lead to proportionally small changes in the resulting ellipsoid barcodes. This provides a theoretical guarantee analogous, if somewhat weaker, to the classical stability results for Rips and \v{C}ech filtrations.
We also conduct extensive experiments and compare Rips-type ellipsoid barcodes with standard Rips barcodes. Our findings indicate that Rips-type ellipsoid complexes are particularly effective for estimating the homology of manifolds and spaces with bottlenecks from samples. In particular, the persistence intervals corresponding to ground-truth topological features are longer compared to those obtained using the Rips complex of the data. Finally, we demonstrate that Rips-type ellipsoid barcodes outperform Rips barcodes in classification tasks.
\end{abstract}

\tableofcontents

\section*{Introduction}

\subfile{sections/introduction}

\section{Preliminaries: Filtrations and Persistent Homology}

\subfile{sections/preliminaries}

\section{Rips-Type Ellipsoid Complexes and Their Properties}
\subfile{sections/ellipsoid}

\section{Persistent Homology via Ellipsoids: the Algorithm}

\subfile{sections/algorithm}

\section{Stability of Rips-Type Ellipsoid Complexes}

\subfile{sections/stability}

\subfile{sections/experiments}

\section{Conclusion and Future Directions}
Previous experiments~\cite{Breiding2018} and theoretical results~\cite{Kalisnik2024} support the statement that using shapes elongated along tangent directions one can reduce the size of the sample from a manifold while still capturing its shape. 

In this paper we present code (available at \url{https://github.com/a-zeg/ellipsoids}) for computing persistent homology with such elongated shapes, the so called Rips-type ellipsoid complexes, where simplices are included based on intersections of ellipsoids, not balls. These complexes can be constructed for general point clouds (unlike in \textcite{Breiding2018}) and are stable under small~$\delta$-perturbations and mild assumptions on the point cloud, meaning that the associated persistence barcodes vary continuously with the data. We also present the results of extensive experiments where we compare Rips-type ellipsoid barcodes with Rips barcodes. In particular, we show that:
\begin{itemize}
\item 
Working with ellipsoids is particularly suitable when the underlying spaces has bottlenecks (as demonstrated in~Subsection~\ref{subsn:dogbone}) or is a manifold.
\item 
Since ellipsoids better approximate the underlying manifold structure of data than balls, their barcodes lead to better classification results in sparsely sampled point clouds and, in general, allow the user to work with smaller samples. 
\item 
By performing experiments similar to those from~\textcite{Turkes22a} we show that Rips-type ellipsoid barcodes outperform alpha barcodes for classification purposes (see Subsection~\ref{subsec:point-cloud-classification}) in all categories except one.
\end{itemize}

These points demonstrate the strengths of working with Rips-type ellipsoid complexes. The slower computational time is partly offset by the much smaller sample size needed to still capture homology groups compared to the Rips complex.

Several questions remain open. For instance, one could optimize the implementation for more efficient computation of Rips-type ellipsoid complexes, or incorporate ellipsoids into alpha complexes to reduce filtration size. Another possible direction is to replace PCA in the tangent-space estimation step to potentially achieve Gromov–Hausdorff stability of the filtration. In the implementation covered in this paper, the user needs to specify the neighborhood size used in calculating ellipsoid orientations. Removing this requirement is a part of the ongoing work.

Additionally, our experiments focused on data that is well-approximated by an embedded or immersed low-dimensional manifold. Adapting ellipsoids to work well in settings, where the underlying spaces is not necessarily a manifold (but say, a manifold with boundary) is an important direction for future work. There are indications that tools from diffusion geometry, such as those developed in \cite{diffusion_axes_selection}, may help define suitable local orientations even when no manifold model is assumed.

\section{Acknowledgements}

We thank Henry Adams for providing us with the pentagons dataset and Clayton Shonkwiler for creating it.
A.Z.\ would like to thank Marco G\"ahler and Jan Sch\"ussler for their help with the programming part of this project. S.K. would like to thank Shahar Kowalsky for a good discussion during her trip to UNC. We also thank Sergio Serrano de Haro Iváñez for suggesting a speed-up of the initial code.

\section{Statement of AI use}
Claude was used to correct the proof of Lemma 4.19, which in the previous version did not account for the displacement of the ellipsoid centers under perturbation.  The authors have checked the corrected proof and take full responsibility for it. All other mathematics in this paper was written by the authors.

\anachange{Sergio Serrano de Haro Iváñez optimized the code for ellipsoid intersections computations with the help of AI.}
\anachange{We used AI to generate additional point clouds with bottlenecks for the experiments described in Section~\ref{subsec:point-cloud-classification}.}

\newpage

\section{Appendix}

\subfile{sections/appendix}

\printbibliography

\end{document}

%% file: sections/introduction.tex
Methods from computational topology have received increased attention due to their ability to capture characteristic properties of data at multiple scales, while being less reliant on the underlying metric or coordinates~\cite{topodata}. Of these, persistent homology is the most prominent~\cite{CdSO14, CDGO16}.
Given an unstructured dataset in the form of a point cloud, the first step of any analysis based on persistent homology involves building a \emph{simplicial complex} on the data. To approximate the underlying shape of the dataset, a common strategy is to calculate the \v{C}ech, Rips, alpha or witness complex on the dataset~\cite{Dantchev12a, Zomorodian10a}.
While this is a useful strategy in general, many real-world high-dimensional datasets  
actually cluster along low-dimensional manifolds. This statement is known as the manifold hypothesis and it forms a cornerstone of modern data science and machine learning~\cite{manifoldhyp}. 

With this in mind, we build a different type of geometrically informed simplicial complex, a \emph{Rips-type ellipsoid complex}, that is specifically geared  to handle samples from manifolds. In particular, we rely on the insight that \emph{ellipsoids} elongated in tangent directions better approximate the dataset than balls centered at sample points.
This statement is also supported by previous work. For instance, experiments carried out in~\textcite{Breiding2018} demonstrate that given a sample from a variety, complexes distorted in tangent directions combined with persistent homology result in a stronger ``topological signal.'' The drawback of this method is that it requires polynomials that determine the variety to approximate the tangent space. 
In this paper, we do away with this restriction and define Rips-type ellipsoid complexes for general point clouds. Tangent spaces are estimated with the help of PCA directly from the sample~\cite{Jolliffe2016PrincipalCA}. Additionally, we establish a \emph{stability result} for these complexes. In particular, we show that for $\delta$-perturbations~(Definition~\ref{def:delta}) of a $k$-generic point cloud~(Definition~\ref{def:genericity}), the neighborhood structure, the singular values of the local PCA matrices, and the orientations of their singular vectors vary continuously. As a consequence, the associated Rips-type ellipsoid filtrations are interleaved, implying the stability of the corresponding barcodes. This extends the classical stability framework in persistent homology~\cite{CdSO14, CDGO16} to the case of Rips-type ellipsoid complexes. 
The restriction to $k$-generic point clouds and to $\delta$-perturbations (rather than arbitrary perturbations in the Gromov–Hausdorff sense) is needed to rule out pathological geometric configurations, as discussed in more detail in Subsection~\ref{sub:obstruction}.

Although the condition on point clouds may appear restrictive, the set of $k$-generic point clouds with $N$ points in $\mathbb{R}^d$ is open and dense in $(\mathbb{R}^d)^N$ endowed with the product Euclidean topology; in particular, any point cloud can be perturbed by an arbitrarily small amount to a $k$-generic point cloud~(Proposition~\ref{prop:open_and_dense}).

We also provide algorithms and code~(available at \url{https://github.com/a-zeg/ellipsoids}) to compute ellipsoid barcodes and carry out extensive experiments comparing them to Rips complexes. We demonstrate that:
\begin{itemize}
\item 
Working with ellipsoids is particularly suitable when the underlying space is a manifold or has bottlenecks (see Subsection~\ref{subsn:dogbone}). \emph{Persistence barcodes} arising from ellipsoid complexes exhibit a larger signal-to-noise ratio; more specifically, the persistence intervals corresponding to a ground-truth topological feature are longer~(as compared to the intervals obtained when using the Rips complex of the data).
\item 
%Working with ellipsoids requires smaller samples while still allowing one to acquire meaningful information about the shape.
Ellipsoid barcodes lead to better classification results in sparsely sampled point clouds and, in general, allow the user to work with smaller samples, which is consistent with the theoretical results from~\cite{Kalisnik2024}.
\item 
Using datasets similar to those introduced in~\textcite{Turkes22a}, we show that ellipsoid barcodes outperform Rips barcodes and also outperform alpha barcodes generated using Distance-to-Measure as the filtration
function~(see Subsection~\ref{subsec:point-cloud-classification}) in classification tasks in all categories except one.
\end{itemize}

We note that the computation time of ellipsoid complexes is comparable to that of the traditional Vietoris-Rips complex, as discussed in Section~\ref{subsec:computation_time}.

Our paper is organized as follows.
Section~\ref{subsec:prelims} recalls the necessary background on simplicial complexes built on point clouds, filtrations, and persistent homology.
In Section~\ref{sec:ellipsoid} we introduce \emph{Rips-type ellipsoid complexes} in both smooth and discrete settings and relate them to classical Rips complexes.
Section~\ref{sec:alg} describes the algorithmic construction based on local PCA and explains how ellipsoid intersections are computed.
In Section~\ref{sec:stability} we prove a stability statement for Rips-type ellipsoid complexes.
Finally, Section~\ref{sec:Experiments} presents various numerical experiments with ellipsoid barcodes.

%% file: sections/preliminaries.tex
\label{subsec:prelims}

In this section we review the definitions of simplicial complexes and filtrations of point clouds, explain how one constructs persistence modules based on point clouds and briefly explain how persistent homology works and what information about the underlying point cloud it provides.

Persistent homology is an adaptation of homology~\cite{hatcher} to the setting of point clouds, i.e., finite metric spaces that arise from applications. The concept appeared independently in the works of Barannikov~\cite{barannikov1994framed}, Frosini~\cite{Frosini}, Robins~\cite{Robins}, and Edelsbrunner et al.~\cite{elz-tps-02}. For an in-depth introduction to persistent homology, see~\textcite{topodata, carlsson2014}.
The goal of persistent homology is to provide a bridge between discrete and non-discrete topological spaces: point clouds, being discrete topological spaces,  have no non-trivial topological features. To obtain topological features, one needs to turn the point cloud into a topological space. One way to accomplish this is to assign for every parameter $\varepsilon >0$ a topological space, more specifically, a simplicial complex, to the point cloud and track the \emph{evolution} of the topological features as the parameter~$\varepsilon$ varies. 
\begin{definition}
	An \define{abstract simplicial complex} $(\Sigma,V)$ is given by a set  $V$  whose elements we call \define{vertices} and a set   $\Sigma$ of non-empty finite subsets of $V$. This data satisfies the following properties: we have that (1) $\{v\}\in \Sigma$ for all $v\in V$, and (2) if $\sigma\in \Sigma$ and $\tau\subset\sigma$, then $\tau\in\Sigma$. If $\sigma\in \Sigma$ has cardinality $p+1$, we say that $\sigma$ is a \define{$p$-simplex}, or a simplex of \define{dimension $p$}. A \define{simplex} is a $p$-simplex for some $p\in \mathbb{N}$.
\end{definition}%
One common way of assigning a simplicial complex to a point cloud is to take the \v{C}ech complex. In this paper all metric spaces are subspaces of Euclidean space $\RR^d$.
\begin{definition} Let $(X, d)$ be a metric space.
The \define{\v{C}ech complex of $X$, $\text{\v{C}}_\varepsilon(X)$, at scale $\varepsilon$} is the abstract simplicial complex with the vertex set $X$, where $\sigma=[v_0, v_1, \ldots, v_n]$ is an $n$-simplex in $\text{\v{C}}_\varepsilon(X)$ if and only if $B_\varepsilon(v_0)\cap \ldots \cap B_\varepsilon(v_n) \neq \emptyset$.   In other words, $\text{\v{C}}_\varepsilon(X)$ is an abstract simplicial complex with the vertex set $X$, where $v_0, v_1, \ldots, v_n$ form an $n$-simplex precisely when the balls of radius $\varepsilon$ centered at these points have a non-empty intersection.
\end{definition} 

The \v{C}ech complex $\text{\v{C}}_\varepsilon(X)$ at scale $\varepsilon$ has the same homotopy type as the union of balls grown around the data points with radius $\varepsilon$. This follows directly from a result referred to as the nerve lemma \cite{nerve}. 
In applications the so-called Vietoris--Rips~(or just Rips) complex is more popular because it is easier to store. 
\begin{definition}
Given a metric space $(X, d)$ and a real number $\varepsilon \geq 0$, the \define{Rips complex of $X$ at scale $\varepsilon$} is:
 \[
R_\varepsilon (X) = \{\sigma \subseteq X \,|\, d(x,y)\leq 2\varepsilon,\forall x,y \in \sigma\}.
 \]
 When $X$ is clear from the context, we write just $R_\varepsilon$ instead of $R_\varepsilon (X)$.
\end{definition} 

\begin{figure}[ht]
    \centering
    \subcaptionbox{$\varepsilon = 0$}{%
        \includegraphics[width=0.25\linewidth]{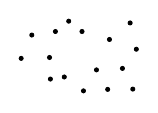}%
    }%
    \subcaptionbox{$\varepsilon = 1.5$}{%
        \includegraphics[width=0.25\linewidth]{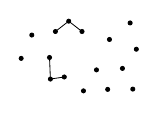}%
    }%
    \subcaptionbox{$\varepsilon = 1.75$}{%
        \includegraphics[width=0.25\linewidth]{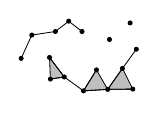}%
    }%
    \subcaptionbox{$\varepsilon = 2.0$}{%
        \includegraphics[width=0.25\linewidth]{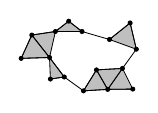}%
    }
    \caption{%
        Four stages of a Rips complex construction for a point cloud, showing simplices up to dimension~$2$.
    }
    \label{fig:Rips complex example}
\end{figure}

For each $\varepsilon \leq \varepsilon'$ we have an inclusion $\text{\v{C}}_\varepsilon(X) \hookrightarrow  \text{\v{C}}_{\varepsilon'}(X)$ as well as an inclusion $R_\varepsilon (X) \hookrightarrow R_{\varepsilon'} (X)$. \Cref{fig:Rips complex example} illustrates the Rips complex construction for the point cloud depicted in \Cref{fig:Rips complex example}(a).
Taking a family of Rips or \v{C}ech complexes indexed over real $\varepsilon \geq 0$ yields a filtered simplicial complex.
 \begin{definition}
A \define{filtered simplicial complex} is a collection $K=\{K_\varepsilon\}_{\varepsilon\in \mathbb{R}_{\geq 0}}$ of simplicial complexes indexed by non-negative real numbers with the property that $K_\varepsilon\subset K_{\varepsilon'}$ whenever $\varepsilon\leq \varepsilon'$.
\end{definition} 
Applying the homology functor $H_k$ in degree $k$  to a filtered simplicial complex, we obtain what is called a `persistence module'~\cite{carlsson2014}.
\begin{definition}\label{definition persistence module}
	A \define{persistence module} $\mathbf{V}$ is a collection of indexed vector spaces $\{ V_t|\,t\in\RR \}$ and linear maps $\{ v_a^b|\,v_a^b\colon V_a\to V_b,\,a\leq b \}$ such that the composition has the properties $v_b^c\circ v_a^b=v_a^c$ whenever $a\leq b\leq c$ and $v^a_b$ is the identity map whenever $a=b$. 
\end{definition} 

To quantify how close two persistence modules are to being isomorphic, we use the notion of a $\delta$-interleaving~\cite{CCsGG2009}; in particular, a $0$-interleaving is an isomorphism.

\begin{definition}[$\delta$-interleaving]
Let $\mathbf{U}$ and $\mathbf{V}$ be two persistence modules, and $\delta \geq 0$. A family of linear maps $(f_t : U_t \to V_{t+\delta})_{t\in \R}$ such that $v_{t+\delta}^{s+\delta}\circ f_t = f_s\circ u_t^s$ is called a $\delta$-\define{morphism} from $\mathbf{U}$ to $\mathbf{V}$. The space $\mathrm{Hom}^\delta(\mathbf{U}, \mathbf{V})$ denotes the set of $\delta$-morphisms from $\mathbf{U}$ to $\mathbf{V}$. We write
\[
    1_{\mathbf{U}}^\delta = (u_t^{t+\delta})_{t\in \R}.
\]
Note that the image of $1_{\mathbf{U}}^\delta$ is just the persistence module $\mathbf{U}$ shifted by $\delta$.

Two persistence modules $\mathbf{U}$ and $\mathbf{V}$ are said to be $\delta$\define{-interleaved} if there exist maps $\Phi \in \mathrm{Hom}^\delta(\mathbf{U}, \mathbf{V})$, ${\Psi \in \mathrm{Hom}^\delta(\mathbf{V}, \mathbf{U})}$ such that 
\[
\Psi \Phi = 1_{\mathbf{U}}^{2\delta} \quad \text{and} \quad \Phi \Psi = 1_{\mathbf{V}}^{2\delta}.
\]
\end{definition}
\begin{remark}\label{rem:interleaving}
The notion of $\delta$-interleaving can be extended from persistence modules to geometric filtrations, such as Vietoris--Rips or Čech complexes. 
In that case, a $\delta$-morphism is defined as a family of simplicial maps 
\[
f_t\colon K_t \to L_{t+\delta}
\]
satisfying the compatibility condition $j_{t+\delta}^{s+\delta}\circ f_t=f_s\circ i_t^s$, where $i_t^s\colon K_t\hookrightarrow K_s$ and $j_t^s\colon L_t\hookrightarrow L_s$ denote the inclusions.

Two filtrations $(K_t)_{t\in\R}$ and $(L_t)_{t\in\R}$ are then $\delta$-interleaved if there exist such families $(f_t)$ and $(g_t)$ with 
$g_{t+\delta}\circ f_t = i_t^{t+2\delta}$ and 
$f_{t+\delta}\circ g_t = j_t^{t+2\delta}$, 
where $i_t^{t+2\delta}$ and $j_t^{t+2\delta}$ are the natural inclusions in the respective filtrations.
\end{remark}

The basic building blocks in the theory of persistence modules are interval modules.
\begin{definition}\label{definition interval module}
For an interval $[b, d)$ we denote by $\mathbb{I}_{[b, d)}$ the persistence module 
\[
 (\mathbb{I}_{[b, d)})_t=\begin{cases}\Bbbk & \textrm{ for } t \in [b, d)\\ 0 & \textrm{ otherwise}\end{cases} \qquad \textrm{ and }\qquad  i_s^t = \begin{cases}\id_{\Bbbk} & \textrm{ for } s \leq t, \textrm{ and }s, t \in [b, d)\\ 0 & \textrm{ otherwise}\end{cases}.
\]
The \emph{lifespan} of $\mathbb{I}_{[b, d)}$ is $d-b$.
\end{definition}
The celebrated \define{decomposition theorem} guarantees that persistence modules that arise from Rips and \v{C}ech complexes and similar filtrations built on point clouds can be expressed as direct sums of `interval modules'~\cite{carlsson2014}.

\begin{theorem}\label{thm:decomposition}
Let $(X, d)$ be a finite metric space and $\{R_\varepsilon\}_{\varepsilon\in \mathbb{R}_{\geq 0}}$ the Rips filtration associated to $X$. Then the persistence module $\{H_k(R_\varepsilon)\}_{\varepsilon\in \mathbb{R}_{\geq 0}}$ over $\Bbbk$ can be decomposed as
 \[
\{H_k(R_\varepsilon)\}_{\varepsilon\in \mathbb{R}_{\geq 0}} \isom \bigoplus_{l\in L} \mathbb{I}_{[b_l, d_l)}
\]
for some $b_l \in [0,+\infty), d_l \in [0,+\infty]$, with $b_l<d_l$ for all $l$.
%The $k$-dimensional \define{barcode} associated to $X$ is
%\[
%\{[b_1, d_1), [b_2, d_2), \ldots, [b_l, d_l)\}.
%\]
\end{theorem}

Instead of Rips or \v{C}ech complexes we can also construct Rips-type ellipsoid complexes (which we define in Section~\ref{subsec:ellipsoid}) to obtain barcodes. By Theorem~\ref{thm:decomposition}, or more generally by~\cite{CrawleyBoevey15}, we can associate to a finite point cloud $(X,d)$ a collection of intervals.
 We can represent this collection as a \define{barcode} or alternatively, as a \define{persistence diagram}. We will use both interchangeably.

\begin{definition}\label{def:persistence_diagram}
Let $\mathbf{V}$ be such that $\mathbf{V} \isom \bigoplus_{l\in L} \mathbb{I}_{[b_l, d_l)}$. The \define{barcode} is the plot obtained by drawing, for each \( \mathbb{I}_{[b_l, d_l)} \), a horizontal line from \(t=b_l \) to \(t=d_l \) (or a ray if \(d_l=+\infty\)). 
The \textbf{persistence diagram} of $\mathbf{V}$, denoted by \(\dgm(\mathbf{V})\),
is the multiset in the extended plane $\overline{\mathbb{R}}^2$, where $\overline{\mathbb{R}}:=\mathbb{R}\cup\{\pm\infty\}$, consisting of the points \(\left\{(b_{l}, d_{l})\right\}_{l\in L}\subset\overline{\RR}^{2}\) (counted with multiplicity) and the diagonal \(\Delta:=\left\{(x,x)\mid x\in\overline{\RR}\right\}\) (where each point on \(\Delta\) has infinite multiplicity).
\end{definition}

\begin{figure}[ht]
    \centering
    \subcaptionbox{Persistence barcode}{%
        \includegraphics[width=0.45\linewidth]{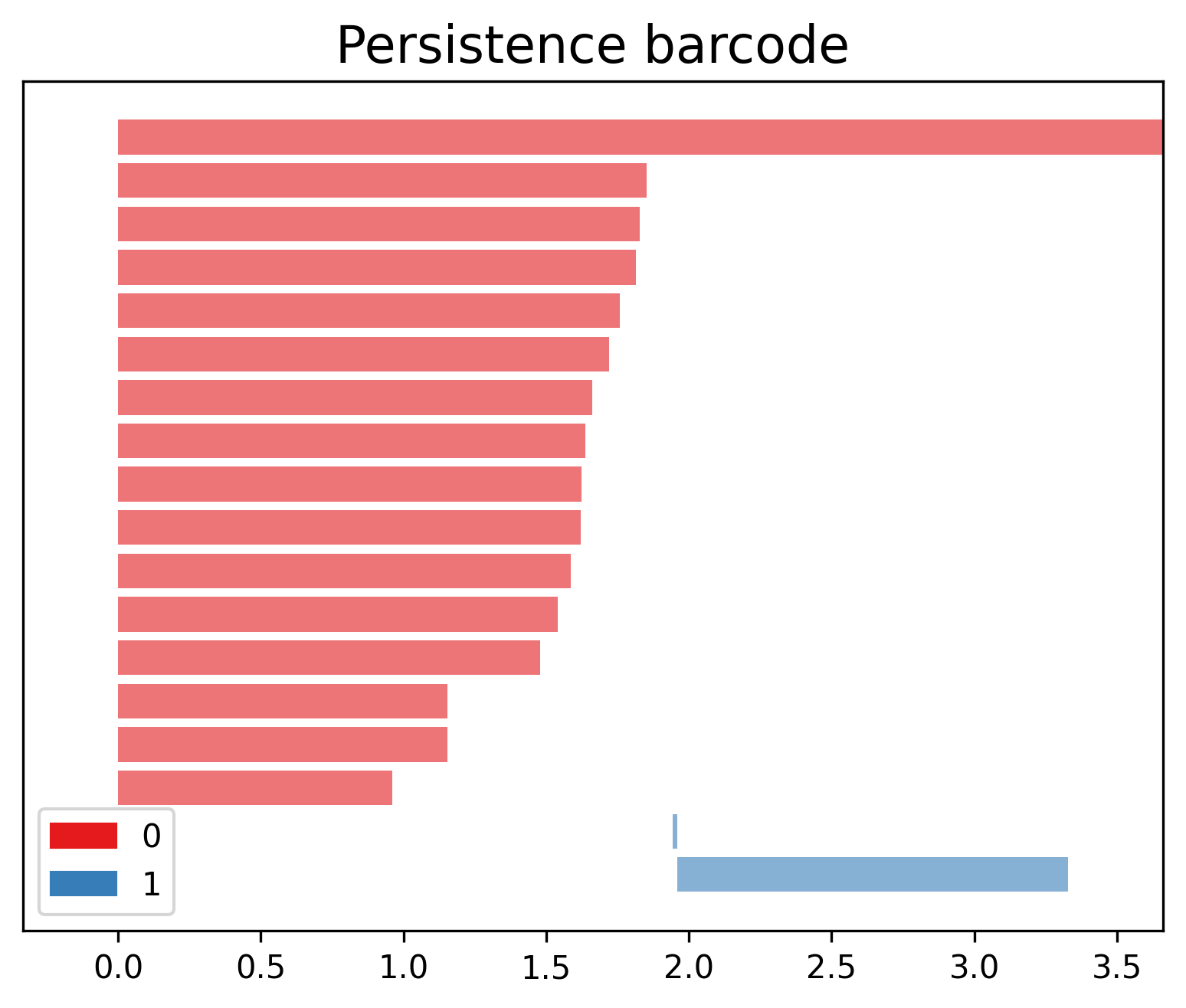}%
    }%
    \hfill
    \subcaptionbox{Persistence diagram}{%
        \includegraphics[width=0.49\linewidth]{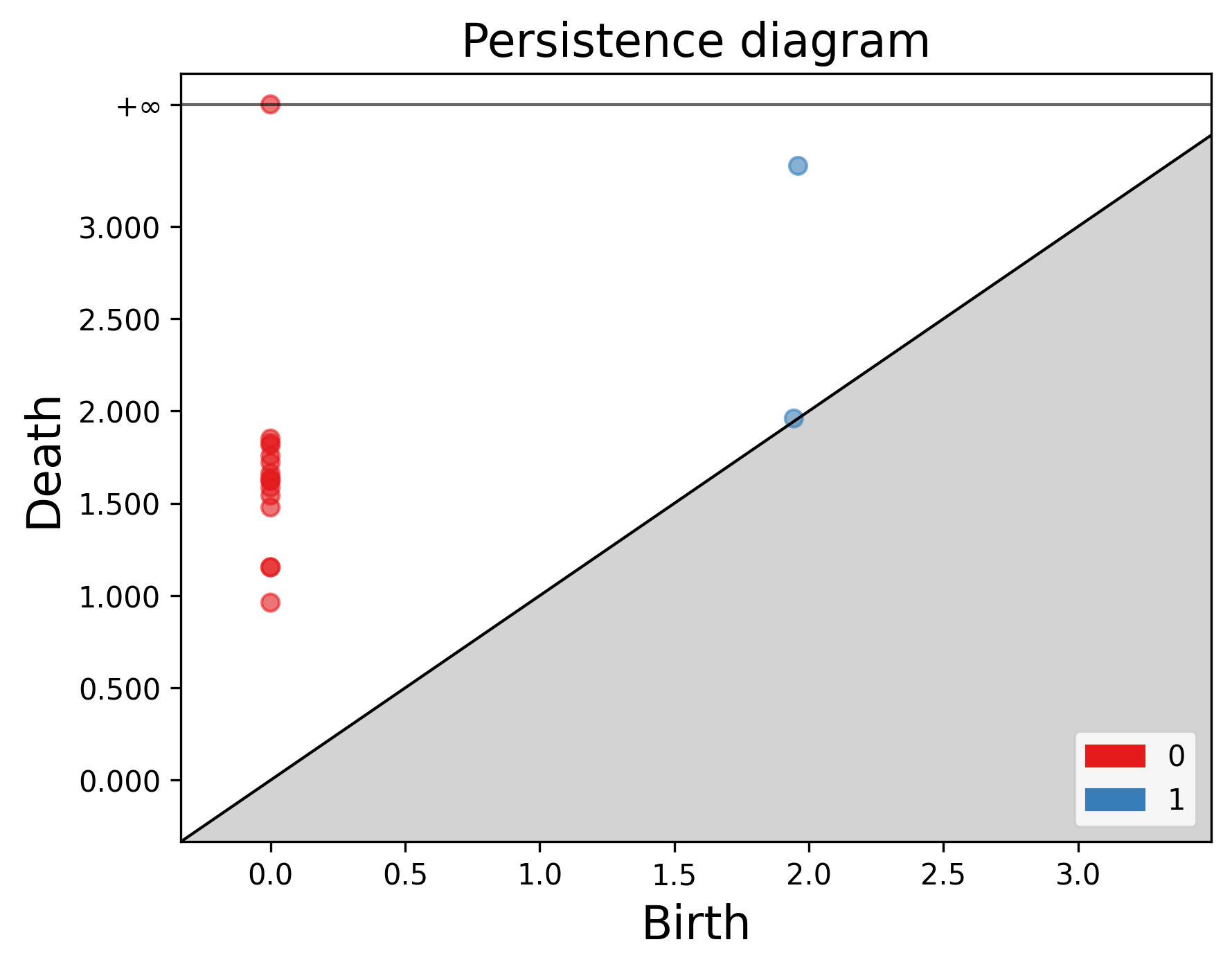}%
    }%
    \caption{%
        Example of a persistence barcode and a persistence diagram. The underlying point cloud dataset is shown in Figure~\ref{fig:Rips complex example}.
    }
    \label{fig:Rips barcode example}
\end{figure}

Each interval in a barcode (equivalently, each point in a persistence diagram) obtained from a point cloud via a simplicial complex filtration corresponds to a topological feature that appears at the parameter value given by the interval’s left endpoint and disappears at the value given by its right endpoint. See Figure~\ref{fig:Rips barcode example} for an illustration.

%
% \begin{figure}[ht]
%     \centering
%     \includegraphics[width=0.30\linewidth]{figures/example_barcode}
%     \caption{%
%         The Rips barcode for the point cloud depicted in \protect\Cref{fig:Rips complex example}(a).
%     }
%     \label{fig:Rips barcode example}
% \end{figure}
%

To measure the similarity between barcodes, resp. persistence diagrams, we use the bottleneck distance.
\begin{definition}[Bottleneck distance]\label{def:bottleneck}
Let $\mathbf{U}, \mathbf{V}$ be pointwise finite-dimensional persistence modules and let 
\[
\Pi = \{ \pi \colon \dgm(\mathbf{U}) \to \dgm(\mathbf{V}) \mid \pi \text{ is bijective} \}
\]
be the set of all bijections between $\dgm(\mathbf{U})$ and $\dgm(\mathbf{V})$. Then, the \define{bottleneck distance} between $\dgm(\mathbf{U})$ and $\dgm(\mathbf{V})$ is
\[
d_b(\dgm(\mathbf{U}), \dgm(\mathbf{V})) := \inf_{\pi \in \Pi} \; \sup_{x \in \dgm(\mathbf{U})} \; \| x - \pi(x) \|_\infty.
\]
\end{definition}

Stability is crucial in applications, as it ensures that persistence diagrams vary continuously with the data~\cite{Cohen-Steiner2007, CCsGG2009, CDGO16}.
The following theorem~\cite{CDGO16} makes this precise: if two persistence modules are $\delta$-interleaved, then their persistence diagrams differ by at most $\delta$ in bottleneck distance.

\begin{theorem}\label{thm:Chazal}
If $\mathbf{U}, \mathbf{V}$ are pointwise finite-dimensional persistence modules that are $\delta$-interleaved, then 
%there exists a $\delta$-matching between the multisets $\dgm(U)$ and $\dgm(V)$.  
%Thus,
the bottleneck distance between the diagrams satisfies the bound
\[
d_b\bigl(\dgm(\mathbf{U}), \dgm(\mathbf{V})\bigr) \leq \delta.
\]
\end{theorem}

%% file: sections/ellipsoid.tex
\label{sec:ellipsoid}
Ellipsoid complexes were first used together with persistent homology by Breiding et al.~\cite{Breiding2018} in the setting of algebraic varieties. That work introduced the \emph{ellipsoid-driven complexes} for point clouds sampled from a variety, with tangent information obtained from the variety’s defining polynomials.

In this section we define \emph{Rips-type ellipsoid complexes} that can be constructed for \emph{any} finite subset of Euclidean space, thereby generalizing \cite{Breiding2018}, and we investigate their properties.

\subsection{Defining Rips-Type Ellipsoid Complexes}

\label{subsec:ellipsoid}

%\anatodo{maybe add explanations about different kinds of ellipsoid complexes: Rips-type, alpha-type, spherisized, PCA-axes?}

{\bf The Topological Setting.} Inspired by~\textcite{Kalisnik2024}, we first provide definitions for Rips-type ellipsoid complexes for the ideal setting, where we have a finite sample $X$ from a known $\mathcal{C}^1$-submanifold \label{subsec:persistence}$\gmnf$ of~$\RR^\ad$. A \df{tangent-normal coordinate system} at $\gp \in \gmnf$ is an $\ad$-dimensional orthonormal coordinate system with the origin in~$\gp$, the first $\md$ coordinate axes tangent to~$\gmnf$ at~$\gp$ and the last $\ad-\md$ axes normal to~$\gmnf$ at~$\gp$.

% previous notation:
	% \begin{definition} \label{def:ellipsoid}
	% 		Let $\gmnf$ be a $\mathcal{C}^1$-submanifold of~$\RR^\ad$ and $\varepsilon \in \RR_{> 0}$. The \df{tangent-normal} \df{$(a, b)$ ellipsoid at scale $\varepsilon$ at point $\gp \in \gmnf$} is the closed ellipsoid in~$\RR^\ad$ with the center in~$\gp$, the tangent semi-axes of length $a$ and the normal semi-axes of length~$b$. Explicitly, in a tangent-normal coordinate system at~$\gp$ the tangent-normal closed ellipsoids are given by
	% 		\[
 %   \begin{array}{rcl}
	% 			E^{a,b}_\varepsilon(x)& \dfeq& \left\{(x_1, \ldots, x_\ad) \in \RR^\ad\,|\,\frac{x_1^2 + \ldots + x_\md^2}{a^2} + \frac{x_{\md + 1}^2 + \ldots + x_\ad^2}{b^2} \leq \varepsilon^2 \right\},
	% 		\end{array}\]
	% 		where $\md$ denotes the dimension of~$\gmnf$ at~$\gp$.
	% 		Observe that the definitions of ellipsoids depend only on the submanifold itself.
	% 	\end{definition}

 \begin{definition} \label{def:ellipsoid}
    Let $\gmnf$ be a $\mathcal{C}^1$-submanifold of~$\RR^\ad$, $q\geq 1$ and $\varepsilon \in \RR_{> 0}$. The \df{tangent-normal} \df{$q$-ellipsoid at scale $\varepsilon$ at point $\gp \in \gmnf$} is the closed ellipsoid in~$\RR^\ad$ with the center in~$\gp$, the tangent semi-axes of length $\varepsilon$ and the normal semi-axes of length~$b \coloneqq \varepsilon / q$. Explicitly, in a tangent-normal coordinate system at~$\gp$ the tangent-normal closed ellipsoids are given by
    \begin{equation} \label{eq:ellipsoid_defn}
   \begin{array}{rcl}
        E^q_\varepsilon(x)& \dfeq& \left\{(x_1, \ldots, x_\ad) \in \RR^\ad\,|\,\frac{x_1^2 + \ldots + x_\md^2}{\varepsilon^2} + \frac{x_{\md + 1}^2 + \ldots + x_\ad^2}{b^2} \leq 1 \right\},
   \end{array}
   \end{equation}
    where $\md$ denotes the dimension of~$\gmnf$ at~$\gp$.
    Observe that the definitions of ellipsoids depend only on the submanifold itself. 
\end{definition}

  If we know $\gmnf$ and have access to its tangent space, then for each point $x$ from the sample $X$ and $\varepsilon>0$ we have an ellipsoid. One way to produce a simplicial complex is to construct a Rips-like complex, in which edges are determined by intersections of ellipsoids. One could also use a \v{C}ech-like construction, including a simplex precisely when the corresponding ellipsoids intersect~(following \textcite{Kalisnik2024}) to keep the theoretical guarantees from the nerve lemma~\cite{nerve}, however, for computational purposes doing that is too expensive.

% old notation
%   \begin{definition} \define{(Ellipsoid Complex)}
%   Let $\gmnf$ be a $\mathcal{C}^1$-submanifold of~$\RR^\ad$ and
% let $(X, d)$ be a finite metric subspace of Euclidean space $\mathbb{R}^d$, where $X \subset \gmnf$. For $x\in X$ let
% $E^{a,b}_\varepsilon(x)$ be the ellipsoid from Definition~\ref{def:ellipsoid}. The \define{(a,b) ellipsoid complex of $X$ at scale $\varepsilon$} is 
% \[
% E_\varepsilon(X)= \{\sigma \subseteq X \,|\, E^{a,b}_\varepsilon(x)\cap E^{a,b}_\varepsilon(y) \neq \emptyset,\forall x,y \in \sigma\}.
% \]
% With other words, $E_\varepsilon(X)$ is an abstract simplicial complex with the vertex set $X$, where $x$ and $y$ are connected by an edge precisely when $E^{a,b}_\varepsilon(x)\cap E^{a,b}_\varepsilon(y) \neq \emptyset$. A higher dimensional simplex is included if and only if all of its edges are in $E_\varepsilon(X)$.
%   \end{definition}

\begin{definition} \label{def:ellipsoid_complex} \define{(Rips-Type Ellipsoid Complex and Filtration)}
    Let $\gmnf$ be a $\mathcal{C}^1$-submanifold of~$\RR^\ad$ and
    let $(X, d)$ be a finite metric subspace of Euclidean space $\mathbb{R}^\ad$, where $X \subset \gmnf$. For $x\in X$ let
    $E^q_\varepsilon(x)$ be the ellipsoid from Definition~\ref{def:ellipsoid}. The \define{Rips-type $q$-ellipsoid complex of $X$ at scale $\varepsilon$} is 
    \[
    E^q_\varepsilon(X)= \{\sigma \subseteq X \,|\, E^q_\varepsilon(x)\cap E^q_\varepsilon(y) \neq \emptyset,\forall x,y \in \sigma\}.
    \]
    In other words, $E^q_\varepsilon(X)$ is an abstract simplicial complex with the vertex set $X$, where $x$ and $y$ are connected by an edge precisely when $E^q_\varepsilon(x)\cap E^q_\varepsilon(y) \neq \emptyset$. A higher-dimensional simplex is included if and only if all of its edges are in $E^q_\varepsilon(X)$.
    Thus, the Rips-type ellipsoid complex is a \emph{flag complex}, i.e., it is fully determined by its edges. 
\end{definition}

\begin{example}
Consider a sample $X$ from a circle depicted in the leftmost image in Figure~\ref{fig:ellipsoids}. The remaining images show the $2$-ellipsoids as well as the Rips-type $2$-ellipsoid complexes built on $X$ at various scales. 
    \begin{figure}[h]
    \centering
    \includegraphics[width=\linewidth]{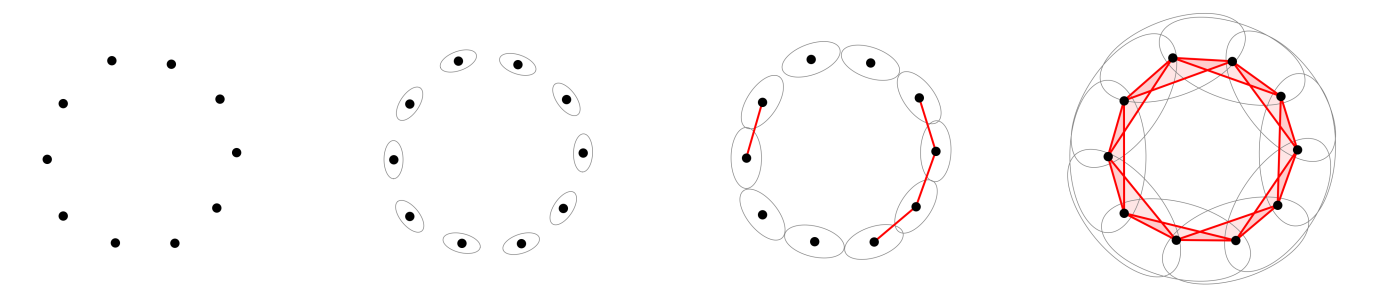}
    \caption{%
        Four stages of a Rips-type ellipsoid complex for a point cloud, showing simplices up to dimension~2.
    }
    \label{fig:ellipsoids}
\end{figure}
\end{example}

{\noindent\bf The Discrete Setting}. The main idea in passing from the topological setting to the discrete setting is that we no longer have access to the underlying manifold and its tangent spaces, but that we have to estimate them directly from the sample.
\label{sec:fitEllipsoid}
To choose the orientation of each of the ellipsoids, Algorithm~\ref{alg:construct_ellipsoids} is used.
\begin{algorithm}[h]
\begin{algorithmic}[1]
\State $L \gets [\,]$\Comment{Initialize empty ellipsoid list}
\ForEach{point $p$ in the point cloud}
    \State Fix user-selected number $k$ of neighbors
    \State Run principal component analysis on the $k$-nearest neighbors of $p$.
    \State Create ellipsoid $E$ centered at $p$ and aligned with PCA eigenvectors.
    \State Set axes lengths of $E$ according to:
\Statex \qquad user-defined axes ratios (ordered by PCA eigenvalues), and
\Statex \qquad the current filtration level.
\State Append $E$ to $L$.
\EndFor
\end{algorithmic}
\caption{Ellipsoid construction algorithm}
\label{alg:construct_ellipsoids}
\end{algorithm}

This algorithm can be implemented efficiently using spatial data structures such as \emph{\mbox{$k$-d} trees}. Building such a data structure for~$N$ points in~$n$ dimensions has a worst-case complexity of $\mathcal{O}(N \log^2 N)$. Calculating a proper ellipsoid for each point then incurs a cost of $\mathcal{O}(k \log N)$ for finding the~$k$ nearest neighbors, followed by $\mathcal{O}(n^2k)$ for calculating principal components \cite[\S 8.6.3, Figure~8.6.1]{GolubVanLoan2013}, with the final alignment step taking constant time.
The total runtime of this algorithm is thus   $\mathcal{O}(Nk\log N + Nn^2k)$.
In lower dimensions and for sufficiently small values of~$k$, this runtime is dominated by finding the~$k$ nearest neighbors, and we  assume that the local PCA calculations effectively run in constant time.

\paragraph{Choice of Tangent and Normal Directions.}
Throughout the paper we assume that the ambient space is $\mathbb{R}^n$.
In our default implementation we take the first $(n-1)$ singular vectors
(corresponding to the largest singular values obtained from the local PCA)
as \emph{tangent} directions, and the remaining singular vector as the
\emph{normal} direction when constructing ellipsoids.  

The number of tangent directions can be selected differently by the user: given any $m \in \{1,\ldots,n-1\}$,
the $m$ largest singular values determine the tangent axes, with the remaining $(n-m)$
axes treated as normal. This allows the construction to adapt to different intrinsic
dimension assumptions.

\begin{notation}\label{def:rips-type}
Given a finite metric subspace \(X \subset \RR^n\), we denote by \(E^q_\varepsilon(X)\) 
the Rips-type ellipsoid complex computed via PCA with scale \(\varepsilon\) 
and axis ratio \(q\), and by \(E^q(X)\) the resulting filtration.
\end{notation}

%Given $X$, a finite metric subspace of $\RR^n$, we denote the Rips-type ellipsoid complex computed via PCA at scale $\varepsilon$ with the ratio $q$ by $E^q_\varepsilon(X)$. The filtration we denote by $E^q(X)$. 
%
   Since $X$ is finite, the filtration $E^q(X)$ consists of finitely many finite simplicial complexes, so $H_k(E^q(X))$ is pointwise finite-dimensional and hence decomposes into interval modules~\cite{CrawleyBoevey15}. Unless otherwise specified we always work with Rips-type ellipsoid complexes as described in the discrete setting subsection. 
We often refer to Rips-type ellipsoid complexes simply as \emph{ellipsoid complexes} and to the corresponding barcodes as \emph{ellipsoid barcodes}.

\subsection{Relation between Rips-Type Ellipsoid and Rips Complexes}

In this subsection we show that Rips-type ellipsoid complexes can be `interleaved' between Rips complexes.\footnote{This is not an interleaving in the sense of Remark~\ref{rem:interleaving}.}
%and by $E^q_\varepsilon$ the ellipsoid complex consisting of ellipsoids whose tangent semi-axes have length $\varepsilon$ and whose ratio of lengths of tangent semi-axes and normal semi axis is $q$, 

\begin{proposition}\label{rips}
Let $X$ be a finite metric subspace of~$\RR^\ad$ (with the metric inherited from $\RR^\ad$). Using Algorithm~\ref{alg:construct_ellipsoids} we construct the Rips-type ellipsoid complex $E^q_\varepsilon(X)$  whose tangent semi-axes have length $\varepsilon$ and whose ratio of lengths of tangent semi-axes and the normal semi-axes is $q$.
We denote by $R_\varepsilon(X)$ the Rips complex at scale $\varepsilon$.
Then the following relation holds
\begin{equation} \label{eq:RE_inclusions}
    R_{\varepsilon / q}(X) \subset E^q_\varepsilon(X) \subset R_\varepsilon(X).
\end{equation}
\end{proposition}
%
% q = a/b
% a = eps
% b = q/eps
%
\begin{proof}
%Let us denote by $b \coloneqq \frac{\varepsilon}{q}$. If two points $p_1$ and $p_2$ are at the distance $2b$ or less, then the ellipsoids $E^q_\varepsilon$ centered at $p_1$ and at $p_2$ will intersect and thus form the 1-simplex $[p_1,p_2] \in E^q_\varepsilon$. In this case, it also holds $[p_1,p_2] \in R_b$ and of course $[p_1,p_2] \in R_\varepsilon$, as the inclusion $R_b \subset R_\varepsilon$ always holds when $\varepsilon \geq b$.
%
%If the distance between $p_1$ and $p_2$ is greater than $2 \varepsilon$, then $[p_1,p_2] \not \in E^q_\varepsilon$, but also $[p_1,p_2] \not \in R_\varepsilon$ and $[p_1,p_2] \not \in R_b$.
%
%If, by contrast, the distance between $p_1$ and $p_2$ is between $2b$ and $2\varepsilon$, then we always have $[p_1,p_2] \in R_\varepsilon$. On the other hand, whether $E^q_\varepsilon$ contains the simplex $[p_1,p_2]$ depends on the construction of the ellipses. The equation~\eqref{eq:RE_inclusions} thus follows.
%
%Figure~\ref{fig:RE_inclusions} shows the relation between ellipsoids used in the construction of the ellipsoid complex and the balls used for Rips complex.

Subsequently, let us denote by $b \coloneqq \frac{\varepsilon}{q}$. 
We first prove that $R_b(X)\subset E^q_\varepsilon(X)$. Let $\sigma \in R_b(X)$. This means that $\forall x, y \in \sigma$, $d(x,y)\leq 2b$. This, in particular, implies that $B_b(x)\cap B_b(y)\neq \emptyset$.  Since $B_b(x)\subset E^q_\varepsilon(x)$ and $B_b(y) \subset E^q_\varepsilon(y)$ it follows that $E^q_\varepsilon(x)\cap E^q_\varepsilon(y) \neq \emptyset$ for all $x, y\in \sigma$. Therefore $\sigma \in E^q_\varepsilon(X)$.

Now we prove that $E^q_\varepsilon(X)\subset R_\varepsilon(X)$. Let $\sigma \in E^q_\varepsilon(X)$. This implies that $E^q_\varepsilon(x)\cap E^q_\varepsilon(y) \neq \emptyset$ for all $x, y\in \sigma$. Let $z\in E^q_\varepsilon(x)\cap E^q_\varepsilon(y)$. Since $E^q_\varepsilon(x)\subset B_\varepsilon(x)$ and $E^q_\varepsilon(y)\subset B_\varepsilon(y)$, it follows by triangle inequality that 
\[
d(x,y) \leq d(x,z)+d(z,y)\leq \varepsilon+\varepsilon=2\varepsilon.
\]
This implies that $d(x,y)\leq 2\varepsilon$ for all $x, y\in \sigma$ and therefore $\sigma\in R_\varepsilon(X)$.
\end{proof}

As an illustration of the previous proof, Figure~\ref{fig:RE_inclusions} shows the relation between ellipsoids used in the construction of the Rips-type ellipsoid complex and the balls used for the Rips complex.

\begin{figure}[h]
\centering
\begin{tikzpicture}
    \draw[dashed] (0,0) ellipse (1 and 1);
    \draw[dashed] (0,0) ellipse (2 and 2);
    \draw[red, very thick] (0,0) ellipse (2 and 1);
    \draw[red, dotted] (-2,0) -- (2,0);
    \draw[red, dotted] (0,-1) -- (0,1);
    \draw [decorate,decoration={brace,amplitude=5pt,mirror,raise=1.1cm}]
  (-2,0) -- (0,0) node[midway,yshift=-1.5cm]{$\varepsilon$};
    \draw [decorate,decoration={brace,amplitude=5pt,mirror,raise=0.1cm}]
  (-2,1) -- (-2,0) node[midway,xshift=-0.8cm]{$b \coloneqq \frac \varepsilon q$};
    \draw [ultra thin, gray] (-2,1) -- (0,1);
    \draw [ultra thin, gray] (-2,-1) -- (-2,0);
    \node at (-1,0) [label={[label distance=0cm, fill=white]90:$B_b(x)$}] {};
    \node at (0,2) [label={[label distance=-0.8cm, fill=white]30:$B_\varepsilon(x)$}] {};
    \node at (2,0) [label={[label distance=0cm, color=red]0:$E^q_\varepsilon$}] {};
    \node[circle, fill, inner sep=0.5pt] at (0,0) [label={[label distance=0cm]0:$x$}] {};
%
   % \begin{scope}[xshift=2.6cm, yshift=2.6cm, rotate=40]
%\draw[dashed] (0,0) ellipse (1 and 1);
    %\draw[dashed] (0,0) ellipse (2 and 2);
   % \draw[red, very thick] (0,0) ellipse (2 and 1);
   % \draw[orange, dotted] (-2,0) -- (2,0);
%    \draw[orange, dotted] (0,-1) -- (0,1);
%    \node[circle, fill, inner sep=0.5pt] at (0,0) [label={[label distance=0cm]0:$y$}] {};
   % \end{scope}
\end{tikzpicture}
\caption{A graphical representation of the nesting property between balls (dashed circles) and ellipsoids (red) that implies the nesting relation between Rips and Rips-type ellipsoid complexes.}
\label{fig:RE_inclusions}
\end{figure}

%% file: sections/algorithm.tex
\label{sec:alg}
In this section we describe the algorithm to compute persistent homology via
Rips-type ellipsoid complexes.
To store the Rips-type ellipsoid complex and calculate its persistent homology,
we use a \emph{simplex tree} data structure~\cite{Boissonat14a} based on the
GUDHI framework~\cite{gudhi:FilteredComplexes}. More specifically, we do the
following:

\begin{algorithm}[H]
  \begin{algorithmic}[1]
    % TODO: Which symbol are we now using for this?
    \Require Point cloud, user-chosen axes ratios of ellipsoids \State Use
    Algorithm~\ref{alg:construct_ellipsoids} to obtain a list of ellipsoids.
    \State $S \gets \emptyset$\Comment{Initialize empty simplex tree}
    \ForEach{point $p$ in the point cloud} \ForEach{point $q$ in the point
      cloud} \State
    $r \gets \texttt{find\_intersection\_radius}(E^q_\bullet(p), E^q_\bullet(q))$\Comment{Find
      radius at which ellipsoids intersect} \State $S \gets S \cup (\{p,q\}, r)$
    \EndFor \EndFor \State $S$.\texttt{expansion()}\Comment{Expand flag complex}
    \State $S$.\texttt{persistence()}\Comment{Calculate barcode}
  \end{algorithmic}
  \caption{Calculating barcodes from Rips-type ellipsoid complexes}
  \label{alg:Ellipsoid barcode}
\end{algorithm}

% \begin{enumerate}
          %     \item Load the point cloud data.
          %     \item Fix user-chosen axes ratios of ellipsoids.
          %     \item Create a list of ellipsoids, each of which is centered at
          %     a point from the point cloud data and oriented according to the
          %     output of the \texttt{fitEllipsoid} algorithm described in
          %     Section~\ref{sec:fitEllipsoid}.
          %     \item Create an empty \texttt{SimplexTree} (part of the Gudhi
          %     library \cite{gudhi:FilteredComplexes}) instance called
          %     \texttt{simplexTreeEllipsoids}.
          %     \item For every two ellipsoids that intersect according to the
          %     \texttt{ellipsoidIntersection} algorithm described in
          %     Section~\ref{sec:ellipsoidIntersection}, add an edge to
          %     \texttt{simplexTreeEllipsoids}.
          %     \item Use the \texttt{expansion()} method on
          %     \texttt{simplexTreeEllipsoids} to obtain higher dimensional
          %     simplices.
          %     \item Use the \texttt{persistence()} method on the expanded
          %     \texttt{simplexTreeEllipsoids} to get the persistence.
          %   \end{enumerate}

          %   \subsubsection{\texttt{fitEllipsoid} algorithm}
          %           \label{sec:fitEllipsoid}

          %           To choose the orientation of each of the ellipsoids, the
          %           following algorithm is used.

          %           For each point $p$ in the point cloud:

          %           \begin{enumerate}
          %             \item Fix user-chosen number $k$ of neighbours.
          %             \item Run principal component analysis on the
  % $k$-nearest neighbours to $p$.
  %\item Align the ellipsoid axes to the eigenvectors obtained from the principal component analysis.
%\end{enumerate}

The complexity of the algorithm depends on the complexity of the algorithm used to first create the ellipsoids, which we earlier determined to be   $\mathcal{O}(Nk\log N + Nn^2k)$. We recall that $N$ refers to the number of sample points, $n$ to their dimension, and $k$ to the number of neighbors used for the tangent-space approximation.
\Cref{alg:Ellipsoid barcode} is thus \emph{prima facie} dominated by the nested loops, which check for all intersections between ellipsoids, for which we use a pre-existing algorithm~\cite{alger2021detect,gilitschenski2012robust} that we outline below for the reader's convenience.
Assuming that this step has \emph{constant} complexity, checking all pairwise intersections has a complexity of $\mathcal{O}(N^2)$.
More precisely, by Subsection~\ref{subseq:intersections}, computing the filtration value of a single edge takes $\mathcal{O}(n^3)$ time, so computing all edge filtration values takes $\mathcal{O}(N^2n^3)$ time, compared with $\mathcal{O}(N^2n)$ for the Vietoris--Rips complex. For a fixed ambient dimension both are $\mathcal{O}(N^2)$. 
The expansion of the flag complex~(executed in the penultimate line) has \emph{output-sensitive complexity} and is trivially upper-bounded by $\mathcal{O}(N^{D+1})$ for expansion up to dimension $D$; see \textcite{Boissonat14a} for a more detailed analysis.
Finally, the barcode calculation takes at most $\mathcal{O}(m^\omega)$ time, where $m$ denotes the size of the resulting flag complex, and $\omega < 2.376$ denotes the matrix multiplication exponent~\cite{Milosavljevic11a}.

Altogether, apart from the last two steps, which are the same as for the Vietoris--Rips complex, \Cref{alg:Ellipsoid barcode} runs in $\mathcal{O}(Nk\log N + Nn^2k + N^2n^3)$ time. 
%
%Since our algorithm shares the last two steps with standard persistent-homology algorithms, improvements of its running time in practice require replacing the intersection checks. We leave this for future work, noting that classical results on improving the performance of rigid-body simulations~\cite{Baraff92a} could potentially be gainfully combined with improved flag complex expansion algorithms~\cite{Zomorodian10a}.
%In comparison to the standard Vietoris--Rips complex expansion based on Euclidean balls,
%our current implementation of the ellipsoids complex is thus less computationally scalable.
%However, we believe that the smaller complexes that one typically obtains with our algorithm are a suitable tradeoff.

\subsection{Intersection of Ellipsoids}\label{subseq:intersections}

We first recall the definition of an ellipsoid given by equation~\eqref{eq:ellipsoid_defn}. In this definition it is assumed that the axes of the ellipsoid are aligned with the coordinate axes and that the ellipsoid is centered at the origin.
We can rewrite equation~\eqref{eq:ellipsoid_defn} as:
\begin{equation}\label{eq:Ellipsoid2}
E^q_\varepsilon(x) = \left\{(x_1, \dots, x_n) \in \mathbb R^n \mid (x_1, \dots, x_n)^T \, \Lambda \, (x_1,\dots,x_n) \leq 1\right\}, 
% \qquad 
% \Lambda = 
% \begin{pmatrix}
%     \varepsilon &        &            &   &        & & \\
%                 & \ddots &            &   &        & &\\
%                 &        &\varepsilon &   &        & & \\
%                 &        &            & b &        & & \\
%                 &        &            &   & \ddots & & \\
%                 &        &            &   &        & b & \\
% \end{pmatrix}
\end{equation}
for $\Lambda=\operatorname{diag}(\lambda_1, \dots, \lambda_n)$, where $\lambda_1 = \dots = \lambda_m = \frac 1 {\varepsilon^2}$ and $\lambda_{m+1}=\dots=\lambda_n = \frac 1 {b^2}$.
To rotate such an ellipsoid so that its axes lie along the orthonormal basis $\{v_1, \dots, v_n\}$, we can apply a rotation matrix $P$ sending the coordinate axes $e_1, \dots, e_n$ to $v_1, \dots, v_n$. In other words, the matrix $P$ is given by 
\[
    P = 
    \begin{pmatrix}
        \vert & & \vert \\
        v_1 & \dots & v_n\\
        \vert & & \vert 
    \end{pmatrix},
\]
so that $P e_i = v_i$.
The equation describing an ellipsoid centered at a point $x \in \mathbb R^{n}$ with the axes given
by the vectors $v_{i} \in \mathbb R^{n}$, $i \in \{1, \dots,n\}$ can thus be
written as follows:
\begin{equation}
  \{y \in \RR^{n} \mid (P^{-1}(y - x))^{T} \, \Lambda \, P^{-1} (y-x) \leq 1 \} \stackrel{(\spadesuit)}{=} \{y \in \RR^{n} \mid (y - x)^{T}  \, P \, \Lambda \, P^T \, (y-x) \leq 1 \}.
\end{equation}
In equality $(\spadesuit)$, we used the fact that rotation matrices are orthogonal, i.e. $P^{-1} = P^T$.

For convenience, we introduce the following notation.
For a symmetric positive definite matrix $A \in \mathbb{R}^{n \times n}$ and a center
$c \in \mathbb{R}^n$, we define
\[
E(A,c) := \{\, y \in \mathbb{R}^n \mid (y - c)^{T} A (y - c) \le 1 \,\}.
\]

%In determining whether two ellipsoids intersect, it is important to keep track of their orientations. 
Thus $E(P\Lambda P^T, x)$ denotes the ellipsoid $E^q_\varepsilon(x)$ with axes lying along the unit vectors $Pe_1, \dots, Pe_n$, where $P$ is a rotation matrix and $e_1, \dots, e_n$ are the coordinate axes. As above, the matrix $\Lambda=\operatorname{diag}(\lambda_1, \dots, \lambda_n)$, where $\lambda_1 = \dots = \lambda_m = \frac 1 {\varepsilon^2}$ and $\lambda_{m+1} = \dots = \lambda_n = \frac {q^2} {\varepsilon^2}$.

%
% % oLD
% We first recall that an ellipsoid centered at the point $p \in \mathbb R^{n}$ with the axes given
% by the vectors $v_{i} \in \mathbb R^{n}$, $i \in \{1, \dots,n\}$ can be
% described as follows:
% \begin{equation}
%   \{x \in \RR^{n} \mid (x - p)^{T} \Sigma (x-p) \leq 1 \},
% \end{equation}
% where $\Sigma$ is a symmetric positive definite matrix whose eigenvectors are
% $v_{i}$ for $i \in \{1,\dots, n\}$ and the eigenvalues are the reciprocals of
% the squares of lengths of the semi-axes. We denote such an ellipsoid by
% $E(\Sigma, \mu)$.
% ---
We use the following result to determine whether two ellipsoids intersect:
\begin{proposition}[{\cite[Proposition~2]{gilitschenski2012robust}}]
  Let $E(A,c)$ and $E(B,d)$ be two ellipsoids (here we use the notation described in the previous paragraph). Denote $v = d - c$ and define
  \begin{equation} \label{eq:K_def}
    K \colon (0,1) \to \RR, \qquad \lambda \mapsto 1 - v^{T} \left( \frac{1}{1-\lambda} B^{-1} + \frac 1 \lambda A^{-1} \right)^{-1} v.
  \end{equation}
  The ellipsoids $E(A,c)$ and $E(B,d)$ intersect if and only if for all
  $\lambda \in (0,1)$ we have $K(\lambda) \geq 0$.
\end{proposition}
We use Proposition~3.1 to compute the filtration value of each edge directly. Let $x,y\in X$ and write $E^q_\varepsilon(x)=E(\varepsilon^{-2}A,x)$ and $E^q_\varepsilon(y)=E(\varepsilon^{-2}B,y)$, where $A$ and $B$ are the matrices of the two ellipsoids at scale $\varepsilon=1$, and set $v=y-x$. Substituting into~\eqref{eq:K_def} gives
\[
    K_\varepsilon(\lambda)=1-\frac{f(\lambda)}{\varepsilon^2},
    \qquad
    f(\lambda)=v^T\Bigl(\frac{1}{1-\lambda}B^{-1}+\frac{1}{\lambda}A^{-1}\Bigr)^{-1}v,
\]
so by Proposition~3.1 the ellipsoids $E^q_\varepsilon(x)$ and $E^q_\varepsilon(y)$ intersect if and only if $\varepsilon^2\ge\max_{\lambda\in(0,1)}f(\lambda)$. Hence the edge $\{x,y\}$ enters the filtration at
\[
    \varepsilon(x,y)=\Bigl(\max_{\lambda\in(0,1)}f(\lambda)\Bigr)^{1/2}.
\]
If $\Phi$ simultaneously diagonalizes $A^{-1}$ and $B^{-1}$, with $\Phi^TB^{-1}\Phi=I$ and $\Phi^TA^{-1}\Phi=\operatorname{diag}(\mu_1,\dots,\mu_n)$, and $w=\Phi^Tv$, then
\[
    f(\lambda)=\sum_{i=1}^n w_i^2\,\frac{\lambda(1-\lambda)}{\lambda+\mu_i(1-\lambda)}.
\]
Since $K_\varepsilon$ is convex~\cite{gilitschenski2012robust}, $f$ is concave, so we can use the golden-section search~\cite{Kiefer1953} to find the maximum. Computing $\Phi$ takes $\mathcal{O}(n^3)$ time~\cite[\S 8.7]{GolubVanLoan2013}, and each evaluation of $f$ takes $\mathcal{O}(n)$ time. With a fixed number of golden-section steps, the filtration value of an edge is computed in $\mathcal{O}(n^3)$ time. For the Vietoris--Rips complex the complexity of this step is $\mathcal{O}(n)$. We thank Sergio Serrano de Haro Iv\'a\~nez for his help in significantly speeding up the algorithm.

%Thus, given two ellipsoids $E(A,c)$ and $E(B,d)$, we find the minimum of the
%function $K$ as defined in equation~\eqref{eq:K_def}. If the minimum is smaller
%than 0, the two ellipsoids intersect.
%
%Since the objective function is convex, the problem is feasible and convergence is guaranteed, enabling the use of efficient optimization procedures~\cite{SciPy}. 

  \subsection{Computation Time Compared to Rips Complexes}

\label{subsec:computation_time}

\begin{figure}[h]
  \centering
  \includegraphics[width=0.7\textwidth]{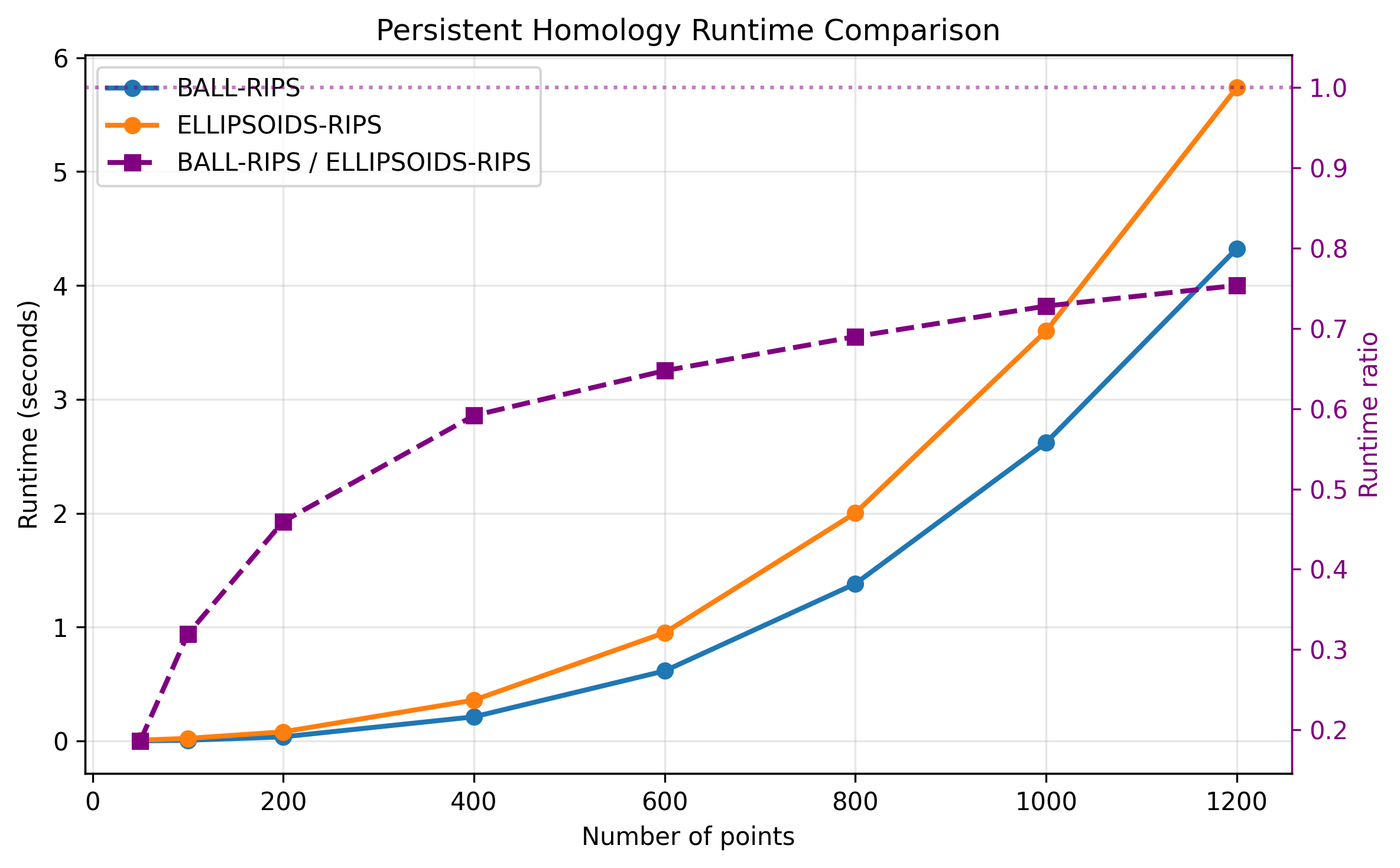}
  \caption{Runtime comparison of persistent homology computations using standard
    Vietoris--Rips complexes (denoted by \texttt{BALL-RIPS} in the legend) and
    ellipsoid complexes (denoted \texttt{ELLIPSOIDS-RIPS}). The horizontal axis
    shows the number of points in the point clouds used for the computations and
    the vertical axes show the runtime (on the left-hand side) and the runtime
    ratio \texttt{BALL-RIPS / ELLIPSOIDS-RIPS} (on the right-hand side).}
  \label{fig:performancevsrips}
\end{figure}

  In Figure~\ref{fig:performancevsrips}, the computation times of Rips-type ellipsoid complexes and standard Vietoris--Rips complexes are shown for a number of differently sized point clouds.
  Both methods have a comparable time complexity, with ellipsoid complexes being slightly slower due to the additional computations needed to orient ellipsoids and compute their intersections.

%% file: sections/stability.tex
\label{sec:stability}
%Fix $X$ to be a finite metric space. For simplicity, we shall assume that $X$ is a subspace of $\RR^n$ with the Euclidean norm. 

Let \(X\) be a finite metric subspace of~\(\mathbb{R}^n\), equipped with the metric inherited from~\(\mathbb{R}^n\). 
The aim of this section is to show that if we perturb \(X\) slightly to obtain a new point cloud \(\tilde{X}\), 
then the persistence diagrams associated to \(X\) and \(\tilde{X}\) remain close in the bottleneck distance. 
To make this precise, we use the notion of \(\delta\)-perturbations, which provide a way to measure the proximity 
between finite subsets of~\(\mathbb{R}^n\).

\begin{definition}[$\delta$-perturbations]\label{def:delta}
    Let $n\in \NN$, $X\subset \RR^n$ be a finite set, and $\delta \geq 0$. We say that an injective function $\mathfrak{p}\colon X\to \RR^n$ is a $\delta$\define{-perturbation} if $\max_{x\in X} \norm{x - \mathfrak{p}(x)} \le \delta$.
\end{definition}
Subsequently, we prove the following theorem.
\begin{theorem}\label{thm:stability}
Let $k, n\in \mathbb{N}$ with $k \ge n$, $X\subset\RR^n$ be a $k$-generic subset (Definition~\ref{def:genericity}), and $\mathfrak{p}\colon X\to \RR^n$ be a $\delta$-perturbation with $\delta \in [0, \delta_u(X))$ for $\delta_u(X)$ defined in Proposition~\ref{prop:stability_basismap}. Then
\[
    \distance_b\!\big(\dgm(E^q(X)), \dgm(E^q(\mathfrak{p}(X)))\big) \leq  Cq\diam(X)\delta %\le Cq\diam(X)\distance_H(X, \mathfrak{p}(X)),
\]
where $C = C(X, \eta_X, k, m, q)$ is from Lemma~\ref{prop:perturb_ell_intersec}.
\end{theorem}
In Proposition \ref{prop:open_and_dense}, we prove that $k$-genericity is not a particularly restrictive condition: $k$-generic point clouds form an open dense subset of $(\mathbb{R}^d)^N$, and any point cloud can be perturbed arbitrarily slightly to satisfy it.

%In this section, we want to show that the ellipsoid complexes are stable. As in \cite{CdSO14}, we will show the interleaving property for the ellipsoid complex. Let $\tilde{X}$ be a slightly perturbed version of $X$. In our case, similar as for the Rips complex, proving the interleaving property then amounts to showing that if two ellipsoids intersect for the point cloud $X$, then the corresponding ellipsoids in $\tilde{X}$ intersect as well, up to slightly increasing the size of the ellipsoids.

\subsection{Notation}

In this subsection we introduce the notation and formally define the objects that appear in Algorithm~\ref{alg:construct_ellipsoids}. 

From now on let
\begin{itemize}
\item
$m\in\{1,\dots,n-1\}$ denote the number of tangent directions that are used to construct the ellipsoids.
\item $\|\cdot\|$ denote the Euclidean norm of a vector and the Frobenius norm of a matrix. 
\end{itemize}

\begin{definition}[Neighborhoods]
    Let $k,n\in \NN$ and let $X\subset \RR^n$ be a finite subset. For each $x\in X$ we can define the preorder $\le_x$ by
    \begin{equation}
        u\le_x v \iff \norm{u - x} \le \norm{v - x}
    \end{equation}
    for all $u, v\in X$. Define $<_x$ in the usual way. Given any $x\in X$ we write 
    \begin{equation}
        N_k(x) \coloneqq \set{y\in X\setminus \{x\}\mid \text{There exist at most $k$ elements $\tilde{x}\in X\setminus \{x\}$ such that $\tilde{x}\le_xy$}}.
    \end{equation}
We call this the \define{set of} $k$\define{-neighbors of }$x$\define{ in }$X$.
\end{definition}

To run the principal component analysis, we have to turn the neighborhoods into matrices.

\begin{definition}[Neighborhood Matrix]\label{def:neighborhood_matrix}
    Let $k, n\in\mathbb{N}$ and let $X\subset \R^n$ be a finite subset. Let $x\in X$ and let $N_k(x) =\{x_1, \dots, x_k\}$. The \define{neighborhood matrix centered at} $x$ is
    \begin{equation}
        \mathbf{N}_k(x) = \begin{pmatrix}
            | &  & |\\
            x_1 - \bar x & \cdots & x_k - \bar x\\
            | &  & |
        \end{pmatrix},
    \end{equation}
    where $\bar x = \frac{1}{k}\sum_{l=1}^kx_l$.
\end{definition}

Note that we do not care about the ordering of the $x_1, \dots, x_k$ with respect to $\le_x$ here. While in principle this makes the definition of the neighborhood matrix ambiguous, it will become apparent in our proofs that the results do not depend on this order.

Recall that, having found the neighborhood matrix $\mathbf{N}_k(x)$ at a point $x\in X$, the next step of Algorithm~\ref{alg:construct_ellipsoids} is to run the principal component analysis on the points in $N_k(x)$. This amounts to computing the singular value decomposition of the matrix $\mathbf{N}_k(x)$. 

\begin{definition}[Ellipsoid Basis]\label{def:ellipsoid_basis}
    Let $k, n\in\NN$ and let $X\subset \R^n$ be a finite subset. For $x\in X$ suppose $\mathbf{N}_k(x) = U\Sigma V^T$ is a singular value decomposition of the neighborhood matrix, i.e., $U$ and $V$ are orthogonal matrices and $\Sigma$ is a positive-semidefinite diagonal matrix. We denote
    \[
        \mathbf{U}(x) \coloneqq U.
    \]
\end{definition}

We shall note here that we will repeatedly assume that $k\ge n$. This ensures that $U$ in the definition above is of size $n\times n$ and hence defines an orthogonal basis for the ambient space $\R^n$.
The construction of an ellipsoid centered at $x$ can then be completed as in Equation~\eqref{eq:Ellipsoid2} by setting
\begin{equation}
        E_{\eps}^q(x)\coloneqq \set{y\in\RR^n \mid (y - x)^T \mathbf{U}(x)\Lambda \mathbf{U}(x)^T(y - x)\leq 1 }.
\end{equation}

Note that $\mathbf{U}(x)$ is not unique. However, the ellipsoid depends only on the orthogonal projection $\mathbf{U}_m(x)\mathbf{U}_m(x)^T$ onto the span of the first $m$ left singular vectors, since
\[
\mathbf{U}(x)\Lambda\mathbf{U}(x)^T=\varepsilon^{-2}\,\mathbf{U}_m(x)\mathbf{U}_m(x)^T+b^{-2}\bigl(I-\mathbf{U}_m(x)\mathbf{U}_m(x)^T\bigr),
\]
and this projection is well defined whenever $\sigma_m(\mathbf{N}_k(x))>\sigma_{m+1}(\mathbf{N}_k(x))$.

\subsection{Obstructions to Gromov-Hausdorff Stability}\label{sub:obstruction}

To motivate some definitions and choices that we make in this section, we give a short exposition on two types of behavior that can cause instability of the ellipsoid construction under small perturbations of the point cloud $X$ with respect to the standard Gromov-Hausdorff distance.

Intuitively, by stability we mean that the ellipsoids do not undergo abrupt changes in orientation under sufficiently small perturbations. In Chazal, de Silva and Oudot~\cite{CdSO14}, stability is proved with respect to so-called correspondences,\footnote{A multivalued map \(C \colon X \rightrightarrows Y\) (i.e., a subset \(C \subseteq X\times Y\) whose projection to \(X\) is surjective) is called a \emph{correspondence} (Definition 4.1~\cite{CdSO14}) if the canonical projection \(C \to Y\) is also surjective. Or equivalently, if its transpose
\[
C^{T} \;:=\; \{\, (y,x)\in Y\times X \;:\; (x,y)\in C \,\}
\]
is a multivalued map \(Y \rightrightarrows X\).} which allow not only small perturbations of the points in $X$, but also the addition of new points that are close to $X$. Since our construction of ellipsoids is sensitive to the arrangement of the points in space, the addition of sufficiently many nearby points can trigger a sudden swap in the orientation of the ellipsoid. Indeed, this can be observed in Figure~\ref{fig:degeneration_by_addition}. We will refer to this phenomenon as \textit{instability under augmentation}.
\begin{figure}[h]
\centering
\begin{tikzpicture}[scale=1.5]
    \filldraw (0,0) circle [radius=1pt];
   
    \filldraw (38:0.95) circle [radius=1pt];
    \filldraw (25:0.92) circle [radius=1pt];
    \filldraw (30:0.9) circle [radius=1pt];
    \filldraw (20:0.89) circle [radius=1pt];
    
    \filldraw (84:0.9) circle [radius=1pt];

    \draw[dashed] (0,0) circle (1);

    \draw[->] (1.5, 0) -- (2.5, 0) node[midway, above] {$\mathfrak{p}$};
    
    \begin{scope}[xshift=4cm]
    \filldraw (0,0) circle [radius=1pt];
   
    \filldraw (38:1.06) circle [radius=1pt];
    \filldraw (25:1.08) circle [radius=1pt];
    \filldraw (30:1.1) circle [radius=1pt];
    \filldraw (20:1.05) circle [radius=1pt];
    
    \filldraw (98:0.94) circle [radius=1pt];
    \filldraw (84:0.95) circle [radius=1pt];
    \filldraw (92:0.9) circle [radius=1pt];
    \filldraw (80:0.89) circle [radius=1pt];
    \filldraw (106:0.86) circle [radius=1pt];
    
    \draw[dashed] (0,0) circle (1);
    \end{scope}
\end{tikzpicture}
\caption{Schematic depiction showing that adding multiple nearby points to a single point can lead to a swapping of the axes.}\label{fig:degeneration_by_addition} 
\end{figure}
Because of this behavior, we cannot allow general perturbations in the Gromov–Hausdorff sense, but must restrict ourselves to perturbations that move each point of the cloud $X$ by at most a fixed distance and do not allow for the addition of extra points.

A second problematic arrangement that may cause an ellipsoid to suddenly switch orientation is when many points leave the $k$-neighborhood $N_k(x)$ in one direction, while other points enter the neighborhood from another direction. This situation is shown in \Cref{fig:swapping}, and it would result in the constructed ellipsoids flipping by about $60^\circ$. Although such behavior cannot occur if the perturbation is chosen sufficiently small~(see Lemma~\ref{lem:stability_neighborhoods}), one can nevertheless construct sequences of examples in which the maximum admissible perturbation size tends to zero. We refer to this phenomenon as \emph{neighborhood-swap instability}.

\begin{figure}[h]
\centering
\begin{tikzpicture}[scale=1.5]
    \filldraw (0,0) circle [radius=1pt];
   
    \filldraw (38:0.95) circle [radius=1pt];
    \filldraw (25:0.92) circle [radius=1pt];
    \filldraw (30:0.9) circle [radius=1pt];
    \filldraw (20:0.89) circle [radius=1pt];
    
    \filldraw (98:1.06) circle [radius=1pt];
    \filldraw (84:1.08) circle [radius=1pt];
    \filldraw (92:1.1) circle [radius=1pt];
    \filldraw (80:1.11) circle [radius=1pt];
    
    \draw[dashed] (0,0) circle (1);

    \draw[->] (1.5, 0) -- (2.5, 0) node[midway, above] {$\mathfrak{p}$};
    
    \begin{scope}[xshift=4cm]
    \filldraw (0,0) circle [radius=1pt];
   
    \filldraw (38:1.06) circle [radius=1pt];
    \filldraw (25:1.08) circle [radius=1pt];
    \filldraw (30:1.1) circle [radius=1pt];
    \filldraw (20:1.05) circle [radius=1pt];
    
    \filldraw (98:0.9) circle [radius=1pt];
    \filldraw (84:0.92) circle [radius=1pt];
    \filldraw (92:0.9) circle [radius=1pt];
    \filldraw (80:0.89) circle [radius=1pt];
    
    \draw[dashed] (0,0) circle (1);
    \end{scope}
\end{tikzpicture}
\caption{Exemplary situation where the swapping degeneracy occurs after applying a $\delta$-perturbation $\mathfrak{p}$. The dashed circles are a visual aid to distinguish the neighborhood relations.}\label{fig:swapping}
\end{figure}
\par 
One might also wonder whether changing the neighborhood size could be a viable strategy. We will illustrate that this also leads to the swapping of the axes. Indeed, let us give a probabilistic example for this type of behavior. 
\begin{example}
Suppose we are given two points $x_1 \coloneqq (0,0), x_2 \coloneqq (a,0)$ with $a>0$. We want to determine the critical number of points that must be added for one of the PCA axes to flip. To this end, we present an example where we add $l$ points with random $ y$-coordinates around the origin. This will ensure that the axes from PCA are situated on the $x$- and $y$-axis. Initially, without the perturbation, the major axis will be situated on the $x$-axis, since $\mathbf{N}_2(x)\mathbf{N}_2(x)^T=\operatorname{diag}(a^2/2,\,0)$. We shall compute the smallest $l$ for which the major axis of the expected matrix $\mathbb{E}[\mathbf{N}_{l+2}(x)\mathbf{N}_{l+2}(x)^T]$ flips to the $y$-axis.

Suppose we are adding $l\in\NN$ points $p_i \coloneqq(0, w_i)$ with independent and identically distributed $w_i \sim \mathcal{N} (0,\sigma^2)$, which shall model the perturbation of the point $x_1$. The neighborhood matrix is given by the following $2\times (l+2)$-matrix 
\begin{align}
    \mathbf{N}_{l+2}(x) \coloneqq \begin{pmatrix}
        | & | & | &  & |\\
        x_1 - \bar{x} & x_2 -\bar{x} & p_1 - \bar{x} & \cdots & p_l - \bar{x}\\
        | & | & | &  & |
    \end{pmatrix}
\end{align}
where $\bar{x}$ is the sample mean
\[
    \bar{x} = \frac{1}{l+2}\left(x_1 + x_2 + \sum_{i=1}^l p_i\right) = \frac{1}{l+2}\begin{pmatrix}
        a\\\sum_{i=1}^l w_i
    \end{pmatrix}
\]
PCA directions are eigenvectors of $\mathbf{N}_{l+2}(x)\mathbf{N}_{l+2}(x)^T$.
Note that
\begin{equation}
    \Ev[\mathbf{N}_{l+2}(x)\mathbf{N}_{l+2}(x)^T] = \begin{pmatrix}
        \lambda_l^1 & 0\\
        0 & \lambda_l^2
    \end{pmatrix},
\end{equation}
for some $\lambda_l^1,\lambda_l^2\ge 0$. Indeed, we may calculate the off-diagonal term of 
\begin{align*}
\bigl(\mathbf N_{l+2}(x)\mathbf N_{l+2}(x)^T\bigr)_{12}
&=(x_1^1-\bar{x}^1)(x_1^2-\bar{x}^2)
  +(x_2^1-\bar{x}^1)(x_2^2-\bar{x}^2)
  +\sum_{i=1}^l (p_i^1-\bar{x}^1)(p_i^2-\bar{x}^2)\\
&=\Bigl(0-\frac{a}{l+2}\Bigr)\Bigl(0-\frac{\sum_{j=1}^l w_j}{l+2}\Bigr)
 +\Bigl(a-\frac{a}{l+2}\Bigr)\Bigl(0-\frac{\sum_{j=1}^l w_j}{l+2}\Bigr)\\
&\qquad
 +\sum_{i=1}^l \Bigl(0-\frac{a}{l+2}\Bigr)
                 \Bigl(w_i-\frac{\sum_{j=1}^l w_j}{l+2}\Bigr)\\[4pt]
&=\frac{a}{(l+2)^2}\sum_{j=1}^l w_j
  -\frac{a(l+1)}{(l+2)^2}\sum_{j=1}^l w_j
  -\frac{a}{l+2}\sum_{i=1}^l w_i
  +\frac{al}{(l+2)^2}\sum_{j=1}^l w_j\\[2pt]
&=-\frac{a}{l+2}\sum_{i=1}^l w_i.
\end{align*}
Taking the expectation implies:
\[\Ev\Bigl[-\frac{a}{l+2} \sum_{i=1}^l \cdot w_i\Bigl]=-\frac{a}{l+2} \sum_{i=1}^l \Ev[w_i]=0, \]
since $w_i$ are normally distributed around zero.
As mentioned above, we have that $\lambda_0^1 > \lambda_0^2$ and we would like to find the critical $l_{\crit}$ so that for all $l\ge l_{\crit}$ we have $\lambda_l^1 \leq \lambda_l^2$.

\begin{proposition}[Flipping the axes]
    Let $\lambda_l^1, \lambda_l^2$ be the eigenvalues of $\Ev[\mathbf{N}_{l+2}(x)\mathbf{N}_{l+2}(x)^T]$. Then we have 
    \begin{align*}
        \lambda_l^2 \geq \lambda_l^1 \iff l\geq \Bigl \lceil \frac{a^2}{\sigma^2}\Bigl \rceil \eqqcolon l_{\crit}.
    \end{align*}
\end{proposition}
\begin{proof}
    As calculated in the above, the mean is given by:
    \begin{align}
        \bar{x} = \begin{pmatrix}\bar{x}^1\\ \bar{x}^2\end{pmatrix} =\frac{1}{l+2}\begin{pmatrix}a\\ \sum_{i=1}^l w_i\end{pmatrix}.
    \end{align}
   Furthermore, since its expectation is diagonal we may calculate the eigenvalues directly (with the calculation being as in the off diagonal case):
    \begin{align*}
        (\mathbf{N}_{l+2}(x)\mathbf{N}_{l+2}(x)^T)_{11} = \frac{l+1}{l+2}a^2,
    \end{align*}
     and
    \begin{align*}
     (\mathbf{N}_{l+2}(x)\mathbf{N}_{l+2}(x)^T)_{22} = 2(\bar{x}^2)^2+\sum_{i=1}^{l}(w_i - \bar{x}^2)^2 = 2(\bar{x}^2)^2+\sum_{i=1}^{l}\left(w_i^2 - 2\bar{x}^2w_i +(\bar{x}^2)^2\right).
    \end{align*}
    Since each of the $w_i$ has zero mean $\sum_{i = 1 }^l \mathbb{E}[w_i] = 0.$
    Moreover, we have
    \begin{align}
       \Ev[(\bar{x}^2)^2] = \frac{1}{(l+2)^2}\Ev[(\sum_{i=1}^l w_i)^2] = \frac{1}{(l+2)^2}\!\left(
\sum_{i=1}^l \mathbb{E}[w_i^2]
+ 2\!\!\sum_{1\le i<j\le l}\!\! \mathbb{E}[w_i w_j]
\right)=\frac{l\sigma^2}{(l+2)^2}
    \end{align}    
    since 
    \begin{align}
        \mathbb{E}\!\left[\sum_{i=1}^l w_i^2\right] = \sum_{i=1}^l\mathbb{E}[w_i^2] = l \sigma^2.
    \end{align}
    Thus, we can see that
    \begin{align}
    \begin{aligned}
        \lambda_l^2 &= \mathbb{E}[(\mathbf{N}_{l+2}(x)\mathbf{N}_{l+2}(x)^T)_{22}]\\
        &= 2\mathbb{E}[(\bar{x}^2)^2]+\sum_{i=1}^{l}\mathbb{E}[w_i^2] - 2\mathbb{E}[\bar{x}^2w_i] +\mathbb{E}[(\bar{x}^2)^2]\\
        &= \sum_{i=1}^{l}\mathbb{E}[w_i^2]
   - 2\,\mathbb{E}\!\left[\bar{x}^2 \sum_{i=1}^{l} w_i\right]
   + (l+2)\,\mathbb{E}\!\left[(\bar{x}^2)^2\right]\\
        &=\frac{l\sigma^2}{(l+2)^2}\left( (l+2)^2 -2(l+2) +(l+2) \right)\\
        &=\frac{l(l+1)}{l+2}\sigma^2.
    \end{aligned}
    \end{align}
Here we used that
\[
\mathbb{E}\!\left[\bar{x}^2 \sum_{i=1}^l w_i\right]=\frac{1}{l+2}\mathbb{E}\!\left[\left(\sum_{i=1}^l w_i\right)^2\right]=\frac{l\sigma^2}{l+2}.
\]    
    Hence we see that $\lambda_l^2 \geq \lambda_l^1$ if and only if $l\geq \frac{a^2}{\sigma^2}$.
\end{proof}

\begin{figure}
\begin{tikzpicture}[scale=1.0]
    \fill (0,0) circle (2pt) node[below left] {$x_1$};
    \fill (3,0) circle (2pt) node[below left] {$x_2$};
    
    \fill (8,0) circle (2pt) node[below left] {$x_1$};
    \fill (11,0) circle (2pt) node[below left] {$x_2$};
    
    \fill (8, 0.6) circle (2pt) node[above left] {$p_1$};
    \fill (8, 0.4) circle (2pt) node[left] {$p_2$};
    \draw[->, thick, opacity=0.5] (8,0) -- ++(1,0);
    \draw[->, thick, opacity=0.5] (0,0) -- ++(1,0);
    \draw[->, thick, opacity=0.5] (8,0) -- ++(0,1);
\end{tikzpicture}
\caption{Schematic depiction where the arrows represent the axes of the ellipsoids. It shows that adding enough points will make the axis flip, i.e. making the shorter axis the longer one and vice versa.}\label{fig:axis_flipping} 
\end{figure}
\end{example}

As these examples demonstrate, to formulate a stability statement we restrict both the class of point clouds and the admissible correspondences.  Instability under augmentation forces us to restrict to correspondences which we call $\delta$-{\emph perturbations} (Definition~\ref{def:delta}). To avoid the neighborhood-swap instability, we work with \(k\)-generic point clouds.

\begin{definition}[$k$-genericity]\label{def:genericity}
    Let $k, n\in \NN$ and $k\ge n$.    
    We say that $X\subset\mathbb{R}^n$ is \textbf{$k$-generic} if it is finite with $|X|\ge k+2$ and satisfies the following two conditions
    \begin{itemize}
        \item \emph{Separation of points:} For all $x\in X$ we have that
            $N_k(x)$ contains exactly $k$ points and for every two distinct $x_1, x_2\in
            N_k(x)$ it holds $x_1 <_x x_2$ or $x_2 <_x x_1$.
        \item \emph{$m$-spectral gap:} There exists $\eta>0$ such that for every $x\in X$,
\[
\sigma_m\bigl(\mathbf N_k(x)\bigr)-\sigma_{m+1}\bigl(\mathbf N_k(x)\bigr)>\eta,
\]
where $m$ is the number of tangent directions (fixed in Subsection 4.1) and $\sigma_1(\cdot)\ge \cdots \ge \sigma_n(\cdot)\ge 0$ are the singular values.
    \end{itemize}
\end{definition}
%Note that these are the minimal assumptions under which we can study the ellipsoid homology. Since the first assumption guarantees that the $k-$nearest neighbours are well defined, and the second assumption corresponds to the manifold hypothesis.
%Let's not draw attention too much to the manifold hypothesis :D
In Proposition \ref{prop:open_and_dense}, we prove that $k$-genericity is not a particularly restrictive condition: $k$-generic point clouds form an open dense subset of $(\mathbb{R}^d)^N$, and any point cloud can be perturbed arbitrarily slightly to satisfy it.

\subsection{Proof of Stability}

\noindent
We prove stability in four steps. First, in Subsection~\ref{sub:nbhd_stability} we rule out neighborhood-swap phenomena. In particular, we prove that for a $k$-generic finite set $X\subset\RR^n$ there is a threshold $\delta_N(X)>0$ such that every $\delta$-perturbation $\mathfrak p\colon X\to\RR^n$ with $\delta<\delta_N(X)$ preserves $k$-nearest neighborhoods in the sense that $\mathfrak p(N_k(x))=N_k(\mathfrak p(x))$ for all $x\in X$ (Lemma~\ref{lem:stability_neighborhoods}). This lets us compare neighborhood matrices up to a permutation of columns.

We find (after reindexing columns by a suitable permutation matrix $P=P(x)$) that
\[
\|\mathbf{N}_k(x)-\mathbf{N}_k(\mathfrak{p}(x))P\|_{1,2}\le 2\delta.
\]
Intuitively, multiplying by the permutation matrix $P$ reorders the columns of $\mathbf N_k(\mathfrak p(x))$, so it does not change the singular values of the matrix. Combining the elementary bound $\|A\|\le k^{1/2}\|A\|_{1,2}$ and the previous bound on neighborhood matrices yields
\[
\max_i\bigl|\sigma_i(\mathbf N_k(x))-\sigma_i(\mathbf N_k(\mathfrak p(x)))\bigr|\le 2k^{1/2}\delta,
\]
as stated in Proposition~\ref{prop:stability_axesmap}. This is the main result of Subsection~\ref{sub:stability_singular}. In the following subsection we prove a stability statement for the \emph{bases} used in the ellipsoid construction. Writing $\mathbf U(x)$ for the matrix of left singular vectors of $\mathbf N_k(x)$, a spectral gap assumption encoded in $k$-genericity implies that for all $\delta<\delta_u(X)$ one has
\[
        \|P(x)-P(\mathfrak{p}(x))\| \le \frac{12\sqrt{2k}}{\eta_X}\,\delta\quad \forall x\in X.
\]
where $\eta_X$ is the minimal $m$ singular-value gap across $x\in X$ and 
\[
        \mathbf U_m(x)\coloneqq \bigl(\mathbf u_1(x)\ \cdots\ \mathbf u_m(x)\bigr),\qquad
        P(x)\coloneqq \mathbf U_m(x)\mathbf U_m(x)^T 
    \]
(See Proposition~\ref{prop:stability_basismap}). 

We use these results to establish that filtrations $E^q(X)$ and $E^q(\mathfrak p(X))$ are $(C\varepsilon_{\max}(X)\delta)$–interleaved (Theorem~\ref{cor:Ellipsoid_interleaving}). Finally, choosing $\varepsilon_{\max}(X)=q\,\diam(X)$ and invoking the standard stability for interleaved persistence modules yields
\[
\distance_b\bigl(\dgm(E^q(X)),\dgm(E^q(\mathfrak p(X)))\bigr)\le C\,q\,\diam(X)\,\delta,
\]
which is our main stability statement~(Theorem~\ref{thm:ellipsoid_barcode_stability}).

\subsubsection{Neighborhood Stability}\label{sub:nbhd_stability}
The next statement addresses the issue of neighborhood-swap instability. More precisely, we will show that this scenario does not arise for $\delta$-perturbations for sufficiently small $\delta$.
\begin{lemma}[Stability of neighborhoods]\label{lem:stability_neighborhoods}
    Let $k, n\in \mathbb{N}$ and let $X\subset \RR^n$ be a $k$-generic finite subset. There exists a $\delta_N(X) > 0$ such that for all $\delta \in [0,\delta_N(X))$ and any $\delta$-perturbation $\mathfrak{p}\colon X\to \RR^n$ 
    \begin{equation}\label{eq:neighborhood_stable}
        \mathfrak{p}(N_k(x)) = N_k(\mathfrak{p}(x))
\end{equation}
for all $x\in X$.
\end{lemma}
\begin{proof}
    Let $x\in X$ and define the following function on $X$ measuring the distance to membership of the neighborhood $N_k(x)$ as
    \[
        f_x(\tilde{x})\coloneqq \min_{x'\in \set{x}\cup N_k(x)}|\norm{\tilde{x} - x} - \norm{x' - x}|
    \]
    Note that for all $x\in X$ the function $f_x$ is identically $0$ on $\set{x}\cup N_k(x)$ and strictly positive on the complement $X\setminus(\set{x}\cup N_k(x))$ because $X$ is $k$-generic.
    
    We claim that the number
    \begin{align}\label{eq:delta_N}
        \delta_N(X)\coloneqq \frac{1}{4}\min_{x\in X}\min_{\tilde{x}\in X\setminus(\set{x}\cup N_k(x))}f_x(\tilde{x})
    \end{align}
    has the desired properties.

    First, we show that $\delta_N(X) > 0$. Indeed, the numbers in the finite set
    \[
        \set*{\frac{1}{4}\min_{\tilde{x}\in X\setminus(\set{x}\cup N_k(x))}f_x(\tilde{x})\mid x\in X}
    \]
    are all strictly positive and $\delta_N(X)$ is just the minimum of this set. Hence, $\delta_N(X) > 0$.

Let $\delta\in [0,\delta_N(X))$ and let $\mathfrak{p}\colon X\to \RR^n$ be a $\delta$-perturbation.  Let us prove \eqref{eq:neighborhood_stable} by contradiction. So, assume  $\mathfrak{p}(N_k(x)) \neq N_k(\mathfrak{p}(x))$. 

%Then there exist an $x\in X$, $x_1\in X\setminus (\set{x}\cup N_k(x))$ and $x_2\in N_k(x)$ such that $\mathfrak{p}(x_1)\in N_k(\mathfrak{p}(x))$ and $\mathfrak{p}(x_2)\in \mathfrak{p}(X)\setminus N_k(\mathfrak{p}(x))$.
Since $|N_k(\mathfrak p(x))|\le k=|\mathfrak p(N_k(x))|$, there is an $x_2\in N_k(x)$ with
$\mathfrak p(x_2)\notin N_k(\mathfrak p(x))$. Hence at least $k+1$ points of
$\mathfrak p(X)\setminus\{\mathfrak p(x)\}$ lie at distance at most
$\|\mathfrak p(x_2)-\mathfrak p(x)\|$ from $\mathfrak p(x)$. As $|N_k(x)|=k$, one of them is
$\mathfrak p(x_1)$ for some $x_1\in X\setminus(\{x\}\cup N_k(x))$. Thus
    \[
        0\le \norm{\mathfrak{p}(x_2) - \mathfrak{p}(x)} - \norm{\mathfrak{p}(x_1) - \mathfrak{p}(x)}.
    \]
    By the reverse triangle inequality and the fact that $\mathfrak{p}$ is a $\delta$-perturbation we see
    \[
    \begin{aligned}
        \norm{\mathfrak{p}(x_1) - \mathfrak{p}(x)} &= \norm{\mathfrak{p}(x_1) - x_1 + x_1 - x + x - \mathfrak{p}(x)} \ge \norm{\mathfrak{p}(x_1) - x_1 +x_1 - x} - \norm{\mathfrak{p}(x)- x} \\& \ge \norm{x_1 - x} - \norm{\mathfrak{p}(x_1)- x_1} - \norm{\mathfrak{p}(x)- x}\ge \norm{x_1 - x} - 2\delta,
    \end{aligned}
    \]
    similarly the triangle inequality implies:
    \[
    \begin{aligned}
        \norm{\mathfrak{p}(x_2) - \mathfrak{p}(x)} &= \norm{\mathfrak{p}(x_2) - x_2 + x_2 - x + x - \mathfrak{p}(x)} \le \norm{x_2 - x} + 2\delta.
    \end{aligned}
    \]
    Thus, we have found
    \[
        0\le \norm{\mathfrak{p}(x_2) - \mathfrak{p}(x)} - \norm{\mathfrak{p}(x_1) - \mathfrak{p}(x)}\le \norm{x_2 - x} - \norm{x_1 - x} + 4\delta.
    \]
     Now we bound $\norm{x_2 - x} - \norm{x_1 - x}$.  Since $x_1 \in X\setminus N_k(x)$ and $X$ is $k$-generic, we have
\[
    \|x_1 - x\| > \|x' - x\| \quad \text{for all } x'\in \{x\}\cup N_k(x).
\]
Hence, for such $x'$,
\[
    f_x(x_1) = \min_{x'\in  \set{x}\cup N_k(x)} 
    \bigl(\|x_1 - x\| - \|x' - x\|\bigr),
\]
and in particular, for $x_2\in N_k(x)$ we have
\[
    f_x(x_1) \le \|x_1 - x\| - \|x_2 - x\|.
\]
Rearranging yields
\[
    \|x_2 - x\| - \|x_1 - x\| \le -\, f_x(x_1).
\]
Using the definition of $\delta_N(X)$ we deduce that
\[
    f_x(x_1)\ge 4\,\delta_N(X) \textrm{ or equivalently} -f_x(x_1)\leq -4\,\delta_N(X).
\]
Combining the inequalities gives
\[
    \|x_2 - x\| - \|x_1 - x\| 
    \le -\, f_x(x_1)
    \le -\,4\,\delta_N(X),
\]
as claimed.

Putting all the inequalities together we get
    \[
        0\le \norm{\mathfrak{p}(x_2) - \mathfrak{p}(x)} - \norm{\mathfrak{p}(x_1) - \mathfrak{p}(x)}\le \norm{x_2 - x} - \norm{x_1 - x} + 4\delta \leq  4(\delta - \delta_N(X)) < 0,
    \]
    which is a contradiction. We have thus shown \eqref{eq:neighborhood_stable}.
\end{proof}
\begin{remark}
    Note that $\delta_N(X)$ in \eqref{eq:delta_N} can be quite small, which would mean that stability only holds for small perturbations. Indeed, the requirement that the neighborhoods do not change much under the perturbation seems to be quite a strict one. In order to avoid this, it might be feasible to use probabilistic arguments when one assumes that the samples $X$ are sufficiently well distributed.
\end{remark}
\subsubsection{Stability of the Singular Values}\label{sub:stability_singular}
The goal of this subsection is to show that the singular values of the neighborhood matrices depend continuously on the data (Proposition~\ref{prop:stability_axesmap}) and that
\[
|\sigma_i(\mathbf N_k(x))-\sigma_i(\mathbf N_k(\mathfrak p(x)))|
\le 2k^{1/2}\,\delta
\quad\text{for all }i=1,\dots,n\text{ and }x\in X.
\]
Here $\mathfrak p$ is a $\delta$-perturbation, $\delta\in [0, \delta_N(X))$ and $\delta_N(X)$ is defined as in Lemma~\ref{lem:stability_neighborhoods}. 

We first recall some standard definitions.
\begin{definition}[Column $2$-norm]
    Let $k,n\in \NN$. We define the \define{column} $2$\define{-norm} of a matrix $A\in \R^{n\times k}$ by
    \[
        \norm{A}_{1, 2}\coloneqq \max_{j=1,\dots, k}\norm{Ae_j},
    \]
    where $\norm{\bullet}$ denotes the usual Euclidean norm, and $e_j$ is the usual standard basis vector with $j$-th entry $1$ and the rest $0$.
\end{definition}

\begin{comment}
This norm can be related to the usual operator norm via Cauchy-Schwarz.

\begin{proposition}\label{prop:norm_comp}
    Let $k, n\in \NN$ and $A\in \R^{n\times k}$. Then
    \[
        \norm{A}_{1, 2} \le \norm{A}\le k^{1/2} \norm{A}_{1, 2}
    \]
\end{proposition}
\begin{proof}
The first inequality holds by the definition of $\norm{A}_{1,2}$ and $\norm{A}$. For the second, observe that for all $x\in\R^k$
\[
    \norm{Ax} = \norm*{\sum_{j=1}^k x_j Ae_j} \le \sum_{j=1}^n \absval{x_j}\norm{Ae_j} \le \norm{A}_{1,2}\norm{x}_1
\]
and by Cauchy-Schwarz
\[
    \norm{x}_1 \le k^{1/2}\norm{x}.
\]
\end{proof}

\begin{definition}[Hausdorff distance]
    Let $(M,d)$ be a metric space and $X, Y$ two subsets of $M$. We define the Hausdorff distance between $X$ and $Y$ as:
    \[\distance_H(X,Y) \coloneqq \max \{\sup_{x\in X}d(x,Y),\sup_{y\in Y}d(X,y) \}.\]
\end{definition}
\end{comment}

The norm is related to the Frobenius norm as follows.
\begin{proposition}\label{prop:norm_comp}
Let $k,n\in\mathbb{N}$ and $A\in\mathbb{R}^{n\times k}$. Then
\[
\|A\|_{1,2}\le\|A\|\le k^{1/2}\|A\|_{1,2}.
\]
\end{proposition}
\begin{proof}
By definition of the Frobenius norm we have
\[
\|A\|_{1,2}^2=\max_{j=1,\dots,k}\|Ae_j\|^2\le\sum_{j=1}^k\|Ae_j\|^2=\|A\|^2\le k\max_{j=1,\dots,k}\|Ae_j\|^2=k\,\|A\|_{1,2}^2.
\]
Taking square roots gives both inequalities.
\end{proof}

\begin{corollary}[Stability of Neighborhood Matrices]\label{cor:stability_distance_mats}
    Let $k, n \in \NN$ and let $X\subset \RR^n$ be a $k$-generic subset. We define $\delta_N(X)$ as in Lemma~\ref{lem:stability_neighborhoods}. Then for each $x\in X$ and all  $\delta$-perturbations $\mathfrak{p}\colon X\to \R^n$ with $\delta\in [0, \delta_N(X))$ there exists a permutation matrix $P=P(x)\in \R^{k\times k}$ such that we have
    \[
        \norm{\mathbf{N}_k(x) - \mathbf{N}_k(\mathfrak{p}(x))P}_{1, 2}\le 2\delta
    \]
   
\end{corollary}
\begin{proof}
Let $x\in X$, $N_k(x) = \set{x_1, \dots, x_k}$ and $N_k(\mathfrak{p}(x)) = \set{y_1, \dots, y_k}$. By Lemma \ref{lem:stability_neighborhoods} we know that
    \[
        \mathfrak{p}(N_k(x)) = N_k(\mathfrak{p}(x)).
    \]
    Thus, for each $x_i\in N_k(x)$ there exists a $y_j\in N_k(\mathfrak{p}(x))$ such that $y_j = \mathfrak{p}(x_i)$. Denote by $\tau$ the permutation that sends $i\in \{1, \dots ,k\}$ to $j\in \{1, \dots ,k\}$ if $y_j = \mathfrak{p}(x_i)$. Let $P$ denote permutation matrix corresponding to $\tau$ which sends $e_i$ to $e_{\tau(i)}$. Then, for each $m\in\{1, \dots ,k\}$
   \begin{align*}
\|(\mathbf{N}_k(x)-\mathbf{N}_k(\mathfrak{p}(x))P)e_m\|
&=\|x_m-\bar{x}-(y_{\tau(m)}-\bar{y})\|
\le\|x_m-\mathfrak{p}(x_m)\|+\|\bar{x}-\bar{y}\|\\
&\le\delta+\frac{1}{k}\sum_{i=1}^k\|x_i-\mathfrak{p}(x_i)\|\le 2\delta,
\end{align*}
where $\bar{x}=\frac{1}{k}\sum_{i=1}^k x_i$ and $\bar{y}=\frac{1}{k}\sum_{j=1}^k y_j=\frac{1}{k}\sum_{i=1}^k\mathfrak{p}(x_i)$.
\end{proof}

\begin{proposition}[Stability of the Singular Values]\label{prop:stability_axesmap}
    Let $k, n \in \mathbb{N}$, $X\subset \RR^n$ be a $k$-generic subset, and let $\delta_N(X)$ be as in Lemma \ref{lem:stability_neighborhoods}. Then, for all $\delta\in [0,\delta_N(X))$ and any $\delta$-perturbation $\mathfrak{p}\colon X\to \RR^n$ of $X$ we have that
    \[
        \max_{i=1,\dots,  n}\absval{\sigma_i(\mathbf{N}_k(x)) - \sigma_i(\mathbf{N}_k(\mathfrak{p}(x)))} \le 2k^{1/2} \delta\quad \forall x\in X,
    \]
   where $\sigma_i$ denotes the $i$-th singular value of $\mathbf N_k(x)$.
\end{proposition}
\begin{proof}
    Let $x \in X$, $\delta\in [0,\delta_N(X))$ and $\mathfrak{p}\colon X\to \RR^n$ be a $\delta$-perturbation. By Corollary \ref{cor:stability_distance_mats} we can then find a permutation matrix $P\in\R^{k\times k}$ such that
    \[
        \norm{\mathbf{N}_k(x) - \mathbf{N}_k(\mathfrak{p}(x))P}_{1,2}\le 2\delta.
    \]
    By Proposition \ref{prop:norm_comp} it follows
    \[
        \norm{\mathbf{N}_k(x) - \mathbf{N}_k(\mathfrak{p}(x))P}\le 2k^{1/2}\delta
    \]
  and since the spectral norm is bounded by the Frobenius norm, Corollary 7.3.5 of \cite{Horn_Johnson_1985} gives
    \[
        \max_{i=1,\dots, n}\absval{\sigma_i(\mathbf{N}_k(x)) - \sigma_i(\mathbf{N}_k(\mathfrak{p}(x))P)} \le 2k^{1/2}\delta.
    \]
Permutation matrices are orthogonal, so right-multiplication does not change singular values.
Hence $\mathbf N_k(\mathfrak p(x))P$ has the same singular values as $\mathbf N_k(\mathfrak p(x))$ up to a permutation of the right singular vectors. From this we deduce that
    \[
        \max_{i=1,\dots, n}\absval{\sigma_i(\mathbf{N}_k(x)) - \sigma_i(\mathbf{N}_k(\mathfrak{p}(x)))} \le 2k^{1/2}\delta.
    \]
\end{proof}

\subsubsection{Stability of the Ellipsoid Bases}\label{sub:stability_bases}
The crucial ingredient in our stability proof for ellipsoids is the stability of the bases for the ellipsoid construction (Definition~\ref{def:ellipsoid_basis}). Recall  that $\mathbf{U}(x)$ is the matrix of left singular vectors of the neighborhood matrix.  For the stability result of the ellipsoid bases, we need the following theorem.
\begin{theorem}[Theorem~VII.5.9 {\cite{Bhatia96}}]\label{thm:bhatia_svd_stability}
    Let $S_1, S_2$ be two subsets of the positive half-line such that $\dist(S_1, S_2) = \eta > 0$. Let $A, B\in\R^{n\times n}$. Let $E$ and $E'$ be the orthogonal projections onto respectively the subspaces spanned by the right and the left singular vectors of $A$ corresponding to its singular values in $S_1$. Let $F$ and $F'$ be the projections associated to $B$ in the same way, corresponding to the singular values in $S_2$. Then
\[
\bigl(\|EF\|^2+\|E'F'\|^2\bigr)^{1/2}\le\frac{\sqrt{2}}{\eta}\,\|A-B\|.
\]
\end{theorem}

The next proposition tells us that, for $k$-generic sets, the basis spanned by the left singular vectors of $\mathbf{N}_k(x)$ and the one spanned by the left singular vectors of the perturbed $\mathbf{N}_k(\mathfrak{p}(x))$ are close.
\begin{proposition}[Stability of Ellipsoid Bases]\label{prop:stability_basismap}
    Let $k, n\in \mathbb{N}$ with $k \ge n$ and let $X\subset \RR^n$ be a $k$-generic subset of $\RR^n$, so that $\eta_X \coloneqq \min_{x\in X}\Bigl(\sigma_m(\mathbf{N}_k(x))-\sigma_{m+1}(\mathbf{N}_k(x))\Bigr) > 0.$
    We set
    \[
        \mathbf U_m(x)\coloneqq \bigl(\mathbf u_1(x)\ \cdots\ \mathbf u_m(x)\bigr),\qquad
        P(x)\coloneqq \mathbf U_m(x)\mathbf U_m(x)^T .
    \]
    Then there exists a $\delta_u(X) > 0$ such that for all $\delta \in [0,\delta_u(X))$ and any $\delta$-perturbation $\mathfrak{p}\colon X \to \RR^n$ the following estimate holds
  \[
        \norm{P(x)-P(\mathfrak{p}(x))} \le \frac{12\sqrt{2k}}{\eta_X}\,\delta\quad \forall x\in X.
    \]
\end{proposition}
\begin{proof}
   We claim that one can take
   \[
       \delta_u(X) = \min\set{\delta_N(X), \eta_X/(6k^{1/2}), \rho_X/4},
   \]
   where $\delta_N(X)$ is as in Lemma~\ref{lem:stability_neighborhoods} and $\rho_X:=\min_{x\neq y\in X}\|x-y\|$.

   Let $P$ be the permutation matrix from Corollary~\ref{cor:stability_distance_mats} and set $A\coloneqq \mathbf N_k(x)$ and $B\coloneqq \mathbf N_k(\mathfrak p(x))P$. Since $P$ is orthogonal, $B$ has the same singular values and left singular vectors as $\mathbf{N}_k(\mathfrak{p}(x))$. By Corollary~\ref{cor:stability_distance_mats} and Propositions~\ref{prop:norm_comp} and~\ref{prop:stability_axesmap},
   \[
       \norm{A-B}\le 2k^{1/2}\delta,
       \qquad
       \max_i |\sigma_i(A)-\sigma_i(B)|\le \norm{A-B}\le 2k^{1/2}\delta.
   \]
   If $\delta<\eta_X/(6k^{1/2})$, then $2k^{1/2}\delta\le \eta_X/3$. Define
   \[
       S_1\coloneqq [\sigma_m(A)-\eta_X/3,\infty),\qquad
       S_2\coloneqq [0,\sigma_{m+1}(A)+\eta_X/3].
   \]
   Then
   \[
       \dist(S_1,S_2)\ge (\sigma_m(A)-\sigma_{m+1}(A))-\tfrac{2}{3}\eta_X\ge \eta_X/3.
   \]
    The singular values $\sigma_1,\dots,\sigma_m$ of $A$ and of $B$ lie in $S_1$, and $\sigma_{m+1},\dots,\sigma_n$ lie in $S_2$. Theorem~\ref{thm:bhatia_svd_stability} is stated for square matrices, so we apply it to the $k\times k$ matrices obtained from $A$ and $B$ by appending $k-n$ zero rows. This changes neither $\|A-B\|$ nor the projections below. The added coordinates correspond to the singular value $0\in S_2$ and are orthogonal to the left singular subspaces for $S_1$. In the notation of Theorem~\ref{thm:bhatia_svd_stability}, let $E'$ be the orthogonal projection onto the left singular subspace of $A$ corresponding to $S_1$ and let $F'$ be the orthogonal projection onto the left singular subspace of $B$ corresponding to $S_2$. Concretely, $E'=P(x)$ and $F'=I-P(\mathfrak p(x))$. Applying Theorem~\ref{thm:bhatia_svd_stability} yields
   \[
       \|P(x)(I-P(\mathfrak p(x)))\| = \|E'F'\|
       \le \frac{\sqrt2}{\dist(S_1,S_2)}\|A-B\|
       \le \frac{\sqrt2}{\eta_X/3}\,2k^{1/2}\delta
       = \frac{6\sqrt{2k}}{\eta_X}\delta.
   \]
   Interchanging the roles of $A$ and $B$ (with the same $S_1$ and $S_2$) gives similarly
   \[
       \|(I-P(x))P(\mathfrak p(x))\|\le \frac{6\sqrt{2k}}{\eta_X}\delta.
   \]
   Using $P-Q=P(I-Q)-(I-P)Q$  we obtain
   \[
       \|P(x)-P(\mathfrak p(x))\|
       \le \|P(x)(I-P(\mathfrak p(x)))\|+\|(I-P(x))P(\mathfrak p(x))\|
       \le \frac{12\sqrt{2k}}{\eta_X}\delta.
   \]
\end{proof}

\subsubsection{Stability of Rips-Type Ellipsoid Barcodes}

The goal of this section is to extend the results developed in previous subsections to the level of persistence modules.
%The next corollary of Proposition~\ref{prop:perturb_ell_intersec}

The following theorem, which establishes an interleaving (as defined in Remark~\ref{rem:interleaving}) between the ellipsoid filtration of a point cloud $X$ and the perturbed point cloud $\mathfrak{p}(X)$ is the main step towards applying the theory from \cite{CdSO14, CDGO16} and establishing a stability result for Rips-type ellipsoid barcodes. 
\begin{theorem}[$\delta$-Interleaving of the Ellipsoid Filtrations]\label{cor:Ellipsoid_interleaving}
    Let $k, n\in \NN$ with $k \ge n$ and let $X\subset \R^n$ be a $k$-generic subset. Let $\varepsilon_{\max}(X):=q\operatorname{diam}(X)$. Then for all $\delta$-perturbations $\mathfrak{p}\colon X\to \RR^n$ with $\delta\in [0, \delta_u(X))$, the ellipsoid filtrations $E^q(X)$ and $E^q(\mathfrak{p}(X))$, as defined in Notation~\ref{def:rips-type}, are $(C\eps_{\max}(X)\delta)$-interleaved, where $C$ is the constant from Lemma~\ref{prop:perturb_ell_intersec}, and $\delta_u(X)$ is as in Proposition~\ref{prop:stability_basismap}.
\end{theorem}
\begin{remark}
The scaling of the ellipsoids is multiplicative and not additive as was the case with Vietoris-Rips complexes in~\cite{CdSO14}.
\end{remark}

To prove Theorem~\ref{cor:Ellipsoid_interleaving} we rely on the statement that if ellipsoids built on a point cloud intersect, then the ellipsoids constructed on the perturbed data
also intersect if the perturbation is sufficiently small. This we formalize with the following lemma.

\begin{lemma}[Ellipsoids on Perturbed Data Intersect]\label{prop:perturb_ell_intersec}
    Let $k, n \in \mathbb{N}$ with $k\ge n$ and let $X\subset \RR^n$ be a $k$-generic set. Recall that $\eta_X$ denotes the smallest $m$-gap $\eta_X=\min_{x\in X}\bigl(\sigma_m(\mathbf N_k(x))-\sigma_{m+1}(\mathbf N_k(x))\bigr),$
    $m$ is the assumed dimension of the manifold, and $q$ is the desired axis ratio for the ellipsoids. There exists a constant $C = C(X, \eta_X, k, m, q)$ such that if $x_1, x_2 \in X$ are so that 
    \[
        E^q_\eps(x_1)\cap E^q_\eps(x_2)\neq \emptyset,
    \]
    then
    \[
        E^q_{(1 + C\delta)\eps}(\mathfrak{p}(x_1)) \cap E^q_{(1 + C\delta)\eps}(\mathfrak{p}(x_2))\neq \emptyset
    \]
    for all $\delta$-perturbations $\mathfrak{p}$ with $\delta \in [0, \delta_u(X))$, where $\delta_u(X)$ is as in Proposition~\ref{prop:stability_basismap}.
\end{lemma}

\begin{proof}
Fix a $\delta$-perturbation
$\mathfrak{p}\colon X\to\mathbb{R}^n$, with
$\delta\in[0,\delta_u(X))$.

For $l\in\{1,2\}$, let $P^{(l)}$ and $\widetilde P^{(l)}$
denote the orthogonal projections onto the $m$-dimensional
subspaces spanned by the first $m$ left singular vectors at $x_l$
and $\mathfrak{p}(x_l)$, respectively; explicitly,
\[
P^{(l)}
=
\sum_{j=1}^m
\mathbf{u}_j^{(l)}(\mathbf{u}_j^{(l)})^T,
\qquad
\widetilde P^{(l)}
=
\sum_{j=1}^m
\widetilde{\mathbf{u}}_j^{(l)}
(\widetilde{\mathbf{u}}_j^{(l)})^T.
\]
We also set
\[
Q^{(l)}:=I-P^{(l)},
\qquad
\widetilde Q^{(l)}:=I-\widetilde P^{(l)}.
\]
Thus $P^{(l)}$ and $Q^{(l)}$ are the projections onto the
tangent and normal directions of the ellipsoid centered at $x_l$,
while $\widetilde P^{(l)}$ and $\widetilde Q^{(l)}$ are the
corresponding projections for the ellipsoid centered at
$\mathfrak{p}(x_l)$.

By Proposition~\ref{prop:stability_basismap},
\[
\|P^{(l)}-\widetilde P^{(l)}\|
=
\|Q^{(l)}-\widetilde Q^{(l)}\|
\le A\delta,
\qquad
A:=\frac{12\sqrt{2k}}{\eta_X}.
\]

If $x_1=x_2$, the claim follows immediately. Therefore we may assume
$x_1\neq x_2$. Since $X$ is finite, let
\[
\rho_X
:=
\min_{\substack{x,y\in X\\x\neq y}}\|x-y\|>0.
\]

Choose
\[
z\in E_\varepsilon^q(x_1)\cap E_\varepsilon^q(x_2)
\]
and set
\[
v_l:=z-x_l,
\qquad
d_l:=\mathfrak{p}(x_l)-x_l
\]
for $l\in\{1,2\}$.

Then
\[
z-\mathfrak{p}(x_l)=v_l-d_l,
\qquad
\|d_l\|\le\delta.
\]

Set $b:=\varepsilon/q$. Since $z\in E_\varepsilon^q(x_l)$,
\[
\frac{\|P^{(l)}v_l\|^2}{\varepsilon^2}
+
\frac{\|Q^{(l)}v_l\|^2}{b^2}
\le 1.
\]
In particular,
\[
\|P^{(l)}v_l\|\le\varepsilon,
\qquad
\|Q^{(l)}v_l\|\le b,
\]
and hence
\[
\|v_l\|
\le \sqrt{\varepsilon^2+b^2}
=
(\sqrt{1+q^{-2}})\varepsilon.
\]
Therefore
\[
\rho_X
\le \|x_1-x_2\|
\le \|z-x_1\|+\|z-x_2\|
\le 2(\sqrt{1+q^{-2}})\varepsilon,
\]
so that dividing by $\rho_X$ and multiplying by $\delta$ we get
\[
\delta
\le
\frac{2(\sqrt{1+q^{-2}})}{\rho_X}\varepsilon\delta.
\tag{*}
\]

We now compare the tangent and normal components of
$z-\mathfrak{p}(x_l)$ with those of $z-x_l$. Using the triangle
inequality, Proposition~\ref{prop:stability_basismap}, and $(*)$, we obtain
\begin{align*}
\|\widetilde P^{(l)}(z-\mathfrak{p}(x_l))\|
&=
\|\widetilde P^{(l)}(v_l-d_l)\|\\
&\le
\|P^{(l)}v_l\|
+
\|(\widetilde P^{(l)}-P^{(l)})v_l\|
+
\|d_l\|\\
&\le
\|P^{(l)}v_l\|
+
A\delta\|v_l\|+\delta\\
&\le
\|P^{(l)}v_l\|+D\varepsilon\delta,
\end{align*}
where
\[
D:=\sqrt{1+q^{-2}}\,\Bigl(A+\frac{4}{\rho_X}\Bigr).
\]
Exactly the same argument gives
\[
\|\widetilde Q^{(l)}(z-\mathfrak{p}(x_l))\|
\le
\|Q^{(l)}v_l\|+D\varepsilon\delta.
\]

Consequently,
\begin{align*}
&
\frac{
\|\widetilde P^{(l)}(z-\mathfrak{p}(x_l))\|^2
}{\varepsilon^2}
+
\frac{
\|\widetilde Q^{(l)}(z-\mathfrak{p}(x_l))\|^2
}{b^2}
\\
&\qquad\le
\frac{(\|P^{(l)}v_l\|+D\varepsilon\delta)^2}
{\varepsilon^2}
+
\frac{(\|Q^{(l)}v_l\|+D\varepsilon\delta)^2}
{b^2}
\\
&\qquad=
\frac{\|P^{(l)}v_l\|^2}{\varepsilon^2}
+
\frac{\|Q^{(l)}v_l\|^2}{b^2}
+
2D\varepsilon\delta
\left(
\frac{\|P^{(l)}v_l\|}{\varepsilon^2}
+
\frac{\|Q^{(l)}v_l\|}{b^2}
\right)
\\
&\hspace{4cm}
+
D^2\varepsilon^2\delta^2
\left(
\frac{1}{\varepsilon^2}
+
\frac{1}{b^2}
\right)
\\
&\qquad\le
1+
2D\varepsilon\delta
\left(
\frac{1}{\varepsilon}
+
\frac{1}{b}
\right)
+
D^2\varepsilon^2\delta^2
\left(
\frac{1}{\varepsilon^2}
+
\frac{1}{b^2}
\right)
\\
&\qquad=
1+
2D(1+q)\delta
+
D^2(1+q^2)\delta^2.
\end{align*}

Here we used
\[
\frac{\|P^{(l)}v_l\|^2}{\varepsilon^2}
+
\frac{\|Q^{(l)}v_l\|^2}{b^2}
\le 1,
\qquad
\|P^{(l)}v_l\|\le\varepsilon,
\qquad
\|Q^{(l)}v_l\|\le b,
\]
as well as $b=\varepsilon/q$.

Since $\delta<\delta_u(X)$, we have
$\delta^2\le \delta_u(X)\delta$. Hence
\begin{align*}
1+2D(1+q)\delta+D^2(1+q^2)\delta^2
&\le
1+
\left(
2D(1+q)
+
D^2(1+q^2)\delta_u(X)
\right)\delta.
\end{align*}
Thus, setting
\[
C
:=
2D(1+q)
+
D^2(1+q^2)\delta_u(X),
\]
we obtain
\[
\frac{
\|\widetilde P^{(l)}(z-\mathfrak{p}(x_l))\|^2
}{\varepsilon^2}
+
\frac{
\|\widetilde Q^{(l)}(z-\mathfrak{p}(x_l))\|^2
}{b^2}
\le 1+C\delta. \tag{$**$}
\]

Now set
\[
\varepsilon':=(1+C\delta)\varepsilon.
\]
Since $\varepsilon'/q=(1+C\delta)b$, it follows from $(**)$ that
\begin{align*}
&
\frac{
\|\widetilde P^{(l)}(z-\mathfrak{p}(x_l))\|^2
}{(\varepsilon')^2}
+
\frac{
\|\widetilde Q^{(l)}(z-\mathfrak{p}(x_l))\|^2
}{(\varepsilon'/q)^2}
\\
&\qquad\le
\frac{1+C\delta}{(1+C\delta)^2}
\le 1.
\end{align*}
Thus
\[
z\in E_{\varepsilon'}^q(\mathfrak{p}(x_l)),
\qquad l=1,2,
\]
and therefore
\[
E_{(1+C\delta)\varepsilon}^q(\mathfrak{p}(x_1))
\cap
E_{(1+C\delta)\varepsilon}^q(\mathfrak{p}(x_2))
\neq\emptyset.
\]
This proves the claim.
\end{proof}

\begin{proof}[Proof of Theorem~\ref{cor:Ellipsoid_interleaving}]
Let $\mathfrak{p}$ be a $\delta$-perturbation, with $\delta\in [0, \delta_u(X))$, and let $\sigma$ be a simplex in $E^q_\eps(X)$. By definition for all $x_1, x_2\in \sigma$ it holds
\[
E^q_\eps(x_1) \cap E^q_\eps(x_2)\neq \emptyset.
\]
By Lemma \ref{prop:perturb_ell_intersec}
\[
    E^q_{(1+C\delta)\eps}(\mathfrak{p}(x_1))\cap E^q_{(1+C\delta)\eps}(\mathfrak{p}(x_2))\neq \emptyset.
\]
Moreover, we have
\[
    (1 + C\delta)\eps \le \eps + C\eps_{\max}(X)\delta
\]
and hence
\[
    E^q_{\eps+C\eps_{\max}(X)\delta}(\mathfrak{p}(x_1))\cap E^q_{\eps+C\eps_{\max}(X)\delta}(\mathfrak{p}(x_2))\neq \emptyset.
\]
This implies that $\mathfrak p(\sigma)\in E^q_{\varepsilon+C\varepsilon_{\max}(X)\delta}\big(\mathfrak p(X)\big)$. This proves the inclusion for $\varepsilon\le\varepsilon_{\max}(X)$, where we used $(1+C\delta)\varepsilon\le\varepsilon+C\varepsilon_{\max}(X)\delta$. 

For $\varepsilon'\ge\varepsilon_{\max}(X)=q\operatorname{diam}(X)$, both $E^q_{\varepsilon'}(X)$ and $E^q_{\varepsilon'}(\mathfrak{p}(X))$ are full simplices: by Proposition~\ref{rips} they contain the Rips complexes at scale $\varepsilon'/q\ge\operatorname{diam}(X)$, and $\operatorname{diam}(\mathfrak{p}(X))\le\operatorname{diam}(X)+2\delta\le 2\operatorname{diam}(X)$ since $\delta<\delta_N(X)\le\operatorname{diam}(X)/4$ (take $x'=x$ in the definition of $f_x$). Hence the inclusion also holds for $\varepsilon>\varepsilon_{\max}(X)$.

Due to symmetry, we get the same result for the inverse of $\mathfrak{p}$. Hence, we have proved the claim.
\end{proof}

The next statement is the main result of this section, proving that Rips-type ellipsoid barcodes are stable under small $\delta$-perturbations.
\begin{theorem}\label{thm:ellipsoid_barcode_stability}
Let $k, n\in \mathbb{N}$ with $k \ge n$, $X\subset\RR^n$ be a $k$-generic subset, and $\mathfrak{p}\colon X\to \RR^n$ be a $\delta$-perturbation with $\delta \in [0, \delta_u(X))$. Then
\[
    \distance_b\!\left(\dgm(E^q(X)),\,\dgm(E^q(\mathfrak{p}(X)))\right) \le Cq\diam(X)\delta,
\]
where $C$ is the constant from Lemma~\ref{prop:perturb_ell_intersec}, and $\delta_u(X)$ is as in Proposition~\ref{prop:stability_basismap}.
\end{theorem}
\begin{proof}
   Note that a persistence module arising from the Rips-type ellipsoid filtration of a finite subset of $\RR^n$ is always pointwise finite-dimensional, just as in the case of Vietoris–Rips and \v{C}ech filtrations. Since $X$ is finite, $\diam(X) < \infty$. So, we may set
    \[
        \eps_{\max}(X) = q\diam(X)
    \]
    and apply Theorem~\ref{cor:Ellipsoid_interleaving} and Theorem \ref{thm:Chazal} to conclude the proof.
\end{proof}

%% file: sections/experiments.tex
\section{Experiments}
\label{sec:Experiments}

In this section we present a series of experiments whose primary goal is to highlight the differences between Rips complexes and Rips-type ellipsoid complexes.
Our experiments aim to answer when ellipsoid barcodes are more \emph{expressive} than Rips barcodes, i.e., in which situations an ellipsoid barcode uncovers \emph{more} information about a dataset than a Rips barcode and to demonstrate that using ellipsoids one can get valuable information from smaller samples.
To this end we provide a visual analysis of both types of barcodes on synthetic and real-world datasets (conformation space of cyclo-octane), followed by several \emph{classification experiments}. 
For the latter, we draw on previous work~\cite{Turkes22a} to obtain a setting in which the performance of Rips barcodes~(and derived topological descriptors) is already well-studied. All barcodes in this section are computed with GUDHI~\cite{gudhi:FilteredComplexes} over the field $\mathbb{Z}/11$.
%For the latter, we use the standard example of the cyclo-octane dataset.

\subsection{Dog Bone Example}\label{subsn:dogbone}

\begin{figure}[ht]
    \centering
    \includegraphics[width=0.3\linewidth]{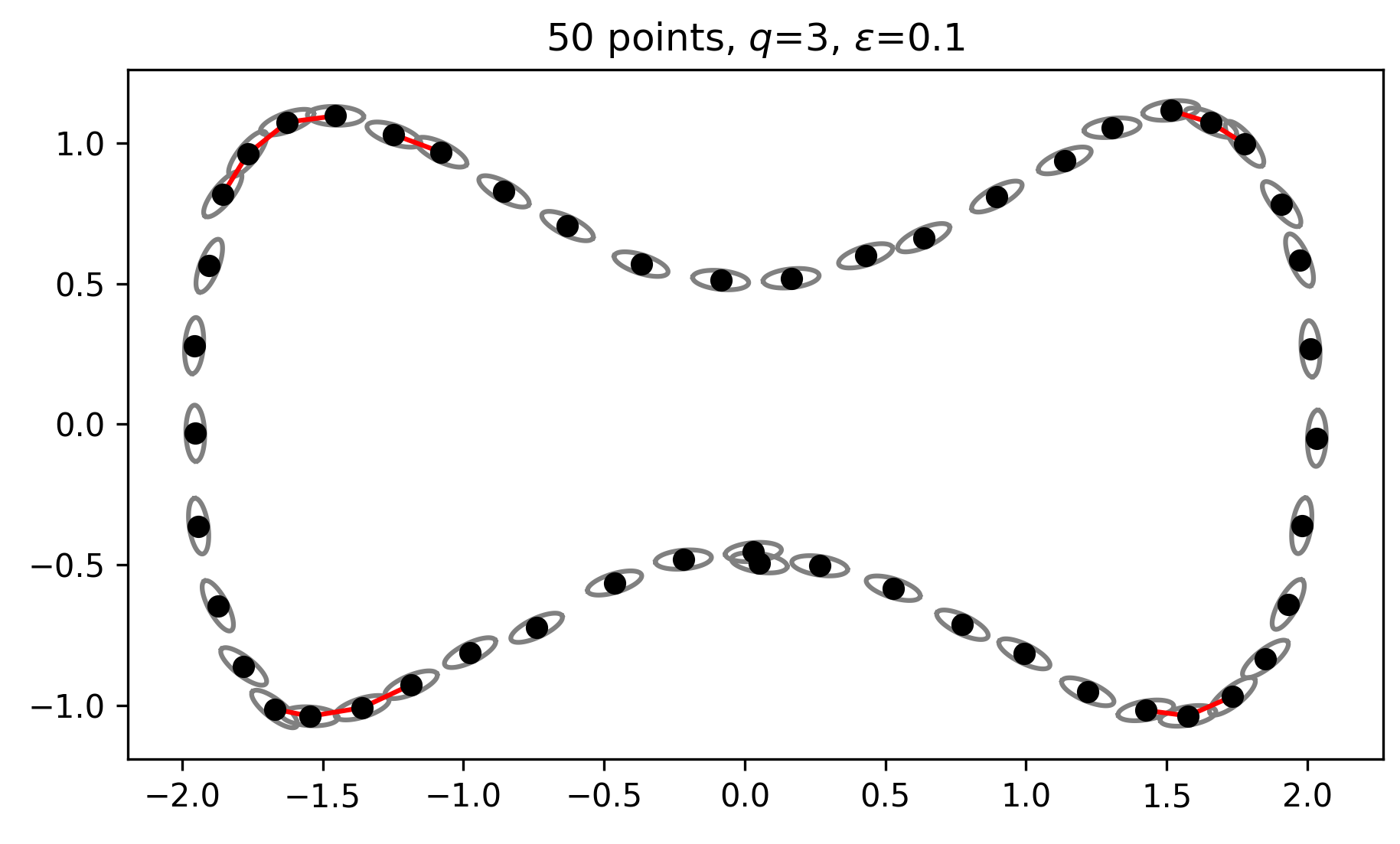}
     \includegraphics[width=0.3\linewidth]{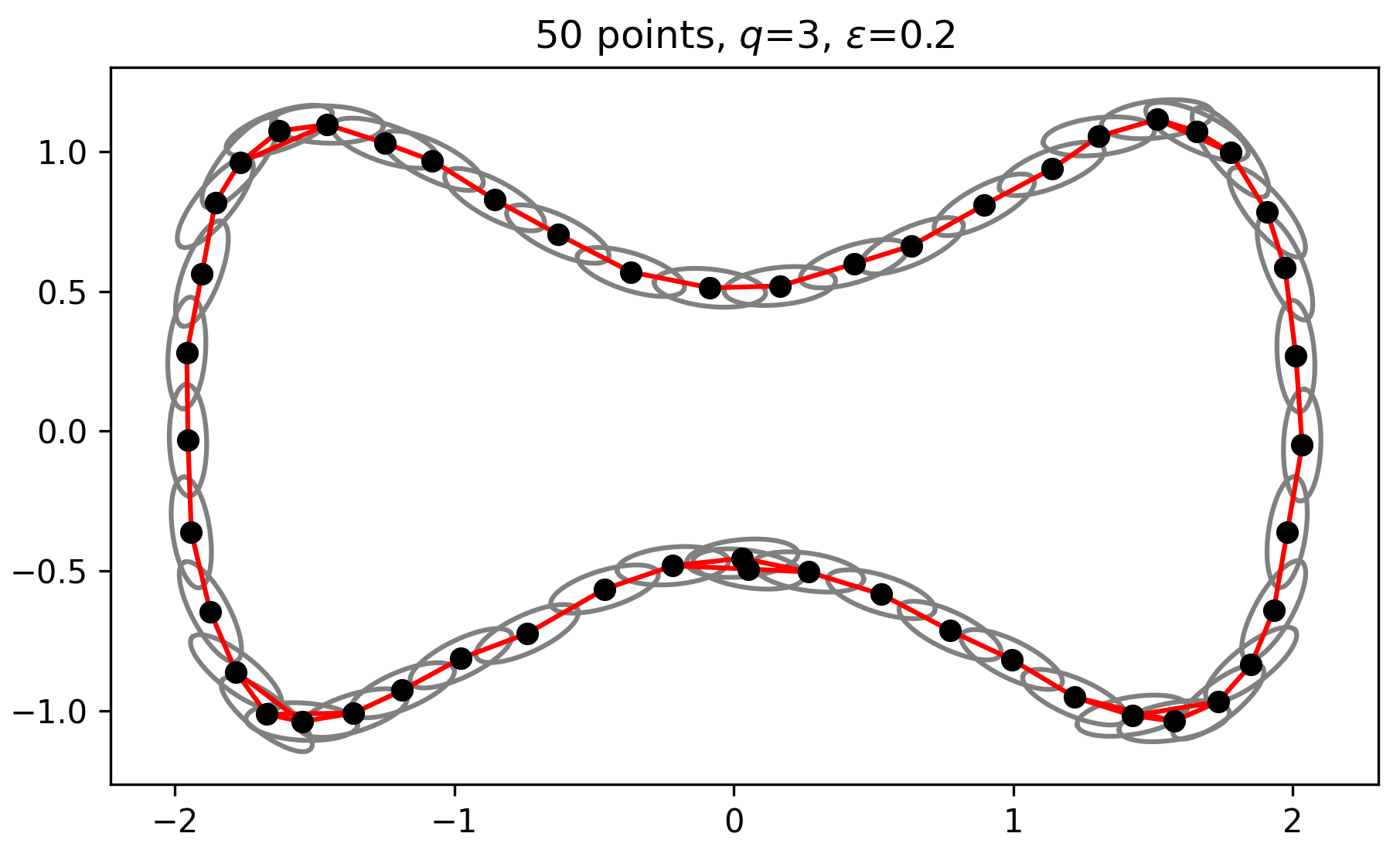}
     \includegraphics[width=0.3\linewidth]{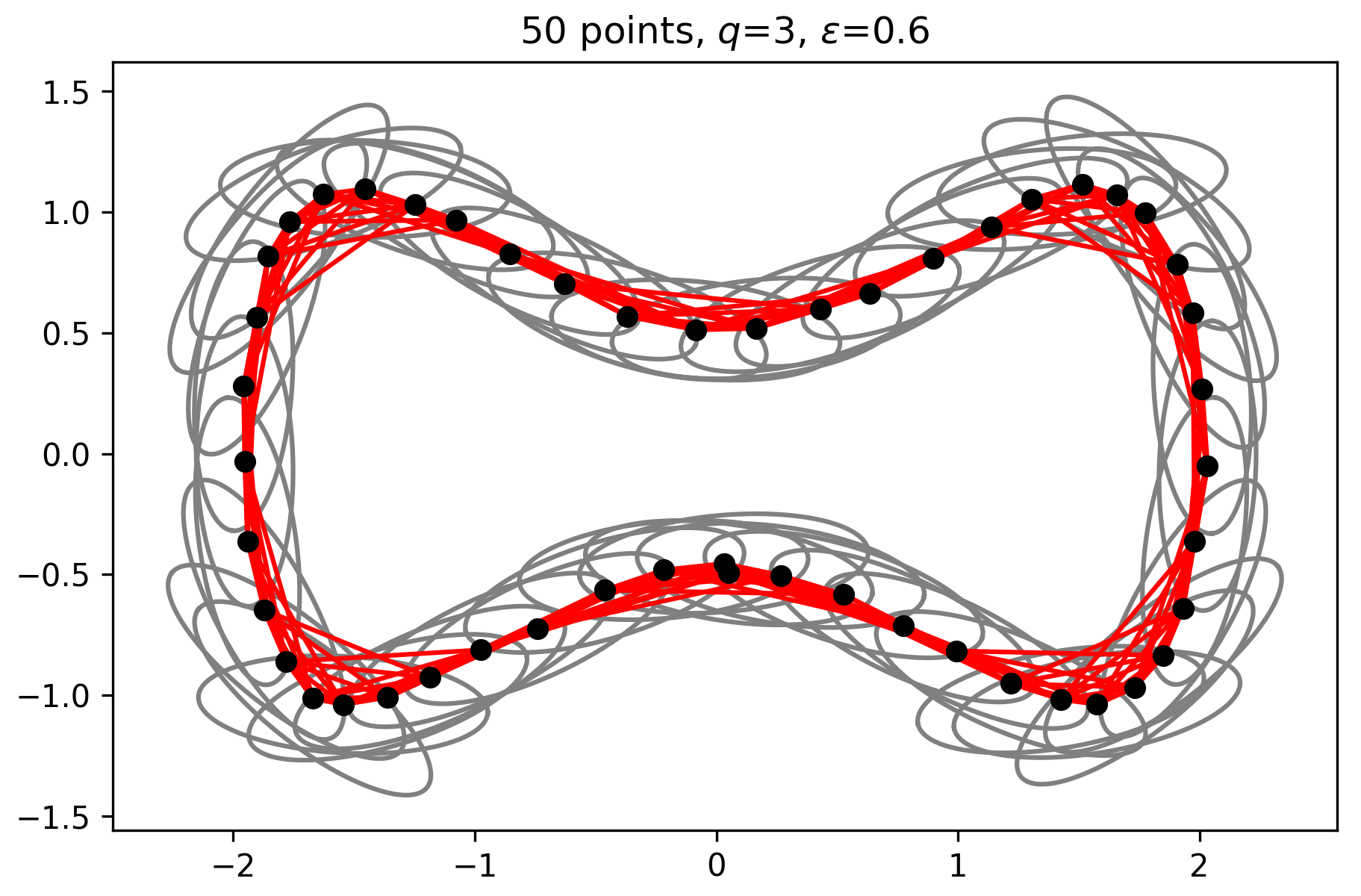}
    \caption{%
        Ellipsoid complexes for $q=3$ at scales $\varepsilon=0.1, \varepsilon=0.2$ and $\varepsilon=0.6$ for a point cloud sampled from a curve resembling a dog bone.
    }
    \label{fig:ellipsoids-dog-bone}
\end{figure}

Examples where ellipsoids are advantageous compared to Rips complexes include spaces with bottlenecks (as already remarked in~\cite{Breiding2018} for ellipsoid-driven complexes). For example, consider a curve in the shape of a dog bone. Figure~\ref{fig:ellipsoids-dog-bone} represents ellipsoids for $q=3$ at different scales: $\varepsilon=0.1, \varepsilon=0.2$ and $\varepsilon=0.6$.

\begin{figure}[ht]
    \centering
    \includegraphics[width=0.65\linewidth]{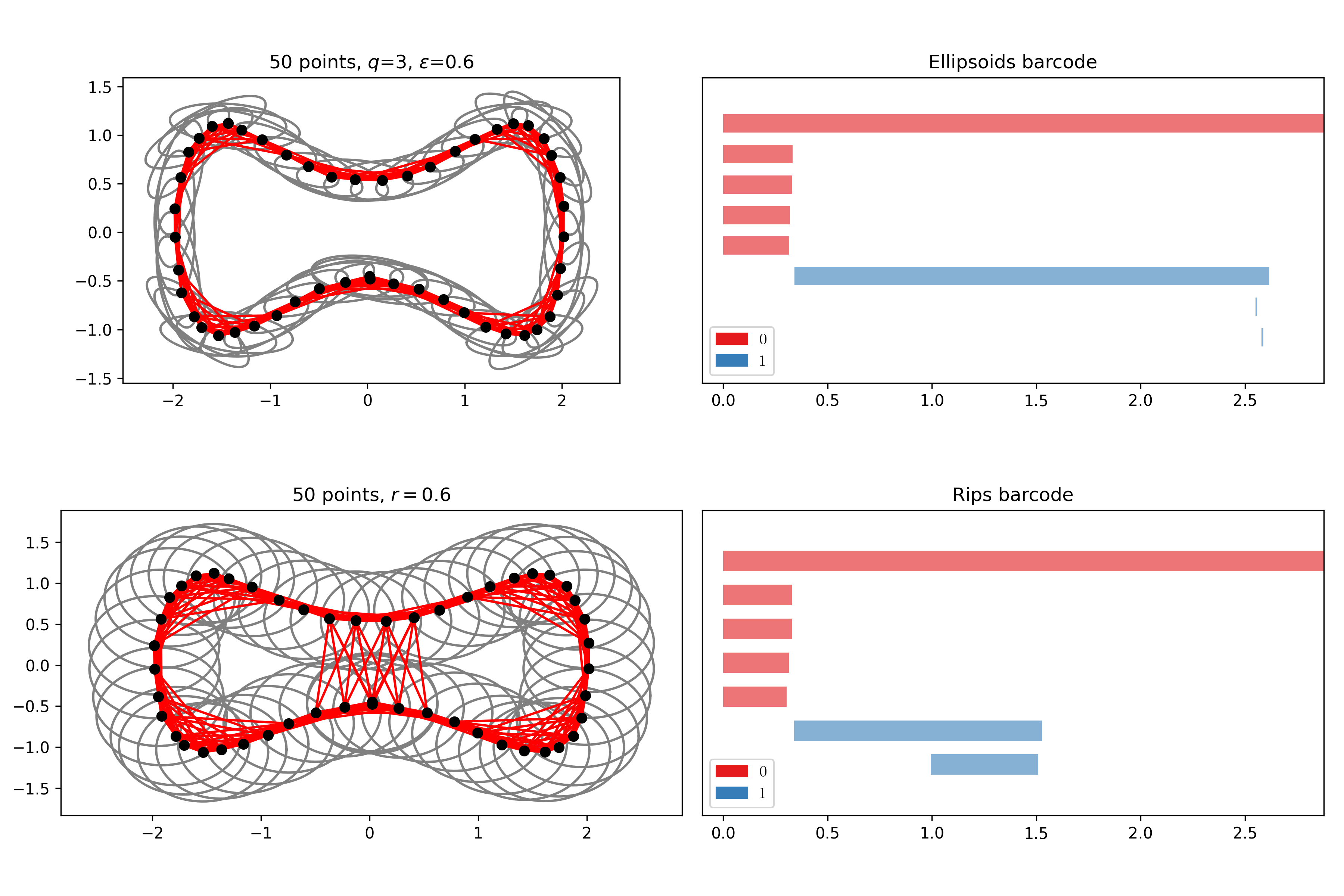}
    \caption{%
        Top left: Ellipsoids for $q=3$ at scale $\varepsilon=0.6$. Bottom left: Rips complex at $\varepsilon=0.6$. Top right: Rips-type ellipsoid barcode. Bottom right: Rips barcode. 
    }
    \label{fig:ellipsoid-vs-rips}
\end{figure}

In cases like this dataset the balls around points on the bottleneck may intersect for $\varepsilon$ smaller than that which is necessary for the full cycle to appear. This is demonstrated in Figure~\ref{fig:ellipsoid-vs-rips}. The ellipsoid barcode has one long bar in 1-dimensional homology, whereas the Rips barcode shows two prominent features.

%%%%%%%%%%%%%%%%%%%%%%%%%%%%%%%%%%%%%%%%%%%%%%%%%%%%%%%%%%%%%%%%%%%%%%%%
\subsection{Point Cloud Classification}\label{subsec:point-cloud-classification}
%%%%%%%%%%%%%%%%%%%%%%%%%%%%%%%%%%%%%%%%%%%%%%%%%%%%%%%%%%%%%%%%%%%%%%%%

%The main feature of the ellipsoid complex is that the ellipsoids used in the construction contain information about the spread of the neighbouring points. In particular, this means that if the point cloud data is sampled from a manifold, the ellipsoids will approximate the tangent spaces of the manifold.

%The experiments show that since ellipsoids approximate the underlying manifold structure of data, their barcodes lead to better classification results in sparsely sampled point clouds.
%genus classification results

To test how the classification based on the Rips-type ellipsoid complex compares to other methods, we run experiments analogous to ones described in \textcite{Turkes22a}. We generate point clouds of 35 different shapes in $\RR^2$ and $\RR^3$ with each shape sampled from a manifold containing either 0, 1, 2, 4 or 9 holes (see Table~\ref{tab:transformations} on the left for examples). We use the 20 shapes introduced in \textcite{Turkes22a}, as well as 15 additional shapes. The latter are meant to highlight the settings in which ellipsoids perform particularly well.

\begin{table}[ht]
\centering
\begin{minipage}[t]{0.42\textwidth}
  \centering
  \setlength{\tabcolsep}{2pt}
  \newcommand{\pc}[1]{\raisebox{\dimexpr.5ex-.5\height\relax}{\includegraphics[width=\dimexpr(\linewidth-1.5cm)/6\relax]{figures/#1}}}%
  \begin{tabular}[t]{@{}cc@{}}
    \toprule
    Holes & Example point clouds \\
    \midrule
    0 & \pc{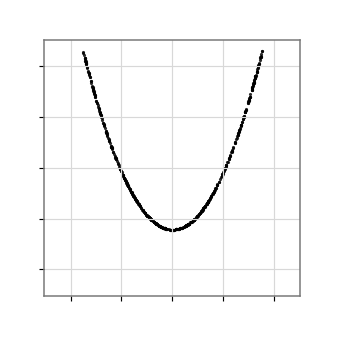}\pc{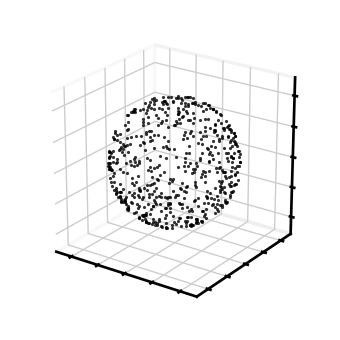}\pc{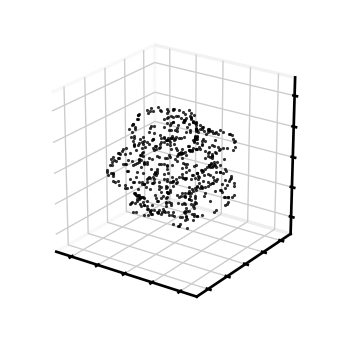}\pc{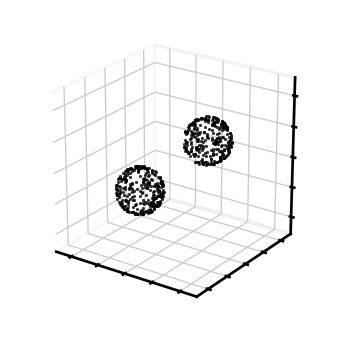}\pc{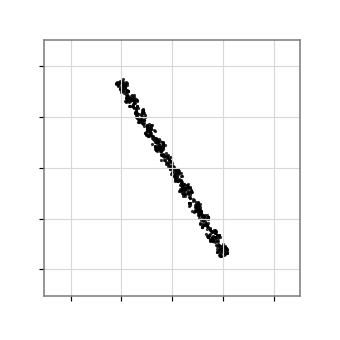}\pc{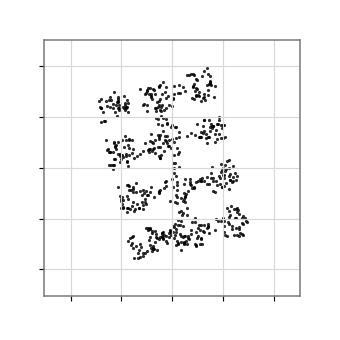} \\
    1 & \pc{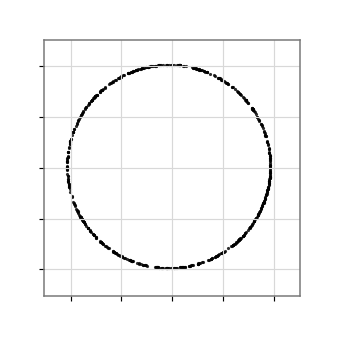}\pc{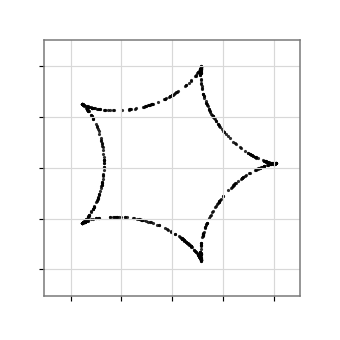}\pc{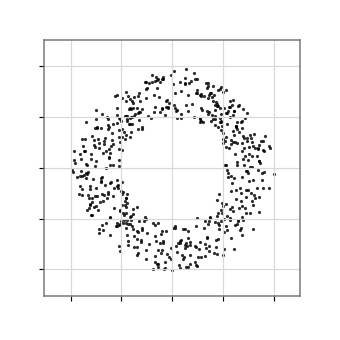}\pc{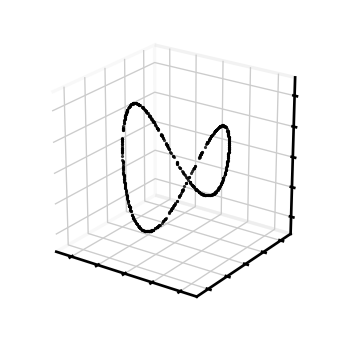}\pc{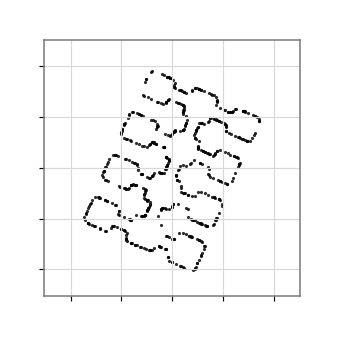} \\
    2 & \pc{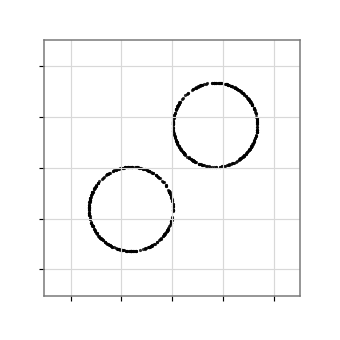}\pc{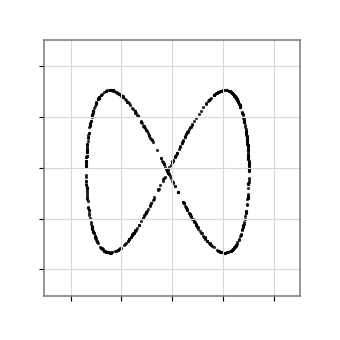}\pc{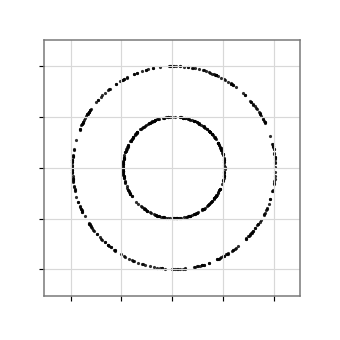}\pc{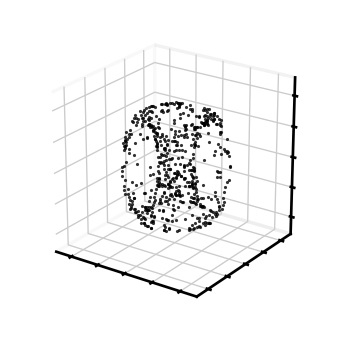}\pc{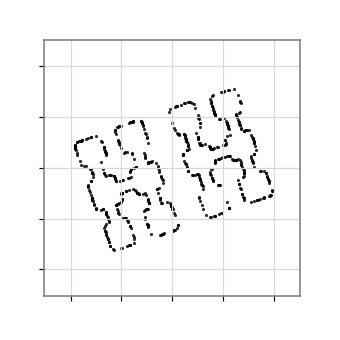}\pc{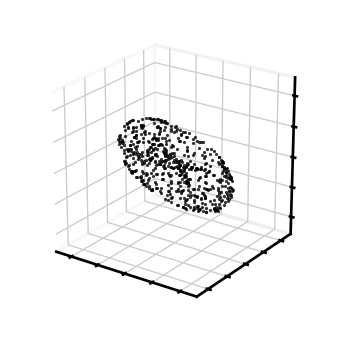} \\
    4 & \pc{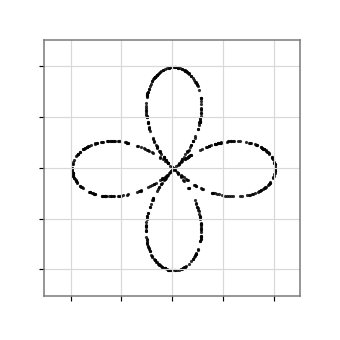}\pc{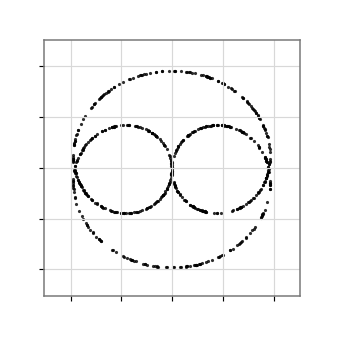}\pc{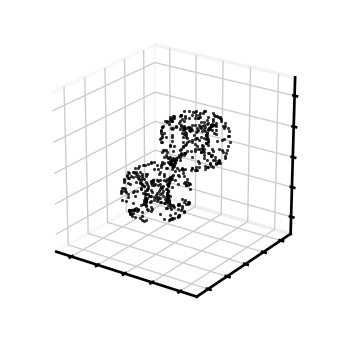}\pc{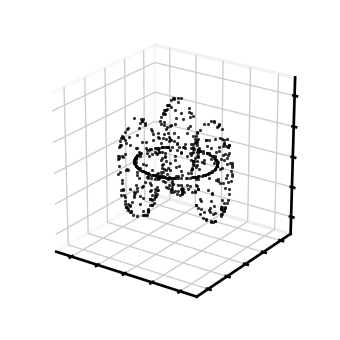}\pc{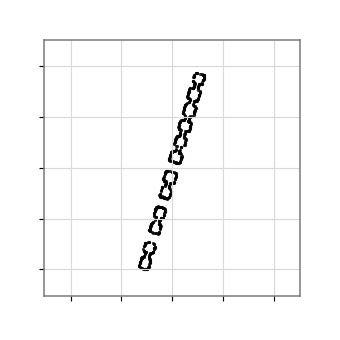}\pc{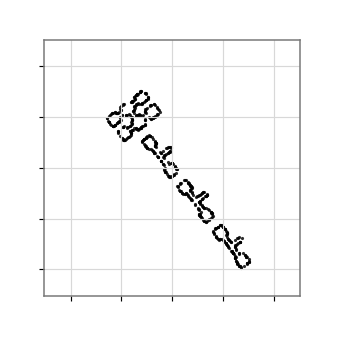} \\
    9 & \pc{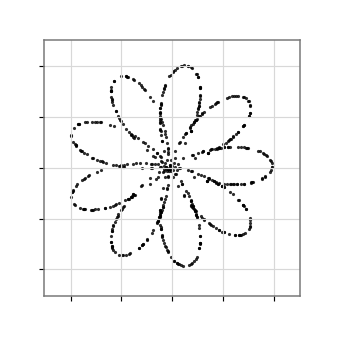}\pc{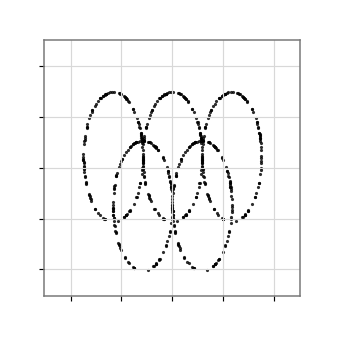}\pc{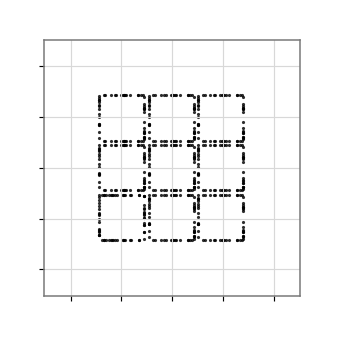}\pc{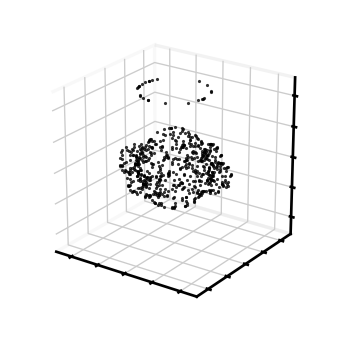} \\
    \bottomrule
  \end{tabular}
\end{minipage}\hfill
\begin{minipage}[t]{0.55\textwidth}
  \footnotesize
  \begin{tabularx}{\linewidth}[t]{@{}l>{\raggedright\arraybackslash}X@{}}
    \toprule
    \textbf{Transformation} & \textbf{Explanation} \\
    \midrule
    original & The original dataset. \\
    \hline
    translation & Translation by random numbers chosen from $[-1, 1]$ for each direction. \\
    \hline
    rotation & Clockwise rotation by an angle chosen uniformly from $[-20, 20]$ degrees. \\
    \hline
    stretch & Scale by a factor chosen uniformly from $[0.8, 1.2]$ in the $x$-direction leaving the other coordinates unchanged. \\
    \hline
    shear & Shear by a factor chosen uniformly from $[-0.2, 0.2]$. A shearing factor of 1 means that a horizontal line turns into a line at 45 degrees. \\
    \hline
    Gaussian noise & Random noise drawn from normal distribution $\mathcal N (0, \sigma^2)$ with the standard deviation $\sigma$ uniformly chosen from $[0, 0.1]$ is added to the point cloud. \\
    \hline
    outliers & A fraction, chosen uniformly from $[0, 0.1]$, of point cloud points are replaced with points sampled from a uniform distribution within the range of the point cloud. \\
    \bottomrule
  \end{tabularx}
\end{minipage}
\caption{Left: Example point clouds of the expanded `holes' dataset~\cite{Turkes22a}. Our experiments assess to what extent predictions of the number of holes also work with fewer points. Right: Explanations of the data transformations used on the datasets.}
\label{tab:transformations}
\end{table}

For each shape, we generate 5 different point clouds, each consisting of 600 points. Note that in \textcite{Turkes22a}, 1000 points were used. Due to the property of ellipsoids to approximate the underlying manifold structure of point clouds, we expect the classification accuracy to remain high even with the lower-resolution datasets. We have therefore decided to reduce sampling to 600 points per point cloud.

\begin{figure}[h]
    \centering
    \includegraphics[width=0.8\textwidth]{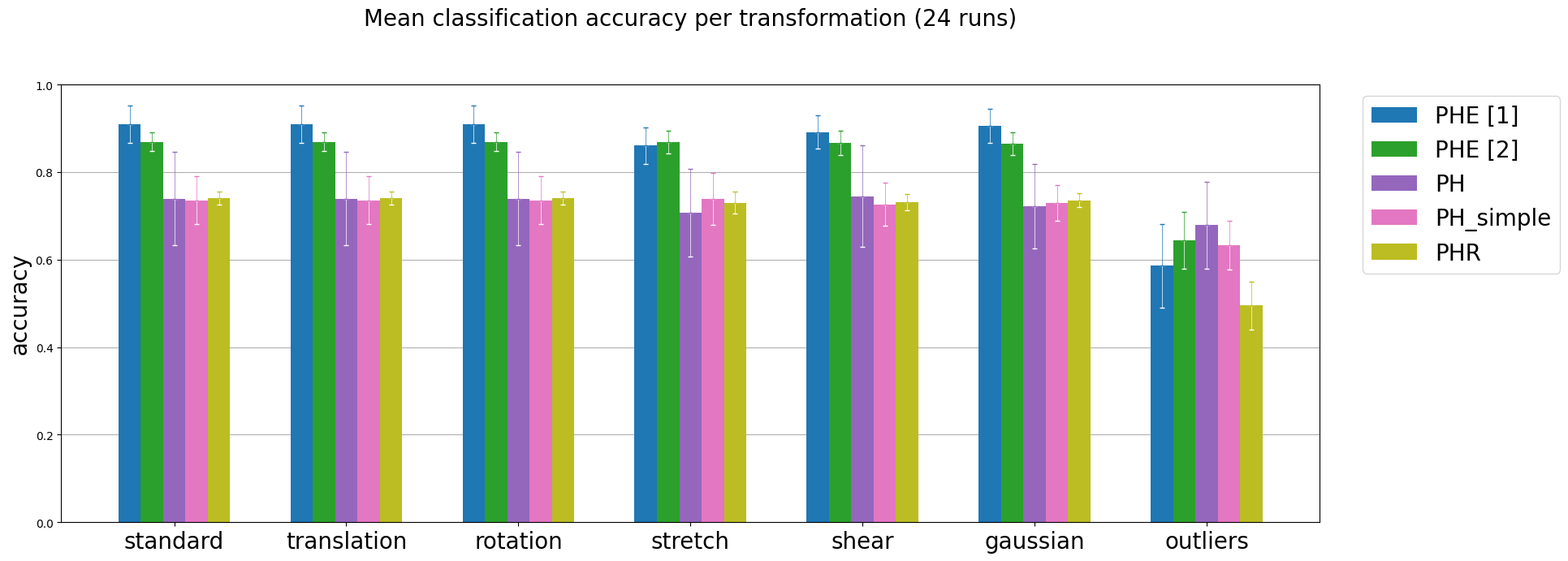}
    \caption{Classification accuracies across 24 runs on 175 point clouds, each consisting of 600 points. PHE refers to pipelines using ellipsoid complexes. More details, as well as an explanation of other pipelines can be found in Section~\ref{subsec:point-cloud-classification}. 
    % The parameter values for PHE pipelines are shown in Table~\ref{table:params}.
    }
    \label{fig:point_cloud_classification}
\end{figure}

The experiments consist of classifying point clouds in $\mathbb R^2$ and $\mathbb R^3$ via different methods:
\begin{enumerate}
    \item Using barcodes coming from the Rips-type ellipsoid complexes. 
    %Borrowing from the language used in \cite{Turkes22a}, 
    We refer to these pipelines as PHE.
    In the results shown in Figure~\ref{fig:point_cloud_classification}, PHE [1] refers to the complexes with axes ratios of $2.5:1$ and where the orientation of the ellipsoid at any given point was determined by performing PCA on 20 neighboring points. These parameter values were chosen as the best ones after a parameter space search was performed on the training data.
    We also include as PHE [2] the results from the complex calculated using PCA on 20 neighboring points, only with axes ratios calculated from the eigenvalues of the PCA. The benefit of the latter one being that there is no need for the user to specify the axes ratios, as they are calculated from the results of the PCA.
    
    \item Using barcodes coming from the standard Vietoris--Rips complexes. 
    We refer to these pipelines as PHR.
    
    \item Using barcodes coming from alpha complexes generated using Distance-to-Measure as the filtration function (PH). As noted in \textcite[Section~2.2]{Turkes22a}, the filtration function used in Rips complex is sensitive to outliers, and to mitigate this limitation, the so-called Distance-to-Measure function is used instead. This function measures the average distance from a number of neighbors on the point cloud. In \cite{Turkes22a}, this method was shown to outperform other methods considered. We therefore use it as the main benchmark against which we compare ellipsoid-based methods and we do not repeat other experiments. We also include a simplified version of this method (PH simple) introduced in the paper. 
    % \item Using only the 10 longest lifespans in the barcodes coming from alpha complexes generated using Distance-to-Measure as the filtration function (PH simple).
    % \item Support vector machine trained on the distance matrices of point clouds (ML).
    % \item Fully connected neural network with a single hidden layer (NN shallow).
    % \item Fully connected neural network with multiple layers (NN deep).
    % \item PointNet~\cite{pointnet1}.
\end{enumerate}

To perform classification based on ellipsoid data, we calculate the barcodes corresponding to the Rips-type ellipsoid complex and then use the remainder of the PH pipeline developed in \textcite{Turkes22a}. In particular, we feed a support vector machine with a signature calculated from the ellipsoid barcode. We choose this signature from the following:
\begin{enumerate}[(a)]
    \item signature containing 100 longest lifespans (or 10 in the case of PH simple); \label{item:signature1}
    \item persistence images (generated by choosing various different parameters) \cite{adams2017persistence}; \label{item:signature2}
    \item persistence landscapes (generated by choosing various different parameters) \cite{bubenik2015statistical}. \label{item:signature3}
\end{enumerate}
Whichever option between (\ref{item:signature1}), (\ref{item:signature2}) or (\ref{item:signature3}) (with whichever combination of parameters) leads to the highest score, i.e., accuracy, is then used as a signature in the actual classification.
This means that, depending on the datasets, different signatures might be used on the ellipsoid barcodes.

In the classification based on the Rips complex we perform the same steps, except that we use the barcodes coming from the Rips complex.

The experiments test the classification of the original datasets, as well as of the datasets after various transformations have been applied to them: translation, rotation, stretching, shear mapping, adding Gaussian noise, and replacing a certain number of points with outliers. In Table~\ref{tab:transformations} we reproduce the table from \textcite{Turkes22a}, explaining these transformations in more detail.

%\begin{table}[h]
 %   \centering
  %  \begin{tabular}{l p{9cm}}
   %     \textbf{Transformation name} & \textbf{Explanation} \\
    %    \hline 

     %   original & The original dataset. \\
      %  \hline
        
       % translation & Translation by random numbers chosen from $[-1, 1]$ for each direction. \\
       % \hline
        
        %rotation & Clockwise rotation by an angle chosen uniformly from $[-20, 20]$ degrees clockwise. \\
        %\hline
        
%        stretch & Scale by a factor chosen uniformly from $[0.8, 1.2]$ in the $x$-direction. The other coordinates remain unchanged, so that the point cloud is stretched. Stretching factor of $0.8$ results in shrinking the point cloud by $20\%$, and the factor of $1.2$ makes it $20\%$ larger. \\
 %       \hline
        
  %      shear & Shear by a factor chosen uniformly from $[-0.2, 0.2]$. A shearing factor of 1 means that a horizontal line turns into a line at 45 degrees. \\
   %     \hline
        
    %    Gaussian noise & Random noise drawn from normal distribution $\mathcal N (0, \sigma)$ with the standard deviation $\sigma$ uniformly chosen from $[0, 0.1]$ is added to the point cloud. \\
     %   \hline 
        
      %  outliers & A percentage, chosen uniformly from $[0, 0.1]$, of point cloud points are replaced with points sampled from a uniform distribution within the range of the point cloud. %\\
    %\end{tabular}
    %\caption{Explanations of the data transformations used on the datasets}
    %\label{tab:transformations}
%\end{table}

% \begin{table}[ht]
% \begin{center}
% \begin{tabular}{ l l l c c } 
% \textbf{Label} & \textbf{Complex type} & \textbf{Subtype} & \textbf{Neighbourhood size} & \textbf{Axes ratios} \\
% \hline
% PHE [1] & Ellipsoid & Rips & 17 & [2 1] \\ 
% PHE [2] & Ellipsoid & Rips & 25 & PCA \\ 
% PHE [3] & Ball & Rips & \textcolor{gray!50}{N/A} & \textcolor{gray!50}{N/A}  \\ 
% \end{tabular}
% \end{center}
% \caption{Summary of parameters used in the classification experiments shown in Figure~\ref{fig:point_cloud_classification}. ``Complex type'' specifies whether balls or ellipsoids were used in defining the simplex tree. ``Subtype'' indicates how the simplex tree was built. ``Neighbourhood size'' is the number of points used in PCA for calculating the orientations of ellipsoids. ``Axes ratios'' are the ellipsoid axes ratios. There, a vector indicates a fixed ratio for all ellipsoids, whereas ``PCA'' axes ratio means that the scaled PCA eigenvalues were used for the axes lengths.}
% \label{table:params}
% \end{table}

The results shown in Figure~\ref{fig:point_cloud_classification} represent average accuracies over 24 runs of the classification pipeline on the same dataset in its original state, as well as after the transformations have been applied to it. The ratio between the training data and the test data remains fixed, but the test and the training data change between the runs. 
% The different PHE pipelines are explained in Table~\ref{table:params}.
The classification based on ellipsoid data performs best in all cases, except when outliers are introduced to the point cloud.

% \anatodo{Get results for smaller point clouds ($\leq 100$ pts) for the difference between ellipsoid and ball to be more apparent. maybe can find a parameter set where pca performs good in a smaller point cloud?}

% We also perform classification based on Rips barcodes on the reduced datasets containing 100 points and the obtained accuracies are shown in Figure~\ref{fig:results_rips}.
% As expected, Rips barcodes lead to worse results.
% \begin{figure}[ht]
%     \centering
%     \includegraphics[width=0.8\textwidth]{figures/rips_accs_trnsfs_averages_id=0005_20240605_145139.png}
%     \caption{Classification accuracies.}
%     \label{fig:results_rips}
% \end{figure}

%\subsubsection{Reproducibility and computer infrastructure}

The code used in this subsection is available at \url{https://github.com/a-zeg/ellipsoids}.
Computations were performed on an Apple M4 Pro CPU.

\subsection{Pentagons}

As the next example, consider a dataset of 14074 points from the configuration space of equilateral planar pentagons, viewed as living in $\RR^6$. 
More precisely, the dataset consists of a sample of 14074 points from
\[
M=\{(x_1, x_2, x_3)\in\RR^6\,|\, \lVert x_i-x_{i+1}\rVert=1, i=1, 2, 3, 4, 5\},
\]
where $x_4$ and $x_5$ are fixed vectors in $\RR^2$ and where we regard $x_6$ as $x_1$.
The dataset was created by Clayton Shonkwiler and provided to us by Henry Adams.

    \begin{figure}[ht]
    \centering
    \includegraphics[width=1\linewidth]{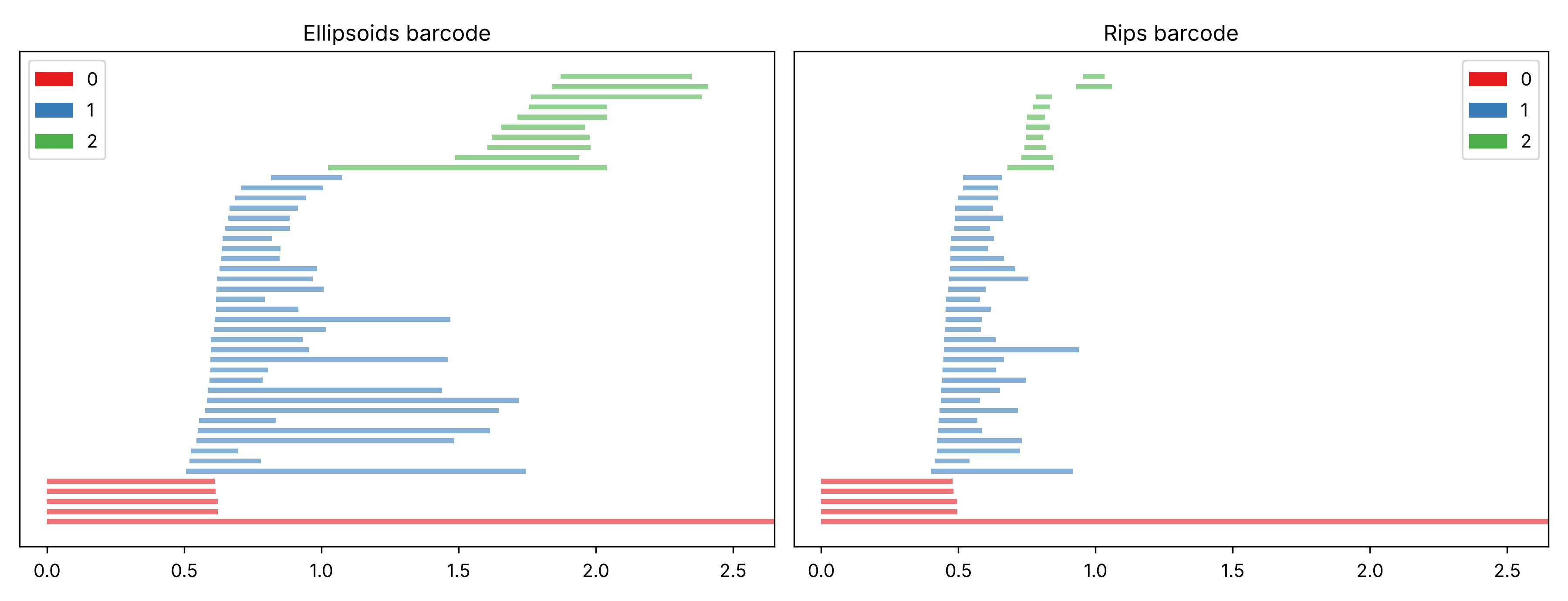}
    \caption{%
        In some cases, even when the pentagons dataset is downsampled to only 130 points, we can see that the ellipsoid barcode captures the correct Betti numbers.
    }
    \label{fig:pentagons_barcodes}
    \end{figure}

It was established in \textcite{Havel} that $M$ is a compact, connected and orientable, two-dimensional manifold of genus 4. We tested if we can recover its homology groups with persistent homology via ellipsoid and Rips complexes %The results are displayed in %Figure~\ref{fig:pentagons_barcodes} %and~\ref{fig:pentagons_barcodes1}. 
and indeed, ellipsoids can detect the `correct homology' with a subsample consisting of as few as 130 points (see Figure~\ref{fig:pentagons_barcodes}).

 %  \textcolor{blue}{Yes :) Below are some pictures for 300 points (where both ellipsoids and Rips work, but ellipsoids works way better) and for 100 points. Depending on the seed used in the downsampling and the neighbourhood size, the ellipsoids barcodes might not always work as well as in the picture below. }

  %  \begin{figure}[ht]
   % \centering
   % \includegraphics[width=1\linewidth]{figures/ellipsoids_data_type=pentagons_n_pts=300_nbhd_size=6_axes_ratios=[3 3 1]__20240816_155703}
    %\caption{%
     %   Ellipsoids and Rips barcode for the pentagons dataset downsampled to 300 points. 
    %}
    %\label{fig:pentagons_barcodes1}
    %\end{figure}

\subsection{Cyclo-octane}
The last example for which we compare the ellipsoid and Rips barcodes is for the conformation space of the cyclo-octane dataset. The cyclo-octane dataset was introduced in \textcite{topcyclooct} and consists of 6040 points in 24 dimensions. It is publicly available as part of the JavaPlex~\cite{Javaplex} software package. 

\begin{figure}[ht]
 \begin{center}
 \begin{minipage}{0.45\linewidth}
\input{cyclooctane-conformations}
   \end{minipage}
   \hspace{1cm}
  \begin{minipage}{0.45\linewidth}
\includegraphics[scale=0.5]{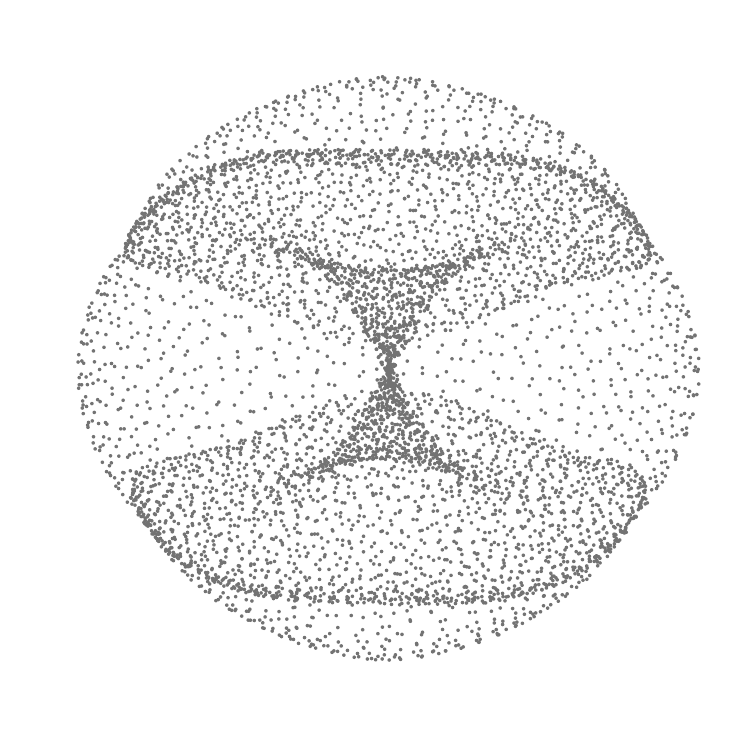}
 \end{minipage}
\caption{Left: Examples of cyclo-octane conformations. Right: The conformation space of cyclo-octane is the union of a sphere with a Klein bottle, glued together along two circles of singularities.}
\label{fig:cyclooctaneconformations}
\end{center}
\end{figure}

 A single molecule of the cyclo-octane consists of eight carbon atoms arranged in a ring, with each carbon atom being bound to two other carbon atoms and two hydrogen atoms. The location of the hydrogen atoms is determined by that of the carbon atoms due to energy minimization. Hence, the conformation space of cyclo-octane consists of all possible spatial arrangements, up to rotation and translation, of the ring of carbon atoms (see the left image in Figure~\ref{fig:cyclooctaneconformations}). Each conformation may therefore be represented by a point in $\RR^{24}$, where we have three spatial coordinates for each of the eight carbon atoms. 
Brown et al.~\cite{Cyclooctane08} and Martin et al.~\cite{topcyclooct} show that the conformation space of cyclo-octane is the union of a sphere with a Klein bottle, glued together along two circles of singularities (see the right image in Figure~\ref{fig:cyclooctaneconformations}).

\begin{figure}[ht]
\begin{center}
\includegraphics[width = 0.45\textwidth]{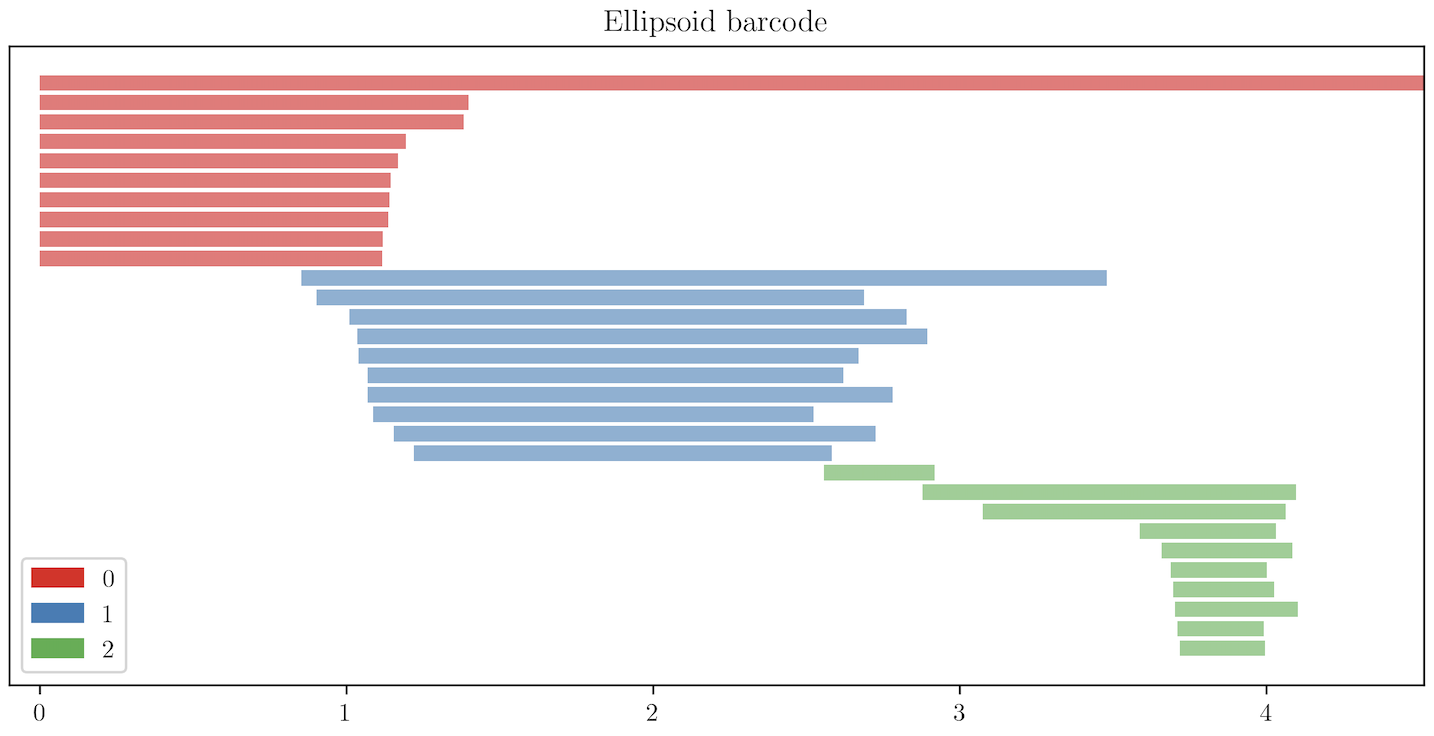}
\includegraphics[width = 0.45\textwidth]{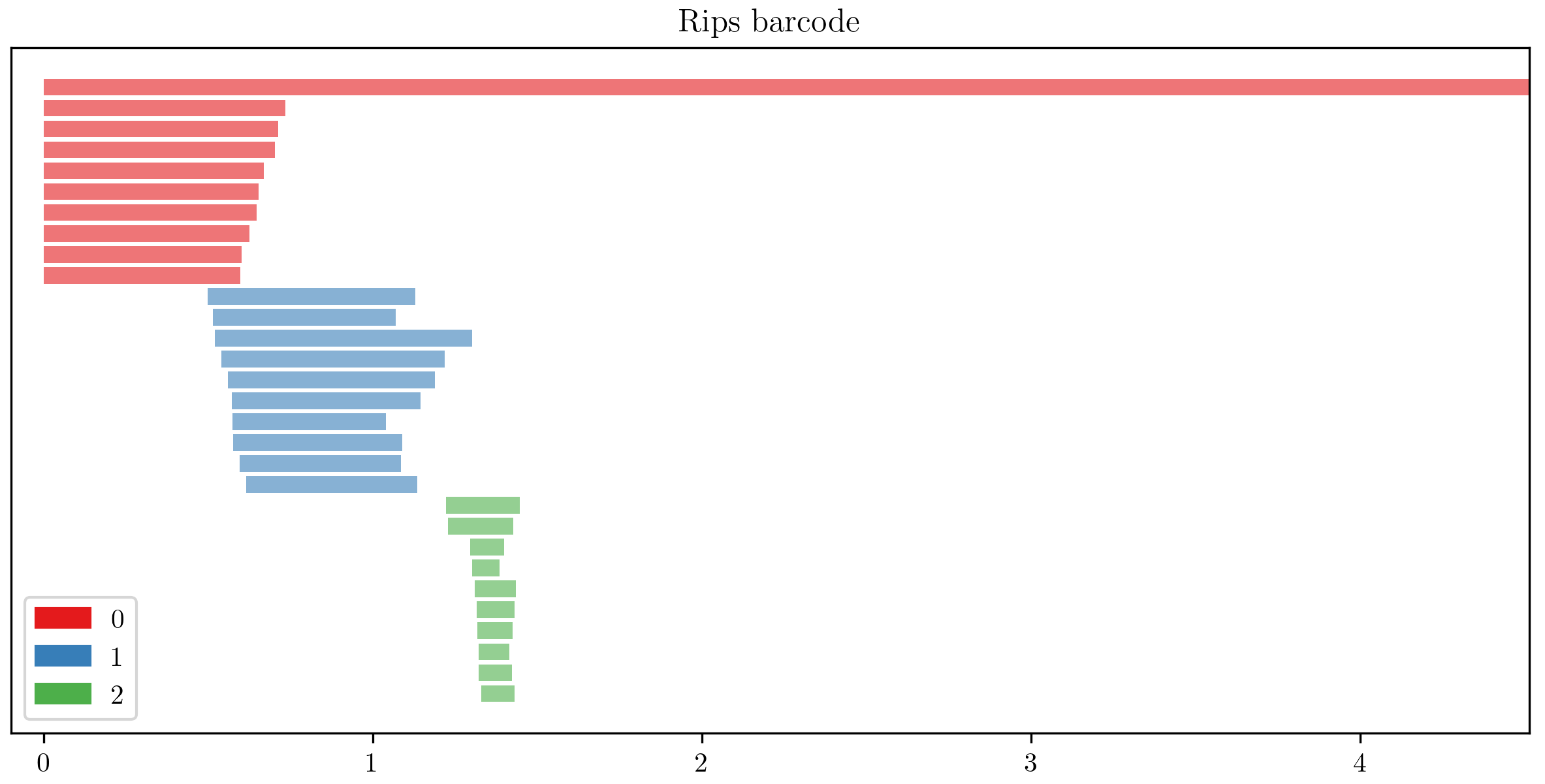}
\caption{Barcodes for a subsample of 500 points from the cyclo-octane dataset. The right plot shows the barcodes for the Vietoris-Rips complex. The left picture shows barcodes for the ellipsoid complex. Only the 10 longest bars per degree are shown.
% the original picture was:
% ellipsoids_dataType=cyclooctane_nPts=500_nbhdSize=26axesRatios=[3, 1]_seed=1_20231129_103034-barcodesDim=0-10_1-10_2-10_20231129_110749
}
\label{fig:cyclooctanebarcodes}
\end{center}
\end{figure}

 The cyclo-octane dataset has been used many times as an example to show that we can recover the homology groups of the conformation space using persistent homology~\cite{zomorodian-textbook, Javaplex}. We confirmed this result using ellipsoid complexes. The results for a 500-point subsample are displayed in Figure~\ref{fig:cyclooctanebarcodes}. The barcodes from the usual Vietoris–Rips complex do not capture the correct homology groups, whereas the ellipsoid barcodes do. In particular, where 2-dimensional Rips barcode only shows noise, the ellipsoid barcode has two prominent bars.

%%% Local Variables:
%%% mode: LaTeX
%%% TeX-master: "../ellipsoid-main"
%%% End:

%% file: cyclooctane-conformations.tex
% Three conformations of cyclo-octane, drawn with TikZ (no text).
% Top left: chair-boat; top right: chair-chair; bottom: boat-boat.
% Carbon positions are taken from the cyclo-octane dataset (JavaPlex,
% pointsCycloOctane.mat); hydrogen atoms are added.
%
% Usage: \input{cyclooctane-conformations}     (needs \usepackage{tikz})
% The picture is \cycloWidth wide (default 0.48\textwidth); its height follows
% from the drawing (0.672 x width). To change it, before \input write e.g.
%   \def\cycloWidth{0.4\textwidth}
\providecommand*\cycloWidth{0.48\textwidth}%
\ifdefined\cycloU\else\newlength\cycloU\fi
\begingroup
\pgfmathsetlength\cycloU{(\cycloWidth)/14.280}%
\begin{tikzpicture}[x=\cycloU, y=\cycloU,
    cyclo C/.style={ball color=black!70},
    cyclo H/.style={ball color=white, draw=black!45, line width=0.015\cycloU},
    cyclo CC/.style={draw=black!55, line width=0.09\cycloU, line cap=round},
    cyclo CH/.style={draw=black!40, line width=0.04\cycloU, line cap=round}]
  \useasboundingbox (-7.140,-4.800) rectangle (7.140,4.800);
  % chair-boat (top left)
  \begin{scope}[shift={(-3.870,2.600)}]
    \draw[cyclo CH] (-0.550,-0.077) -- (-0.883,-0.038);
    \shadedraw[cyclo H] (-0.883,-0.038) circle (0.2);
    \draw[cyclo CH] (0.414,1.073) -- (0.678,1.516);
    \shadedraw[cyclo H] (0.678,1.516) circle (0.2);
    \draw[cyclo CC] (-0.550,-0.077) -- (-1.162,-0.085);
    \draw[cyclo CC] (-0.550,-0.077) -- (-0.068,0.498);
    \draw[cyclo CH] (-0.550,-0.077) -- (-0.012,-1.016);
    \shade[cyclo C] (-0.550,-0.077) circle (0.36);
    \shadedraw[cyclo H] (-0.012,-1.016) circle (0.2);
    \draw[cyclo CC] (0.414,1.073) -- (-0.068,0.498);
    \draw[cyclo CC] (0.414,1.073) -- (1.053,0.883);
    \draw[cyclo CH] (0.414,1.073) -- (-0.101,1.842);
    \shade[cyclo C] (0.414,1.073) circle (0.36);
    \draw[cyclo CH] (1.693,0.694) -- (2.178,-0.115);
    \shadedraw[cyclo H] (2.178,-0.115) circle (0.2);
    \draw[cyclo CH] (-1.774,-0.092) -- (-2.485,0.648);
    \shadedraw[cyclo H] (-2.485,0.648) circle (0.2);
    \draw[cyclo CH] (-1.774,-0.092) -- (-2.256,-1.066);
    \shadedraw[cyclo H] (-2.256,-1.066) circle (0.2);
    \shadedraw[cyclo H] (-0.101,1.842) circle (0.2);
    \draw[cyclo CC] (-1.774,-0.092) -- (-1.633,0.047);
    \draw[cyclo CC] (-1.774,-0.092) -- (-1.162,-0.085);
    \shade[cyclo C] (-1.774,-0.092) circle (0.36);
    \draw[cyclo CH] (1.693,0.694) -- (2.372,1.546);
    \shadedraw[cyclo H] (2.372,1.546) circle (0.2);
    \draw[cyclo CC] (1.693,0.694) -- (1.053,0.883);
    \draw[cyclo CC] (1.693,0.694) -- (1.598,0.478);
    \shade[cyclo C] (1.693,0.694) circle (0.36);
    \draw[cyclo CH] (0.745,-1.042) -- (0.740,-1.580);
    \shadedraw[cyclo H] (0.740,-1.580) circle (0.2);
    \draw[cyclo CC] (-1.493,0.186) -- (-1.091,-0.350);
    \draw[cyclo CC] (-1.493,0.186) -- (-1.633,0.047);
    \draw[cyclo CH] (-1.493,0.186) -- (-0.958,1.133);
    \draw[cyclo CH] (-1.493,0.186) -- (-2.443,0.313);
    \shade[cyclo C] (-1.493,0.186) circle (0.36);
    \draw[cyclo CC] (1.504,0.262) -- (1.598,0.478);
    \draw[cyclo CC] (1.504,0.262) -- (1.125,-0.390);
    \draw[cyclo CH] (1.504,0.262) -- (0.976,1.052);
    \draw[cyclo CH] (1.504,0.262) -- (2.485,0.164);
    \shade[cyclo C] (1.504,0.262) circle (0.36);
    \shadedraw[cyclo H] (-0.958,1.133) circle (0.2);
    \draw[cyclo CC] (0.745,-1.042) -- (1.125,-0.390);
    \draw[cyclo CC] (0.745,-1.042) -- (0.028,-0.964);
    \draw[cyclo CH] (0.745,-1.042) -- (1.287,-1.660);
    \shade[cyclo C] (0.745,-1.042) circle (0.36);
    \shadedraw[cyclo H] (2.485,0.164) circle (0.2);
    \shadedraw[cyclo H] (-2.443,0.313) circle (0.2);
    \shadedraw[cyclo H] (0.976,1.052) circle (0.2);
    \draw[cyclo CH] (-0.689,-0.887) -- (-1.200,-1.842);
    \shadedraw[cyclo H] (-1.200,-1.842) circle (0.2);
    \draw[cyclo CC] (-0.689,-0.887) -- (-1.091,-0.350);
    \draw[cyclo CH] (-0.689,-0.887) -- (-0.677,-0.661);
    \draw[cyclo CC] (-0.689,-0.887) -- (0.028,-0.964);
    \shade[cyclo C] (-0.689,-0.887) circle (0.36);
    \shadedraw[cyclo H] (1.287,-1.660) circle (0.2);
    \shadedraw[cyclo H] (-0.677,-0.661) circle (0.2);
  \end{scope}
  % chair-chair (top right)
  \begin{scope}[shift={(3.870,2.600)}]
    \draw[cyclo CH] (0.044,0.382) -- (0.061,0.801);
    \shadedraw[cyclo H] (0.061,0.801) circle (0.2);
    \draw[cyclo CH] (0.044,0.382) -- (0.028,-0.703);
    \shadedraw[cyclo H] (0.028,-0.703) circle (0.2);
    \draw[cyclo CH] (-1.238,0.821) -- (-1.982,1.049);
    \shadedraw[cyclo H] (-1.982,1.049) circle (0.2);
    \draw[cyclo CH] (1.325,0.780) -- (2.105,0.937);
    \shadedraw[cyclo H] (2.105,0.937) circle (0.2);
    \draw[cyclo CC] (0.044,0.382) -- (-0.597,0.601);
    \draw[cyclo CC] (0.044,0.382) -- (0.685,0.581);
    \shade[cyclo C] (0.044,0.382) circle (0.36);
    \draw[cyclo CC] (-1.238,0.821) -- (-1.533,0.315);
    \draw[cyclo CC] (-1.238,0.821) -- (-0.597,0.601);
    \draw[cyclo CH] (-1.238,0.821) -- (-1.043,1.747);
    \shade[cyclo C] (-1.238,0.821) circle (0.36);
    \draw[cyclo CC] (1.325,0.780) -- (0.685,0.581);
    \draw[cyclo CC] (1.325,0.780) -- (1.577,0.279);
    \draw[cyclo CH] (1.325,0.780) -- (1.164,1.736);
    \shade[cyclo C] (1.325,0.780) circle (0.36);
    \shadedraw[cyclo H] (1.164,1.736) circle (0.2);
    \shadedraw[cyclo H] (-1.043,1.747) circle (0.2);
    \draw[cyclo CH] (-1.828,-0.192) -- (-1.611,-1.194);
    \shadedraw[cyclo H] (-1.611,-1.194) circle (0.2);
    \draw[cyclo CH] (1.828,-0.223) -- (1.611,-1.228);
    \shadedraw[cyclo H] (1.611,-1.228) circle (0.2);
    \draw[cyclo CC] (-1.828,-0.192) -- (-1.577,-0.138);
    \draw[cyclo CC] (-1.828,-0.192) -- (-1.533,0.315);
    \draw[cyclo CH] (-1.828,-0.192) -- (-2.912,-0.081);
    \shade[cyclo C] (-1.828,-0.192) circle (0.36);
    \shadedraw[cyclo H] (-2.912,-0.081) circle (0.2);
    \draw[cyclo CC] (1.828,-0.223) -- (1.577,0.279);
    \draw[cyclo CC] (1.828,-0.223) -- (1.533,-0.146);
    \draw[cyclo CH] (1.828,-0.223) -- (2.912,-0.136);
    \shade[cyclo C] (1.828,-0.223) circle (0.36);
    \shadedraw[cyclo H] (2.912,-0.136) circle (0.2);
    \draw[cyclo CH] (-0.044,-0.847) -- (-0.028,-1.747);
    \shadedraw[cyclo H] (-0.028,-1.747) circle (0.2);
    \draw[cyclo CC] (-1.325,-0.083) -- (-0.685,-0.465);
    \draw[cyclo CC] (-1.325,-0.083) -- (-1.577,-0.138);
    \draw[cyclo CH] (-1.325,-0.083) -- (-1.164,0.969);
    \draw[cyclo CH] (-1.325,-0.083) -- (-2.105,-0.442);
    \shade[cyclo C] (-1.325,-0.083) circle (0.36);
    \draw[cyclo CC] (1.238,-0.068) -- (1.533,-0.146);
    \draw[cyclo CC] (1.238,-0.068) -- (0.597,-0.458);
    \draw[cyclo CH] (1.238,-0.068) -- (1.043,0.989);
    \draw[cyclo CH] (1.238,-0.068) -- (1.982,-0.384);
    \shade[cyclo C] (1.238,-0.068) circle (0.36);
    \shadedraw[cyclo H] (1.043,0.989) circle (0.2);
    \shadedraw[cyclo H] (-1.164,0.969) circle (0.2);
    \draw[cyclo CC] (-0.044,-0.847) -- (-0.685,-0.465);
    \draw[cyclo CH] (-0.044,-0.847) -- (-0.061,-1.173);
    \draw[cyclo CC] (-0.044,-0.847) -- (0.597,-0.458);
    \shade[cyclo C] (-0.044,-0.847) circle (0.36);
    \shadedraw[cyclo H] (-2.105,-0.442) circle (0.2);
    \shadedraw[cyclo H] (1.982,-0.384) circle (0.2);
    \shadedraw[cyclo H] (-0.061,-1.173) circle (0.2);
  \end{scope}
  % boat-boat (bottom)
  \begin{scope}[shift={(0.000,-2.600)}]
    \draw[cyclo CH] (-1.278,0.166) -- (-1.345,-0.216);
    \shadedraw[cyclo H] (-1.345,-0.216) circle (0.2);
    \draw[cyclo CH] (1.286,0.154) -- (1.349,-0.237);
    \shadedraw[cyclo H] (1.349,-0.237) circle (0.2);
    \draw[cyclo CH] (0.008,0.966) -- (0.013,1.748);
    \shadedraw[cyclo H] (0.013,1.748) circle (0.2);
    \draw[cyclo CC] (-1.278,0.166) -- (-0.635,0.566);
    \draw[cyclo CH] (-1.278,0.166) -- (-2.123,0.838);
    \draw[cyclo CC] (-1.278,0.166) -- (-1.345,-0.417);
    \shade[cyclo C] (-1.278,0.166) circle (0.36);
    \draw[cyclo CC] (1.286,0.154) -- (0.647,0.560);
    \draw[cyclo CC] (1.286,0.154) -- (1.349,-0.426);
    \draw[cyclo CH] (1.286,0.154) -- (2.138,0.820);
    \shade[cyclo C] (1.286,0.154) circle (0.36);
    \draw[cyclo CC] (0.008,0.966) -- (-0.635,0.566);
    \draw[cyclo CC] (0.008,0.966) -- (0.647,0.560);
    \draw[cyclo CH] (0.008,0.966) -- (0.008,1.463);
    \shade[cyclo C] (0.008,0.966) circle (0.36);
    \shadedraw[cyclo H] (2.138,0.820) circle (0.2);
    \shadedraw[cyclo H] (-2.123,0.838) circle (0.2);
    \draw[cyclo CH] (-1.413,-1.001) -- (-0.652,-1.743);
    \shadedraw[cyclo H] (-0.652,-1.743) circle (0.2);
    \draw[cyclo CH] (1.413,-1.005) -- (0.652,-1.748);
    \shadedraw[cyclo H] (0.652,-1.748) circle (0.2);
    \draw[cyclo CH] (-1.413,-1.001) -- (-2.382,-1.474);
    \shadedraw[cyclo H] (-2.382,-1.474) circle (0.2);
    \draw[cyclo CH] (1.413,-1.005) -- (2.382,-1.482);
    \shadedraw[cyclo H] (2.382,-1.482) circle (0.2);
    \draw[cyclo CC] (-1.413,-1.001) -- (-1.349,-0.818);
    \draw[cyclo CC] (-1.413,-1.001) -- (-1.345,-0.417);
    \shade[cyclo C] (-1.413,-1.001) circle (0.36);
    \draw[cyclo CC] (1.413,-1.005) -- (1.349,-0.426);
    \draw[cyclo CC] (1.413,-1.005) -- (1.345,-0.817);
    \shade[cyclo C] (1.413,-1.005) circle (0.36);
    \shadedraw[cyclo H] (0.008,1.463) circle (0.2);
    \draw[cyclo CH] (-0.008,0.097) -- (-0.008,1.070);
    \shadedraw[cyclo H] (-0.008,1.070) circle (0.2);
    \draw[cyclo CC] (-1.286,-0.636) -- (-0.647,-0.269);
    \draw[cyclo CC] (-1.286,-0.636) -- (-1.349,-0.818);
    \draw[cyclo CH] (-1.286,-0.636) -- (-2.138,-0.015);
    \draw[cyclo CH] (-1.286,-0.636) -- (-1.349,-1.549);
    \shade[cyclo C] (-1.286,-0.636) circle (0.36);
    \draw[cyclo CC] (1.278,-0.629) -- (1.345,-0.817);
    \draw[cyclo CC] (1.278,-0.629) -- (0.635,-0.266);
    \draw[cyclo CH] (1.278,-0.629) -- (2.123,0.000);
    \draw[cyclo CH] (1.278,-0.629) -- (1.345,-1.536);
    \shade[cyclo C] (1.278,-0.629) circle (0.36);
    \shadedraw[cyclo H] (-2.138,-0.015) circle (0.2);
    \shadedraw[cyclo H] (2.123,0.000) circle (0.2);
    \draw[cyclo CC] (-0.008,0.097) -- (0.635,-0.266);
    \draw[cyclo CC] (-0.008,0.097) -- (-0.647,-0.269);
    \draw[cyclo CH] (-0.008,0.097) -- (-0.013,0.284);
    \shade[cyclo C] (-0.008,0.097) circle (0.36);
    \shadedraw[cyclo H] (-1.349,-1.549) circle (0.2);
    \shadedraw[cyclo H] (1.345,-1.536) circle (0.2);
    \shadedraw[cyclo H] (-0.013,0.284) circle (0.2);
  \end{scope}
\end{tikzpicture}%
\endgroup

%% file: sections/appendix.tex
Let $k,n,N$ be positive integers such that $N\ge k+2$ and $k\ge n$. Denote by $\mathcal{G}_{n,k}^N$ the set of all $k$-generic $N$-element subsets of $\RR^n$. We shall call an element of $\mathcal{G}_{n, k}^N$ a \define{$k$-generic point cloud with $N$ points in $\R^n$}. One may view $\mathcal{G}_{n, k}^N$ as a subset of $(\R^n)^N$ by assigning an arbitrary order to the elements in any subset contained in $\mathcal{G}_{n, k}^N$.
%Furthermore, for $N\in \mathbb{N}$ we define
%\[
%\mathcal{G}^N_{n,k} \coloneqq \underbrace{\mathcal{G}_{n, k}\times \cdots\times \mathcal{G}_{n, k}}_{N\text{ times}}.
%\]
%Intuitively, one should think of an element of $\mathcal{G}_{n, k}^N$ as the partition of a point cloud into $k$-neighborhoods.
%Note that $\mathcal{G}^N_{n,k}$ is a subset of $(\mathbb{R}^n)^N$. 
\begin{proposition}\label{prop:open_and_dense}
    The set of $k$-generic point clouds with $N$ points in $\mathbb{R}^n$ is open and dense in $(\mathbb{R}^n)^N$ endowed with the product Euclidean topology. In particular, any point cloud can be perturbed by an arbitrarily small amount to a $k$-generic point cloud.
\end{proposition}
\begin{proof}
First we prove that $\mathcal{G}^N_{n,k}$ is open in  $(\mathbb{R}^n)^N$. Let $X\in \mathcal{G}^N_{n,k}$ be arbitrary and write $X=(p_1,\dots,p_N)$ with $p_i\neq p_j$ for $i\neq j$.
Define
\begin{align*}
    \rho \coloneqq \min_{i\neq j}\|p_i-p_j\|>0 \quad\text{ and }\quad \varepsilon_0 \coloneqq \frac{\rho}{4}.
\end{align*}

For each $i\in\{1,\dots,N\}$, order the numbers $\{\|p_j-p_i\| : j\neq i\}$ as
\[
0<r_{i,1}<\dots<r_{i,k+1}\le\dots\le r_{i,N-1}.
\]
The first $k+1$ inequalities are strict by the separation condition.

%For each $i\in\{1,\dots,N\}$, the separation condition implies that the numbers
%$\{\|p_j-p_i\|: j\neq i\}$ are pairwise distinct. Order them as
%\[
%0<r_{i,1}<\cdots<r_{i,N-1},
%\]
Set
\begin{align*}
    \Delta_i \coloneqq \min\Bigl(\min_{1\le a<b\le k}(r_{i,b}-r_{i,a}),\ r_{i,k+1}-r_{i,k}\Bigr)>0 \quad\text{ and }\quad \Delta \coloneqq \min_{1\le i\le N}\Delta_i>0
\end{align*}
Furthermore set $\varepsilon_1 \coloneqq \frac{\Delta}{8}$.
Moreover, since $X$ is $k$-generic (with $m$-gap), we have
\[
\eta_X \coloneqq \min_{1\le i\le N}\Bigl(\sigma_m(\mathbf N_k(p_i))-\sigma_{m+1}(\mathbf N_k(p_i))\Bigr)>0.
\]
We also set $\varepsilon_2\coloneqq \frac{\eta_X}{8\sqrt{k}}$ and $\varepsilon \coloneqq \min(\varepsilon_0,\varepsilon_1,\varepsilon_2).$
\newline
\newline
Let $Y=(q_1,\dots,q_N)$ satisfy $\max_{1\le i\le N}\|q_i-p_i\|<\varepsilon$.
Then for $i\neq j$ by the reverse triangle inequality
\begin{align*}
    \|q_i-q_j\|
    &\ge \|p_i-p_j\|-\|q_i-p_i\|-\|q_j-p_j\|
    \ge \rho-2\varepsilon
    \ge \frac{\rho}{2}>0,
\end{align*}
hence the points of $Y$ are pairwise distinct.
\newline
\newline
Fix $i\neq j $. By the reverse triangle inequality,
\begin{align} \label{eq:reverseTriangle}
\bigl|\|q_j-q_i\|-\|p_j-p_i\|\bigr|
\le \|q_j-p_j\|+\|q_i-p_i\|
<2\varepsilon.
\end{align}

Hence, for any $i \neq j, \ell$ using \eqref{eq:reverseTriangle} twice implies:
\begin{align*}
   \|q_j-q_i\|-\|p_j-p_i\| + \|p_\ell-p_i\| - \|q_\ell-q_i\|
<4\varepsilon.
\end{align*}
Assuming that $\|p_\ell-p_i\|-\|p_j-p_i\|\ge 4\varepsilon$ implies:
\[
\|p_j-p_i\|<\|p_\ell-p_i\|
\ \Longrightarrow\
\|q_j-q_i\|<\|q_\ell-q_i\|.
\]
Since $4\varepsilon\le \Delta_i$, the strict ordering of the first $k$ distances from $p_i$
and the strict gap $r_{i,k}<r_{i,k+1}$ persist for $q_i$. Therefore, for each $i$, the set
$N_k(q_i)$ has exactly $k$ points and is strictly ordered by distance; in particular,
the neighbor indices of $q_i$ coincide with those of $p_i$.
\newline
\newline
Now fix $i$ and write $N_k(p_i)=\{p_{j_1},\dots,p_{j_k}\}$ (same indices for $q_i$).
Let $\bar p_i=\frac1k\sum_{\ell=1}^k p_{j_\ell}$ and $\bar q_i=\frac1k\sum_{\ell=1}^k q_{j_\ell}$.
Then $\|\bar q_i-\bar p_i\|\le \varepsilon$, and for each column $\ell$,
\[
\|(q_{j_\ell}-\bar q_i)-(p_{j_\ell}-\bar p_i)\|
\le \|q_{j_\ell}-p_{j_\ell}\|+\|\bar q_i-\bar p_i\|
<2\varepsilon.
\]
Hence by Proposition~\ref{prop:norm_comp} we have:
\[
\|\mathbf N_k(q_i)-\mathbf N_k(p_i)\| < 2\sqrt{k}\,\varepsilon.
\]
Since the spectral norm is bounded by the Frobenius norm, Corollary 7.3.5~\cite{Horn_Johnson_1985} implies that for all $a$
\[
\bigl|\sigma_a(\mathbf N_k(q_i))-\sigma_a(\mathbf N_k(p_i))\bigr|
< 2\sqrt{k}\,\varepsilon.
\]
Therefore,
\begin{align*}
&\sigma_m(\mathbf N_k(q_i))-\sigma_{m+1}(\mathbf N_k(q_i))
\\&\ge \bigl(\sigma_m(\mathbf N_k(p_i))-\sigma_{m+1}(\mathbf N_k(p_i))\bigr)
      -\bigl|\sigma_m(\mathbf N_k(q_i))-\sigma_m(\mathbf N_k(p_i))\bigr|
      -\bigl|\sigma_{m+1}(\mathbf N_k(q_i))-\sigma_{m+1}(\mathbf N_k(p_i))\bigr| \\
&> \eta_X - 4\sqrt{k}\,\varepsilon
\ \ge\ \frac{\eta_X}{2},
\end{align*}
since $|\sigma_a(\mathbf N_k(q_i))-\sigma_a(\mathbf N_k(p_i))|<2\sqrt{k}\varepsilon$ for all $a$ and
$\varepsilon\le \varepsilon_2=\eta_X/(8\sqrt{k})$.
Thus $Y$ satisfies the $m$-spectral gap condition (with uniform constant $\eta_X/2$), and consequently
$Y\in \mathcal{G}^N_{n,k}$. This shows that $\mathcal{G}^N_{n,k}$ is open.

Now we argue that $\mathcal{G}^N_{n,k}$ is dense in $(\RR^n)^N$. For $1\le i<j\le N$ we define
\[
F_{i,j}(X)\coloneqq \|p_i-p_j\|^2,
\]
and for $1\le i\le N$ and $j,\ell\in\{1,\dots,N\}\setminus\{i\}$, $j<\ell$, set
\[
G_{i;j,\ell}(X)\coloneqq \|p_j-p_i\|^2-\|p_\ell-p_i\|^2 .
\]
Furthermore, for $1\le i\le N$ and any subset $I=\{j_1,\dots,j_k\}\subset\{1,\dots,N\}\setminus\{i\}$ of cardinality $k$ we define
\[
\bar p_I \coloneqq \frac1k\sum_{t=1}^k p_{j_t},\qquad
\mathbf N_I(p_i)\coloneqq \begin{pmatrix} \,|& &|\\ p_{j_1}-\bar p_I&\cdots&p_{j_k}-\bar p_I\\ |& &| \end{pmatrix}\in\mathbb R^{d\times k},
\]
%and its Gram matrix $H_{i,I}(X)\coloneqq \mathbf N_I(p_i)^\top \mathbf N_I(p_i)\in\mathbb R^{k\times k}$.
and its Gram matrix $H_{i,I}(X):= \mathbf N_I(p_i) \mathbf N_I(p_i)^\top\in\RR^{n\times n}$. Let
\[
\chi_{i,I}(\lambda;X):=\det(\lambda I_n-H_{i,I}(X))
\]
%Let
%\[
%\chi_{i,I}(\lambda;X)\coloneqq \det(\lambda I_k - H_{i,I}(X))
%\]
and define
\[
D_{i,I}(X)\coloneqq \mathrm{Disc}\bigl(\chi_{i,I}(\cdot;X)\bigr),
\]
the discriminant of $\chi_{i,I}$. Hence, $D_{i,I}(X)=0$ if and only if $H_{i,I}(X)$ has a repeated eigenvalue, equivalently, if $\mathbf N_I(p_i)$ has a repeated
singular value.

Define
\[
\mathcal U \coloneqq (\mathbb R^n)^N \setminus
\Biggl(
\bigcup_{i<j}\{F_{i,j}=0\}\ \cup\
\bigcup_{i}\ \bigcup_{j<\ell}\{G_{i;j,\ell}=0\}\ \cup\
\bigcup_{i}\ \bigcup_{I}\{D_{i,I}=0\}
\Biggr),
\]
where in the middle union $j,\ell$ range over $\{1,\dots,N\}\setminus\{i\}$ and in the last union $I$ ranges over
all $k$-subsets of $\{1,\dots,N\}\setminus\{i\}$. Note that trivially $\mathcal{U}\subset \mathcal{G}^N_{n,k}$.
\newline
\newline
Each set $\{F_{i,j}=0\}$ and $\{G_{i;j,\ell}=0\}$ is the zero set of a nonzero polynomial.
Moreover, for each $(i,I)$, $D_{i,I}$ is not the zero polynomial, hence $\{D_{i,I}=0\}$ is a proper algebraic subset. Indeed, since $k\ge n$, the points $p_{j_1},\dots,p_{j_k}$ can be chosen so that $\mathbf{N}_I(p_i)$ has $n$ distinct singular values.

Therefore the union above is a finite union of proper algebraic sets, so $\mathcal U$ is open and dense in
$(\mathbb R^n)^N$, thus also $\mathcal{G}^N_{n,k}$.

\end{proof}